\pdfoutput=1
\documentclass[a4paper,11pt]{article}

\usepackage[OT1]{fontenc}
\usepackage{amsfonts}
\usepackage{amsmath}
\input epsf
\usepackage{epsfig}
\usepackage{latexsym}
\usepackage{graphicx}
\usepackage{footmisc}

\usepackage[utf8]{inputenc}

\renewcommand{\thesection}{\arabic{section}}  
\DefineFNsymbols*{asterisks}{*{**}{***}{****}{*****}{******}{*******}}
\setfnsymbol{asterisks}

\newtheorem{thm}{\bf Theorem}[section]
\newtheorem{lem}[thm]{\bf Lemma}
\newtheorem{prop}[thm]{\bf Proposition}

\newtheorem{defi}[thm]{\bf Definition}
\newtheorem{rem}[thm]{\bf Remark}

\makeatletter

\@addtoreset{equation}{section}
\makeatother

\newcommand{\R}{\mathbb{R}}
\newcommand{\CC}{\mathbb{C}}
\newcommand{\Z}{\mathbb{Z}}
\newcommand{\NN}{\mathbb{N}}
\newcommand{\Q}{\mathbb{Q}}

\newcommand{\Frac}[2]{{\displaystyle\frac{ #1}{ #2}}} 
\newcommand{\TFrac}[2]{{\textstyle\frac{ #1}{ #2}}} 

\newcommand{\dindice}[2]{{\stackrel{\scriptstyle #1}{\scriptstyle #2}}}

\newcommand{\supl}{\sup\limits}

\newcommand{\grandO}[1]{\mathop{\hbox{${\cal O}$}}\limits_{#1}}
\newcommand{\Ocal}{{\cal O}}
\newcommand{\calO}{{\cal O}}

\newcommand{\I}{{\rm i}}
\newcommand{\E}{{\rm e}}
\newcommand{\eps}{\varepsilon}

\renewcommand{\ker}{\textrm{Ker}\;}

\newcommand{\hsp}[1][1]{\hspace{#1  ex}}
\newcommand{\lba}{\vspace{0.5ex}\\}

\newcommand{\cqfd}{\hfill$\Box$}
\def\ie{i.e.}

\newcommand{\scal}[2]{\left \langle #1, #2 \right   \rangle}
\newcommand{\Norme}[2]%
  {\left |  #2 \right |_{\hspace{-0.9ex}\raisebox{-0.5ex}{ $\scriptstyle
{#1}$}}}
\newcommand{\NNorme}[2]%
{\left \|  #2 \right \|_{\hspace{-0.9ex}\raisebox{-0.5ex}{ $\scriptstyle
#1  $}}}
\newcommand{\norme}[2]%
  {\left |  #2 \right |_{\hspace{-0.9ex}\raisebox{-0.5ex}{ $\scriptscriptstyle
{#1}$}}}
\newcommand{\nnorme}[2]%
{\left \|  #2 \right \|_{\hspace{-0.9ex}\raisebox{-0.5ex}{ $\scriptscriptstyle
#1  $}}}
\newcommand{\sNorme}[2]%
  { |  #2 |_{\hspace{-0.9ex}\raisebox{-0.5ex}{ $\scriptstyle #1$}}}
\newcommand{\sNNorme}[2]%
{ \|  #2  \|_{\hspace{-0.9ex}\raisebox{-0.5ex}{ $\scriptstyle #1$}}}

\newcommand{\bNorme}[2]%
  { \Bigl |  #2 \Bigr |_{\hspace{-0.9ex}\raisebox{-0.5ex}{ $\scriptstyle #1$}}}
\newcommand{\bNNorme}[2]%
{ \Bigl \|  #2 \Bigr   \|_{\hspace{-0.9ex}\raisebox{-0.5ex}{ $\scriptstyle #1$}}}

\newcommand{\+}{\hspace{-0,5ex}+\hspace{-0,5ex}}
\newcommand{\moins}{\hspace{-0.5ex}-\hspace{-0.5ex}}

\def\H{H}									
\def\N{n}
\def\n{n}									
\def\Q{N_\N}							
\def\Reste{\mathcal{R}_\n} 
\def\No{N_0}							

\def\ql{q_1}
\def\qz{q_2}
\def\pl{p_1}
\def\pz{p_2}					
\def\qp{(\ql,\pl,\qz,\pz)}

\def\ie{i.e.}
\def\An{\mathcal{A}}      
\def\C{\mathcal{C}}				

\def\Om{\Omega} 					
\def\Hbf{\mathbf{H}}
\def\Hlb{\mathbf{H_\lb}} 	
\def\lb{\lambda} 					
\def\V{\mathbf{V}_{\Hlb}} 
\def\om{\omega}
\def\P{P}									
\def\a{a}									
\def\Czero{C_0}						
\def\lzero{\ell_0}				

\def\roo{\rho_0}			  	

\def\Hprime{\underline{H}}            
\def\qlp{\underline{q}_1}
\def\plp{\underline{p}_1}
\def\qzp{\underline{q}_2}
\def\pzp{\underline{p}_2}						
\def\xp{\underline{x}}
\def\Izp{\underline{I}_2}
\def\Nprime{\underline{N}}						
\def\Rprime{\underline{R}}						

\def\entiere{E}						
\def\Gam{\Gamma}					

\def\WIP{W^{u}(\P)}               
\def\WIPa{W^{u}(\P^\a)}

\def\WSP{W^{s}(\P)}               

\newcommand{\Hcal}{{\cal H}}
\newcommand{\Rcal}{{\cal R}}

\newcommand{\xtilde}{\widetilde{x}}
\newcommand{\Hcaltilde}{\widetilde{\Hcal}}
\newcommand{\Ncal}{{\cal N}}

\newcommand{\xbf}{\mathbf{x}}
\newcommand{\ybf}{\mathbf{y}}
\newcommand{\pbfun}{\mathbf{p}_{_1}}
\newcommand{\pbfd}{\mathbf{p}_{_2}}
\newcommand{\qbfun}{\mathbf{q}_{_1}}
\newcommand{\qbfd}{\mathbf{q}_{_2}}

\def\omo{{\om_0}}

\newcommand{\Htildeulb}{\underline{\widetilde{H}}_\lb}
\newcommand{\qtildeun}{\widetilde{q}_1}
\newcommand{\ptildeun}{\widetilde{p}_1}
\newcommand{\qtilded}{\widetilde{q}_2}
\newcommand{\ptilded}{\widetilde{p}_2}

\newcommand{\Ncaltildeunlb}{\underline{\widetilde{\cal N}}_{n,\lb}}
\newcommand{\Mcaltildeunlb}{\underline{\widetilde{\cal M}}_{n,\lb}}
\newcommand{\Rcaltildeunlb}{\underline{\widetilde{\cal R}}_{n,\lb}}

\newcommand{\Ptildeunlb}{\underline{\widetilde{H}}_{n,\lb}}
\newcommand{\uctilde}{\underline{\widetilde{c}}}

\newcommand{\Qtildeunlb}{\underline{\widetilde{Q}}_{n,\lb}}
\newcommand{\epstilde}{\widetilde{\eps}}
\newcommand{\thetatilde}{\widetilde{\theta}}

\newcommand{\Htilde}{\widetilde{H}}

\newcommand{\ttilde}{\widetilde{t}}
\newcommand{\Rtilde}{\widetilde{R}}
\newcommand{\Qtilde}{\widetilde{Q}}
\newcommand{\ctilde}{\widetilde{c}}

\def\flot{\phi}     		
\def\flotp{\widetilde{\phi}}				
\def\Hp{\widetilde{H}}						
\def\rot{\texttt{R}}							
\def\Rot{R}								
\def\rotep{\texttt{R}_{\frac{\om(\eps)}{2\eps^2}t}}     
\def\Rotep{R_{\frac{\om(\eps)}{2\eps^2}t}}   

\def\ko{k_0}               
\def\lo{\ell_0} 						

\def\Rotepm{R_{-\frac{\om(\eps)}{2\eps^2}t}}   

\def\Ret{\Psi}         
\def\Tret{T}						

\def\Reta{\Psi^\a}      
\def\nup{\underline{\nuu}}
\def\cl{c_1}
\def\cz{c_2}
\def\co{c_0}
\def\ucl{\underline{c}_1} %
\def\ucz{\underline{c}_2}  

\def\r{r}
\def\teta{\theta}				

\def\Teta{\Theta}				
\def\Iz{I_2}

\newcommand{\NormeC}[2]{\left| #2 \right|_{\C^{#1}}}

\def\Retah{\tilde{\Psi}^\a}   

\def\nuup{\bar{\eps}}         
\def\q{q}
\def\ro{\rho}						
\def\Retap{\hat{\Psi}^\a}  

\def\Sigl{\Sigma_1}				
\def\Sigz{\Sigma_2}				
\def\dlo{\delta_0}			
\def\dll{\delta}			
\def\dlz{\delta}			

\def\B{\mathcal{B}}    

\def\CS{W^{cs}(0)}     
\def\CSep{W_{\ep}^{cs}(0)}
\def\SPa{W^s(\P^\a)}   
\def\gcsep{p_{1,\ep}^{cs}} 
\def\gsig{\xi_1^\Sigma}   
\def\gsigep{\xi_{1,\ep}^{\Sigma}} 

\def\g{g_{\ep}}					
\def\Ca{\mathcal{C}_\a} 
\def\Cap{\mathcal{C}_{\a'}}

\def\ghpa{\pl^{\H}}       
\def\ghpaep{p_{1,\ep}^\H}
\def\ghpah{\tilde{\pl}^{\hspace{-0.5ex}\H}} 
\def\hep{\tilde{h}_{\ep}} 
\def\hnup{\tilde{h}_{\nup}}
\def\h{h}									

\def\eps{\varepsilon}		
\def\epso{\eps_0}				
\def\nuu{\nu}
\def\muu{\mu}
\def\ep{\underline{\eps}}  
\def\nup{\underline{\eps}'}

\newcommand{\crochet}[1]{\left[ #1 \right]}

\def\Fstarep{F_{\ep}^*}          
\def\Fstar{F^*}
\def\Fstaro{F_0^*}
\def\Fep{F_{\ep}}								
\def\F{F}
\def\Fo{F_0}

\def\Ftild{\widetilde{\F}}
\def\Ftildep{\widetilde{\F}_{\ep}}		

\def\xl{x_1}
\def\yl{y_1}
\def\xz{x_2}
\def\yz{y_2}											
\def\xy{(\xl,\yl,\xz,\yz)}
\def\x{x}													
							
\def\Rest{\mathcal{R}}						
\def\all{\alpha_1}
\def\alz{\alpha_2}								

\def\FFo{\mathcal{F}_0}						
\def\FFep{\mathcal{F}_{\ep}}      
\def\FF{\mathcal{F}}             

\def\roop{\rho'_0}			

\def\xil{\xi_1}
\def\xiz{\xi_2}
\def\etl{\eta_1}
\def\etz{\eta_2}									
\def\xiet{(\xil,\etl,\xiz,\etz)}
\def\xii{\xi_i}
\def\eti{\eta_i}

\def\philep{\varphi_{1,\ep}}
\def\psilep{\psi_{1,\ep}}
\def\phizep{\varphi_{2,\ep}}
\def\psizep{\psi_{2,\ep}}
\def\phiiep{\varphi_{i,\ep}}
\def\psiiep{\psi_{i,\ep}}   	
\def\philo{\varphi_{1,0}}
\def\psilo{\psi_{1,0}}
\def\phizo{\varphi_{2,0}}
\def\psizo{\psi_{2,0}} 					

\def\invphilep{\varphi_{1,\ep}^{-}}
\def\invpsilep{\psi_{1,\ep}^{-}}
\def\invphizep{\varphi_{2,\ep}^{-}}
\def\invpsizep{\psi_{2,\ep}^{-}}   

\def\Mo{\mathcal{M}_0}											
\def\M{\mathcal{M}}													

\def\HL{K}										
\def\HLep{K_{\ep}}

\def\flo{f_{1,0}}
\def\flep{f_{1,\ep}}
\def\glo{g_{1,0}}
\def\glep{g_{1,\ep}}
\def\fzep{f_{2,\ep}}
\def\gzep{g_{2,\ep}}							
\def\fiep{f_{i,\ep}}
\def\giep{g_{i,\ep}}

\def\alep{a_{1,\ep}}
\def\azep{a_{2,\ep}}
\def\blep{b_{1,\ep}}
\def\bzep{b_{2,\ep}}
\def\aiep{a_{i,\ep}}

\def\aleptilde{\tilde{a}_{1,\ep}}						%
\def\azeptilde{\tilde{a}_{2,\ep}}						%
\def\bleptilde{\tilde{b}_{1,\ep}}						%
\def\bzeptilde{\tilde{b}_{2,\ep}}						%
\def\aieptilde{\tilde{a}_{i,\ep}}						%
\def\bieptilde{\tilde{b}_{i,\ep}}						%
\def\phileptilde{\tilde{\varphi}_{1,\ep}}		%
\def\psileptilde{\tilde{\psi}_{1,\ep}}			%
\def\phizeptilde{\tilde{\varphi}_{2,\ep}}		%
\def\psizeptilde{\tilde{\psi}_{2,\ep}}      
\def\phiieptilde{\tilde{\varphi}_{i,\ep}}		%
\def\psiieptilde{\tilde{\psi}_{i,\ep}}

\def\philepstar{\varphi_{1,\ep}^*}		%
\def\psilepstar{\psi_{1,\ep}^*}			%
\def\phizepstar{\varphi_{2,\ep}^*}		%
\def\psizepstar{\psi_{2,\ep}^*} 						
\def\phiiepstar{\varphi_{i,\ep}^*}		%
\def\psiiepstar{\psi_{i,\ep}^*}	

\def\philostar{\varphi_{1,0}^*}		%
\def\psilostar{\psi_{1,0}^*}			%
\def\Mep{M_{\ep}}                      
\def\Aep{A_{\ep}}											
\def\Nep{N_{\ep}}

\def\S{S}													
\def\Wep{W_{\ep}}									

\def\PPhi{\mathbf{\Phi}}             
\def\Phil{\Phi_1}
\def\Psil{\Psi_1}
\def\Phiz{\Phi_2}
\def\Psiz{\Psi_2}
\def\Philt{\tilde{\Phi}_1}
\def\Psilt{\tilde{\Psi}_1}
\def\Phizt{\tilde{\Phi}_2}
\def\Psizt{\tilde{\Psi}_2}

\def\PPhiep{\mathbf{\Phi}_{\ep}}			
\def\Philep{\Phi_{1,\ep}}							
\def\Psilep{\Psi_{1,\ep}}
\def\Phizep{\Phi_{2,\ep}}
\def\Psizep{\Psi_{2,\ep}}
\def\Phiiep{\Phi_{i,\ep}}
\def\Psiiep{\Psi_{i,\ep}}

\def\Phiio{\Phi_{i,0}}
\def\Psiio{\Psi_{i,0}}

\def\PPhieptild{\mathbf{\widetilde{\Phi}}_{\ep}}			

\def\Phiieptild{\widetilde{\Phi}_{i,\ep}}
\def\Psiieptild{\widetilde{\Psi}_{i,\ep}}

\def\xii{\xi_{i}}
\def\eti{\eta_{i}}

\def\Schapl{\hat{S}_{1,\ep}}						
\def\Schapz{\hat{S}_{2,\ep}}
\def\Schapi{\hat{S}_{i,\ep}}

\def\fep{f_{\ep}}

\def\phileptch{\check{\varphi}_{1,\ep}}		%
\def\psileptch{\check{\psi}_{1,\ep}}			%
\def\phizeptch{\check{\varphi}_{2,\ep}}		%
\def\psizeptch{\check{\psi}_{2,\ep}} 						
\def\philotch{\check{\varphi}_{1,0}}		%
\def\psilotch{\check{\psi}_{1,0}}			%
\def\phizotch{\check{\varphi}_{2,0}}		%
\def\psizotch{\check{\psi}_{2,0}} 

\def\Siep{S_{i,\ep}}
\def\Slep{S_{1,\ep}}
\def\Szep{S_{2,\ep}}

\def\Sep{S_{\ep}}

\begin{document}

\title{Homoclinic orbits with many loops near a $0^2 i\omega$ resonant fixed point of Hamiltonian systems} 
\author{Tiphaine J\'ez\'equel\footnote{Laboratoire de Math\'ematiques Jean Leray, Universit\'e de Nantes, 2, rue de la Houssini\`ere - BP 92208 F-44322 Nantes Cedex 3, France. E-mail : {\tt tiphaine.jezequel@univ-nantes.fr}} ,
Patrick Bernard\footnote{CEREMADE, Universit\'e Paris - Dauphine,
Place du Mar\'echal De Lattre De Tassigny, 75775 Paris Cedex 16, France. E-mail : {\tt patrick.bernard@ceremade.dauphine.fr}}
 \hspace{1ex} and \'Eric Lombardi\footnote{Institut de Math\'ematiques de Toulouse, Universit\'e Paul Sabatier, 118 route de Narbonne, 31062 Toulouse cedex 9, France. E-mail : {\tt eric.lombardi@math.univ-toulouse.fr}}}

\date{\today}
\maketitle
\begin{center}

\begin{abstract}
In this paper we study the dynamics near the equilibrium point of a family of Hamiltonian systems in the neighborhood of a $0^2i\omega$ resonance. The existence of a family of periodic orbits surrounding the equilibrium is well-known and we show here the existence of homoclinic connections with several loops for every periodic orbit close to the origin, except the origin itself. To prove this result, we first show a Hamiltonian normal form theorem inspired by the Elphick-Tirapegui-Brachet-Coullet-Iooss normal form. We then use a local Hamiltonian normalization relying on a result of Moser. We obtain the result of existence of homoclinic orbits by geometrical arguments based on the low dimension and with the aid of a KAM theorem which allows to confine the loops. The same problem was studied before for reversible non Hamiltonian vectorfields, and the splitting of the homoclinic orbits lead to exponentially small terms which prevent the existence of homoclinic connections to exponentially small periodic orbits. The same phenomenon occurs here but we get round this difficulty thanks to geometric arguments specific to Hamiltonian systems and by studying homoclinic
orbits with many loops.

\vspace{1ex}

\end{abstract}
\end{center}

\section{Introduction}

In this paper, we study the dynamics near an equilibrium of a one-parameter family of real analytic Hamiltonian vector fields with two degrees of freedom. We suppose that these vector fields admit an equilibrium point that we take at the origin.

\begin{figure}[!h]
\unitlength=1cm
\begin{picture}(12,4)    \put(1,2){\line(1,0){2}}    \put(2,0.5){\line(0,1){3}}    \put(2,0.75){\circle*{.1}}    \put(2,3.25){\circle*{.1}}    \put(2.2,0.75){$-\I (\om_0+{\cal O}(\lb))$}    \put(2.2,3.25){$+\I (\om_0+{\cal O}(\lb))$}    \put(1.8,0){$\lb <0$}   
 \put(2,1.5){\circle*{.1}}    \put(2,2.5){\circle*{.1}}    \put(2.2,1.5){$+\I\lb$}    \put(2.2,2.5){$-\I\lb$}
 \put(5,2){\line(1,0){2}}    \put(6,0.5){\line(0,1){3}}    \put(6,0.75){\circle*{.1}}    \put(6,3.25){\circle*{.1}}    \put(6.2,0.75){$-\I\om_0$}    \put(6.2,3.25){$+\I \om_0$}    \put(5.8,0){$\lb=0$}    
 \thicklines    \put(5.8,1.8){\line(1,1){0.4}}     \put(5.8,2.2){\line(1,-1){0.4}}    \thinlines    \put(5.5,1.5){0}   \put(9,2){\line(1,0){2}}    \put(10,0.5){\line(0,1){3}}    \put(10,0.75){\circle*{.1}}    \put(10,3.25){\circle*{.1}}    \put(10.2,0.75){$-\I(\om_0+{\cal O}(\lb))$}    \put(10.2,3.25){$+\I(\om_0+{\cal O}(\lb))$}    \put(9.8,0){$\lb >0$}    
 \put(9.3,2){\circle*{.1}}    \put(10.7,2){\circle*{.1}}    \put(9.1,1.6){$-\lb$}    \put(10.6,1.6){$\lb$}
 \end{picture}\caption{\label{Bif2} Eigenvalues of $DV_{H_\lb}(0)$ in terms of $\lb$ for a $0^2\I\omega$ resonant equilibrium.}
 \end{figure}
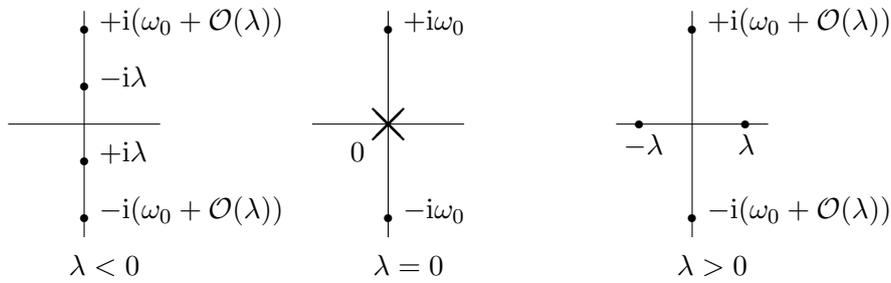

An equilibrium point $0$ of a vector field is called {\it non degenerated} if the linear part of the vector field is invertible. For the real Hamiltonian vector fields with two degrees of freedom, there exist only three types of non degenerated equilibria : the {\it Elliptic equilibria} when there are two pairs of conjugated purely imaginary eigenvalues, the {\it Saddle-Center equilibria} if there are one pair of conjugated purely imaginary eigenvalues and one pair of opposed real eigenvalues, and the {\it Hyperbolic equilibria} if there are two pairs of real opposed eigenvalues or four opposite and conjugated eigenvalues.
We study a family $H_{\lambda}$ of Hamiltonians with a fixed point at the origin, whose linear part  undergoes a transversal bifurcation at $\lambda=0$ through the stratum of degenerate fixed points, from an elliptic fixed point to a saddle center fixed point. We assume that the degenerate fixed point of $H_0$ admits a pair of null eigenvalues with a non-trivial Jordan block. This case, which is generic, is called an $0^2i\omega$ resonance (see Figure \ref{Bif2}).

Although we  are interested in the description of the dynamics associated to the saddle-center fixed point, we must distinguish two cases which are best described by considering the elliptic side of the bifurcation. Either the quadratic part of the Hamiltonian at the elliptic fixed point is definite, or it has index two. We will consider only the definite case. The homoclinic phenomenon  discribed in the present paper does not occur in the
other case.


We only study the dynamic for the "half-bifurcation" $\lb>0$, $\ie$ we study the dynamic in the neighborhood (in space) of a Saddle-Center fixed point in the neighborhood (in term of $\lb$) of a $0^2\I\omega$ resonance. Thus, this paper is related to two types of previous works : on one hand studies of Hamiltonian vector fields near Saddle-Center equilibria (see below part \ref{Intro2Ham} of this introduction), and on the other hand some works on reversible vector fields ($\ie$ vector fields $V$ which anticommute with some symmetry $S$ : $V\circ S=-S\circ V$) near a $0^2\I\omega$ resonance (see part \ref{Intro2Rev} below).

\subsection{$0^2\I\om$ resonance and waterwaves}\label{Intro2hydro}

This study is motivated historically by the waterwaves problem : the $0^2\I\omega$ resonance appears when one looks for two dimensional travelling waves for the Euler equation. The spatial dynamic method introduced by Kirchgassner \cite{Kirch} (see also for instance \cite{GerardMariana}), by a center manifold reduction, leads to the study of a four-dimensional reversible Hamiltonian vector field with an equilibrium passing through many different bifurcations in terms of the fluid parameters (the Bond number $b$ and the Froude number $\lambda$). For $b<1/3$ and $\lambda$ close to 1,   a $0^2\I\omega$ resonance occurs. 

In particular, a natural question is wether there exist solitary waves for this problem : in the reduced system this corresponds to {\bf ask wether there exists homoclinic orbits to the origin}. Results were only obtained close to resonance curves in the parameter plane $(b,\lambda)$. For $b>1/3$ and $\lambda$ close to 1, i.e. when Euler equations are well approximated by the KdV equation, the existence of true solitary waves was obtained in \cite{Waterwaves}. It was also obtained for $b<1/3$ and $\lambda$ close to a critical curve $\lambda=C(b)$ along which a Hamiltonian Hopf bifurcation occurs \cite{IoossPeroueme}.

\begin{figure}[!h]
\centering
\includegraphics[scale=0.4]{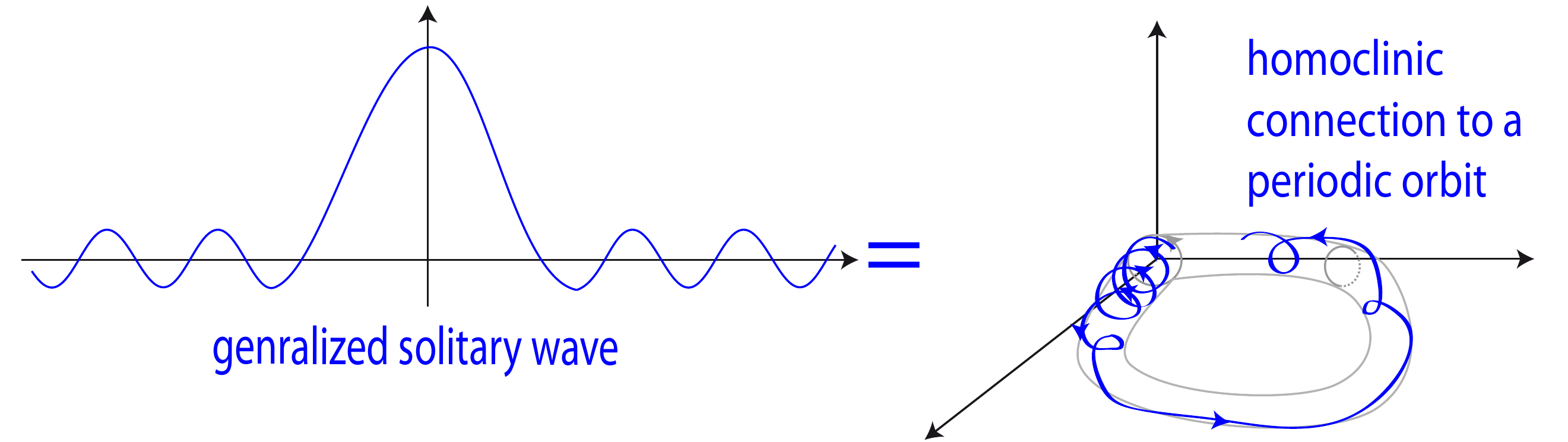}
\caption{Generalized solitary wave and homoclinic connection to a periodic orbit.}
\label{Vagues}
\end{figure}

For the $0^2\I\omega$ resonance the problem is far more intricate because it involves exponentially small issues. The existence of generalized solitary waves, which correspond to homoclinic connections to small periodic orbits (see Figure \ref{Vagues}) was obtained in this case \cite{Lombardi97}. The non existence of true solitary waves for $b<1/3$ and close to $1/3$ and for $\lambda$ close to 1 was proved by S.M. Sun in \cite{Sun2}. However for $b$ not close to $1/3$ the problem remains fully open.

\subsection{Reversible $0^2\I\om$ resonance}\label{Intro2Rev}

The $0^{2+}\I\omega$ resonance was extensively studied in the reversible context  (Iooss and Kirchgassner \cite{Waterwaves}, Lombardi \cite{Eric}, Iooss and Lombardi \cite{NoteCras})  with an additional assumption on the action of the symmetry on the eigenvectors of the linear part of the vector field (this assumption is satisfied in the waterwave problem). In this case, for $\lb>0$ sufficiently small, the Lyapunov-Devaney theorem ensures that the equilibrium is surrounded by a family of periodic solutions of arbitrary small size. Lombardi \cite{Eric} proved that 
there exists two exponentially small functions $\kappa_1(\lb)<\kappa_2(\lb)$ ($\ie$ $\kappa_i(\lb)={\cal O}(\E^{-\frac{c_i}{\lb}})$) such that on one hand the periodic orbits smaller than $\kappa_1(\lb)$ do not admit any homoclinic reversible connection {\it with one bump}, while on the other hand there exist homoclinic connections to each periodic orbit of size greater than $\kappa_2(\lb)$. Observe that in particular there is no homoclinic connection to the origin, but always homoclinic connections to exponentially small periodic orbits.

To obtain those results, the a key tool is the use of a normal form for which there exist homoclinic connections to all the periodic orbits. Then the study of the persistence of this dynamic which involves  exponentially small issues, was perform through a careful study of  the holomorphic continuation of solutions in the complex field.

\subsection{Homoclinic orbits to a Hamiltonian Saddle-Center equilibrium}\label{Intro2Ham}

In the Hamiltonian case, there are many works about the dynamic near a Saddle-Center equilibrium, without considering the neighborhood of a $0^22\I\omega$ singularity. Near an Hamiltonian Saddle-Center, the Lyapunov-Moser theorem ensures that there exists a family of periodic orbits surrounding the origin. Here also in most of the cases, the existence of homoclinic connections to periodic orbits is proved, but not to $0$.

A first group of results describes the {\bf consequences of the existence of on homoclinic connection to $0$}. For instance, in \cite{Ragazzo} and \cite{Ragazzo1}, Grotta Ragazzo studies the dynamic near an homoclinic orbit to the equilibrium, or in \cite{Salomao}, the authors obtain the existence of homoclinic connections to all the periodic orbits and some chaotic behavior for perturbations of a vector field admitting an homoclinic orbit to $0$. 

In another group of works, some results of {\bf existence of homoclinic connections to the periodic orbits near $0$} are obtained without making assumptions of existence of an homoclinic orbit to $0$ (for instance, \cite{Cascades} and \cite{Bernard}). But in fact the proofs often rely on assumptions allowing to prove that the vector field is in some sense a perturbation of a simpler case in which there is an explicit homoclinic orbit to $0$, so that the strategy is close to that of \cite{Salomao}.

\begin{figure}[!h]
\centering
\includegraphics[scale=0.27]{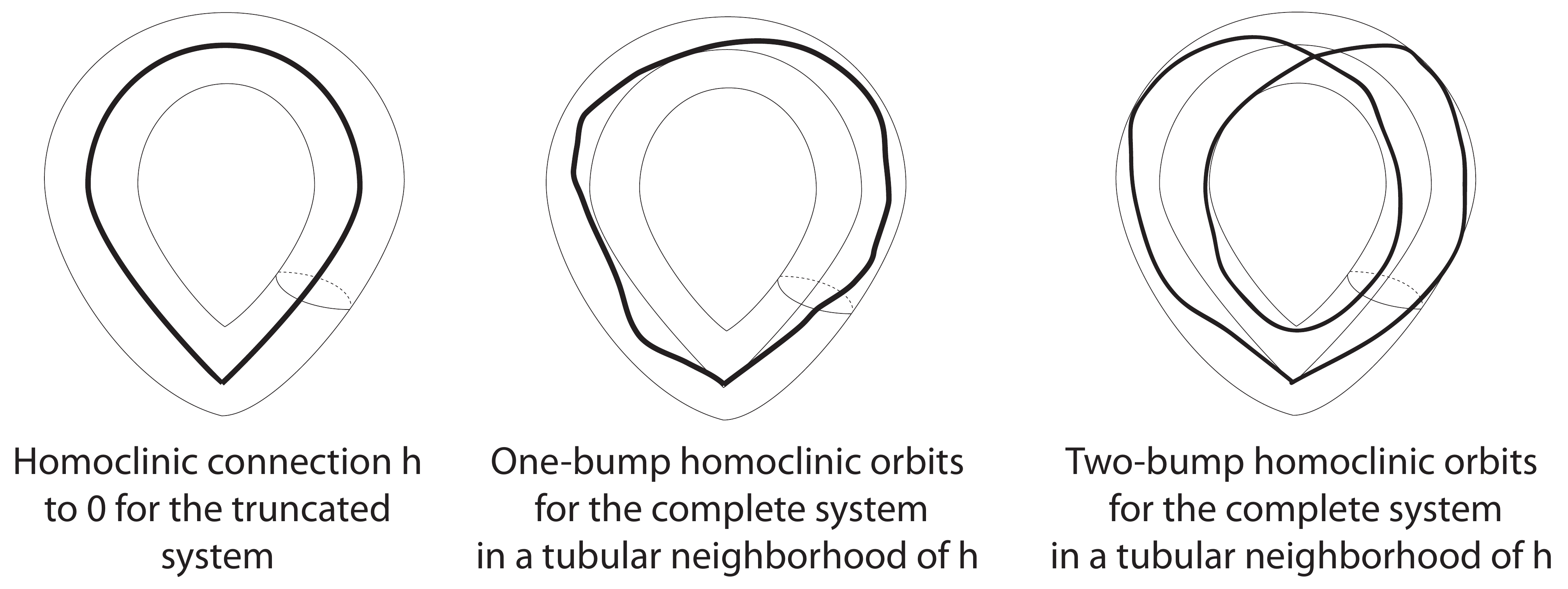}
\label{HomoclinesBoucles}
\end{figure}

Most of those papers use a method introduced by Conley \cite{Conley1} : this is a semi-global strategy which allows to construct homoclinic connections to periodic orbits as perturbations of an homoclinic orbit to $0$. This method is constructive and the homoclinic connections obtained may have "several loops" (see Figure \ref{HomoclinesBoucles}), while in the works on reversible systems described above the homoclinic connections considered were all with one loop. For that purpose, the idea of Conley is to introduced Poincar\'e sections transverse to the homoclinic orbit orbit to $0$ and to study the iterations of the Poincar\'e map. In \cite{Salomao}, the authors use an additional KAM argument to confine the iterations.

\vspace{-2ex}

\subsection{Homoclinic orbits with several loops for the Hamiltonian $0^2\I\omega$ resonance}

In this paper, {\it we consider  the Hamiltonian $0^2\I\omega$ resonance.} 
More precisely, let $\R^4$ be endowed with a symplectic form $\Om$. Consider a $\C^1$ one parameter family of analytic Hamiltonians $\Hlb$, where $\lb$ belongs to an interval $I$ of $\R$. We introduce the following normed space of analytic functions.

\begin{defi}\label{DefAn}
Let $\An(\B_{\R^m}(0,\rho),\R^n)$\label{NotAn} be the set of analytic functions $f:\B_{\R^m}(0,\rho)\rightarrow\R^n$ such that $f$ admits a {\it bounded} analytic continuation $\tilde{f}:\B_{\CC^m}(0,\rho)\rightarrow\CC^n$. We define the norm $\NNorme{\An}{\cdot}$ on $\An(\B_{\R^m}(0,\rho),\R^n)$ by
$$\NNorme{\An}{f}:=\sup_{z\in\B_{\CC^m}(0,\rho)}\NNorme{\CC^n}{\tilde{f}(z)}.\label{NotNormeAn}$$
\end{defi}
 
In the following we suppose that there exists $\roo$ such that 
$$\begin{array}{rcl}
\mathbf{H}:I&\to& \An(\B_{\R^4}(0,\roo),\R)\\
   \lb &\mapsto& \mathbf{H}_{\lb}
\end{array}
$$
is a $\C^1$ map for the norm $\NNorme{\An}{\cdot}$, $\ie$ we assume
\begin{equation}\tag{H0}\label{Hypreg}
\begin{array}{ll}
 &\mathbf{H}\in \C^1(I,\An(\B_{\R^4}(0,\roo),\R)). \\
\end{array}
\end{equation}

We study the associated family of Hamiltonian vector fields $\V := \Om\nabla\Hlb,$ which we suppose to admit a fixed point at the origin, $\ie$
\begin{equation}\tag{H1}\label{Hypzero}
\V(0)=0 \quad \text{ for all }\lb\in I.
\end{equation}

We assume that for $\lb=0$, the fixed point admits a \textbf{$0^2\I\om$ resonance}. This means that there exists a basis $(u_0,u_1,u_+,u_-)$ of $\CC^4$ in which  
\begin{equation}\tag{H2}
\label{Hyplin}
D_x\mathbf{V_{H_0}}(0)=\left(\begin{matrix}
    								0& 1 &  0 & 0\\
           					0 &0&0&0	 \\    
           					0&0&\I\om_0&0\\
           					0&0&0&-\I\om_0
    					\end{matrix}\right).
\end{equation}
We make an additional \textbf{assumption on $D_\lb(D_x\V(0))$} which characterizes the behaviour of the spectrum of $D_x\mathbf{V_{H_0}}(0)$ for $\lb=0$. Denote by $(u_0^*,u_1^*,u_+^*,u_-^*)$ the dual basis of $(u_0,u_1,u_+,u_-)$. We assume that
\begin{equation}\tag{H3}\label{Hypc1}
c_{10}:=\left\langle u_1^*,D^2_{x,\lb}\V(0).u_0\right\rangle\neq 0
\end{equation}
holds. This hypothesis will ensure that the spectrum of $D_x\mathbf{V_{H_0}}(0)$ is as represented in Figure \ref{Bif2} : we do not consider the hyper-degenerated case when the double eigenvalue $0$ stays at $0$ for $\lb\neq0$. Hence, the origin is an elliptic fixed point when $c_{10}\lb\leq0$ and a saddle-center when $c_{10}\lb\geq0$ (thus more precisely we get Figure \ref{Bif2} when $c_{10}>0$ ; if $c_{10}<0$ it holds for $\lb'=-\lb$).

\begin{sloppypar} Furthermore, we make the following \textbf{assumption on the quadratic part $D^2_{x,x}\mathbf{V_{H_0}}(0).[x,x]$ of the vector field at $\lb=0$},\end{sloppypar}
\begin{equation}\tag{H4}\label{Hypc2}
c_{20}:=-\frac{1}{3}\left\langle u_1^*,D^2_{x,x}\mathbf{V_{H_0}}(0).[u_0,u_0]\right\rangle\neq 0.
\end{equation}
This hypothesis ensures that, in some sens, the quadratic part of the vector field is not degenerated : $c_{20}$  will appear in the normal form as the coefficient of the only quadratic term. Thank to this nonzero term, we will be able to show the existence of homoclinic orbits for the normal form while the linearized vector field does not admit any such orbit.

In the following we focus our interest on the existence of homoclinic connections, so we only study the "half bifurcation" \begin{equation}\tag{H5}\label{HypCS}
c_{10}\lb\geq0,
\end{equation}
because in that case we work in the neighborhood of a \textbf{saddle-center fixed point}.

Finally, we assume that
\begin{equation}\label{HypElliptique}\tag{H6}
\omo>0
\end{equation}
holds, which means that for the small $c_{10}\lb<0$, the quadratic part of the Hamiltonian is a definite quadratic form.

Under these hypotheses we prove in this paper the following theorem which ensures that there exist homoclinic connections to all the periodic orbits surrounding the origin provided that the homoclinic connections are allowed to admit any number of bumps. This result is stronger than the one obtained in the reversible case since it ensures that when any number of bumps is allowed, there is no lower bound of the size of periodic orbits admitting an homoclinic connection to them.

\begin{thm}\label{TH}
Under the hypotheses $(H1),\ldots,(H6)$, there exist $\lb_0\geq0$, $\Czero\geq0$ and $\lzero\in\NN$ such that for all $\lb\leq\lb_0$,
\begin{enumerate}
\item the origin is surrounded by a family of periodic orbits $\P^\a_\lb$, labelled by their symplectic area $\a\in[0,\a_0]$ (Lyapunov-Moser);
\item for $\a\in ]0,\Czero\lb^{\lzero}[$, every periodic orbit $\P^\a_\lb$ admits an homoclinic connection.
\end{enumerate}
\end{thm}

Note that the Theorem \ref{TH} only deals with homoclinic connections to periodic orbits of arbitrary small size. It says {\it nothing on homoclinic connections to 0}. Their existence when several bumps are allowed remains fully open. The proof suggests (see Section \ref{StrategieDemonstration}) that the number of bumps of the homoclinic connection increases when the area of the periodic orbits decreases, if this really happens it might prevents the existence of an homoclinic connection to 0. However the theorem does not give any link between the number of bumps and the area.

Observe also that unfortunately, two obstacles prevent the result of Theorem \ref{TH} to hold in the case of the waterwaves problem. The first obstacle is the Hypothesis \eqref{Hypreg} : the $0^2\I\omega$ resonance appears in the waterwaves problem after a centermanifold reduction, and this reduction  does not preserve the analyticity of the initial equation. The second is Hypothesis \ref{HypElliptique} : in the waterwaves problem $\omo$ is negative, $\ie$ the quadratic part of the Hamiltonian is not definite. These two hypothesis are crucial in the proof : the analyticity to get a local linearizing change of coordinates without any remainder (the result of Appendix \ref{AppendixA}), and $\omo>0$ to confine the flow.

\vspace{2ex}

The {\it proof} of this Theorem relies on three main tools: normal forms, the method of Conley described above (in Part \ref{Intro2Ham}) and KAM theory to build an invariant curve which confine the iterations.

We prove in this paper a very general {\it Hamiltonian Normal Form theorem} (stated in Appendix \ref{AppendixA})
on the model proposed by Elphick \& al. \cite{NFcritere} for standard ODE : unlike the Birkhoff Normal Forms \cite{Birkhoff} which exist in the neighborhood of semisimple Elliptic fixed point, this theorem holds also for non Elliptic and non semisimple fixed points. 

The method of Conley relies on the construction of a return map. To build this map close to the equilibrium, we need to perform a second normalization, which "almost linearizes" the dynamics locally. This normalization is the one proposed by Moser in \cite{Moser2}. However we had to prove (see Appendix \ref{constr}) that this normalization does not blow up when $\lb$ goes to 0. Indeed the main difficulty in the $0^2\I\omega$ resonance is the competition between two scales (slow in the hyperbolic directions and fast in the elliptic directions), which leads to {\it exponentially small splitting} of the homoclinic connection. In the reversible case the proof requires a very careful description of the holomorphic continuation of the solution to be able to compute the exponentially small term. Here, after a scaling, the difficulty becomes a singularity in terms of $\lb$ in the elliptic directions (a fast rotation), and this requires a careful study of the change of coordinates of Moser \cite{Moser2} used in the method of Conley to be sure that it does not blow up when $\lb$ goes to 0.

\subsection{Plan of the paper}

The proof of this theorem is in section \ref{StrategieDemonstration} : but in fact this section contains only the mains steps of the strategy while the technical results are stated in propositions whose proofs are in the next sections \ref{PartieFormeNormale}, \ref{Poincare}, \ref{PartieReta} and \ref{SectionKAM}. 

In the Appendix \ref{AppendixA}, we state and prove a general Normal Form Theorem used in part \ref{PartieFormeNormale}. In Appendix \ref{tech}, we state some technical lemmas concerning an order relation on the formal power series, useful for the next appendix. The Appendix \ref{constr} is devoted to the proof of estimates about the dependence of the Moser's normalization in term of the parameters (used in part \ref{StrategieDemonstration}).


\label{PartStatement}

\subsubsection*{List of notations}
$$\begin{array}{lcc}
\An(\B_{\R^m}(0,\rho),\R^n)\ldots\ldots\ldots\ldots\ldots \pageref{NotAn} &\hspace{4ex}&\NNorme{\An}{\cdot}\ldots\ldots\ldots\ldots\ldots\ldots\ldots\ldots\ldots\ldots\pageref{NotNormeAn}\\
\roo \dotfill \pageref{Notroo}&& \Iz \dotfill \pageref{NotIz}\\
\ep,\nuu,\muu \dotfill \pageref{Notep} && \Sigl, \dll \dotfill \pageref{NotSigl}\\
\roop \dotfill \pageref{Notroop} && \FFep, \Mo \dotfill \pageref{NotFFep}\\
\HLep \dotfill \pageref{NotHLep} && \Sigz \dotfill \pageref{NotSigz}\\
\Ret \dotfill \pageref{PropDefRetp} && \P^\a_{\ep} \dotfill \pageref{DefPaep}\\
\gcsep,\ghpa \dotfill \pageref{PropReta}&& \C^\a_s \dotfill \pageref{PropReta}\\
\Ret^\a_{\ep} \dotfill \pageref{PropReta} && \Gamma^\a_\eps \dotfill \pageref{propGamma} \\
\C^\a_u\dotfill \pageref{NotCai}&&\gsigep \dotfill \pageref{LemGrapheSigl}\\
\Rot_\theta, \rot_\theta \dotfill \pageref{LemFlotp} && k_0\dotfill     \pageref{LemFlotp} \\
\nup \dotfill \pageref{Notnup} && \hep \dotfill \pageref{LemGrapheHPa} \\
\prec \dotfill \pageref{DefPrec} && [\cdot] \dotfill \pageref{NotCrochet} 
\end{array}$$


\setlength{\parindent}{0cm}

\section{Structure of the proof of Theorem \ref{TH}}\label{StrategieDemonstration}
This section is devoted to the proof of Theorem \ref{TH}. Some technical steps of this proof are stated in propositions whose proofs are postponed in the next sections of the paper.

\subsection{Normalization and scaling, dynamics of the normal forms of degree 3 and $n$}\label{SubNormalisation}\label{subsubFN3}\label{subsubFNn}\label{subsubStrategie}\label{subsubchange}

\subsubsection*{Normalization and scaling}

Proposition \ref{PropNF} below gathers the results of normalization an scaling of the Hamiltonian. Point $(i)$ is the change of coordinates given by the normal form Theorem \ref{ThmNF}. Then $(iii)$ is a scaling in space and time : after this scaling, the homoclinic orbits of the truncated normalized system have a size of order 1 (see the next subsections), which allows a perturbative proof in the neighborhood of these homoclinic orbits. We must perform the change of parameter $(ii)$ to have a $\C^1$ smoothness of scaling $(iii)$.

\begin{prop} \label{PropNF}

Under hypotheses (H0),$\cdots$,(H6), for all $n\geq 3$, there exist $\eps_1>0$, $\rho_1>0$ and
\begin{enumerate}
\item a $\C^1$ one parameter family of canonical analytic transformations of $\B(0,\rho_1)$, $$\xbf=\phi_{n,\lb}(\xtilde)=\xtilde+{\cal O}(|\xtilde|^2), \qquad \phi_{n,\lb}\in\An(\B_{\R^4}(0,\rho_1),\R^4);$$
\item a change of parameter $\C^1$, $\lb=\theta(\eps^4)$ with its inverse $\eps^4 = c_{10}\lb+o(\lb)$ defined for $\eps\in]-\eps_1,\eps_1[$,
\item a scaling in space and time $\xtilde=\sigma_\eps(\xp)$, $t =\eps^2 \ \widetilde{t}$;
\end{enumerate}

\vspace{1ex}

such that in a neighborhood of the origin, for all $\eps\in]-\eps_0,0[ \cup ]0,\eps_0[$, the normalized Hamiltonian $\Hprime_{\eps}(\xp) :=\frac{8 c_2(\eps)^2}{\eps^{12}} \Hbf_{\theta(\eps^4)}\left(\phi_{n,\theta(\eps^4)}\bigl(\sigma_\eps(\xp)\bigr)\right)$ is of the form
\begin{equation}
\Hprime_\eps(\xp)=\frac{1}{2}\left(\plp^2-\qlp^2\right)+2\sqrt{2}(\qlp)^3+\frac{\omega(\eps)}{2\eps^2}\Izp+\eps^2 \Nprime_{n,\eps}(\qlp,\Izp)+ \eps^{4n-8} \Rprime_{n,\eps}(\xp); \nonumber
\end{equation}
where
$$\xp=(\qlp,\plp,\qzp,\pzp), \qquad \Izp = \qzp^2+\pzp^2,$$
and $\om$ and $c_2$ are $\C^1$ functions of $\eps\in]-\eps_0,\eps_0[$ such that $\om(\eps)=\omo+{\cal O}(\eps^4)$ et $c_2(\eps)=c_{20}+{\cal O}(\eps^4)$ with $c_{20}\neq 0$. 

The normal form $\Nprime_{n,\eps}$ is a real polynomial of degree less than $n$ in $(\qlp, \qzp,\pzp)$ such that
$$\Nprime_{n,\eps}(\qlp,\Izp)={\cal O}\left( |\qlp||\Izp|+\eps^2(|\qlp|^2+|\Izp|)^2\right ),$$
and the coefficients of this polynomial are $\C^1$ functions in $\eps\in]-\eps_0,\eps_0[$. The remainder $\Rprime_{n}:]-\eps_0,\eps_0[\ \to \An(\B_{\R^4}(0,\rho_1);\R) : \eps \mapsto \Rprime_{n,\eps}$ is $\C^1$ one parameter family of analytic Hamiltonians satisfying $\Rprime_{n,\eps}(\xp)={\cal O}(|\xp|^{n+1})$.
\end{prop}

This proposition is proved in Section \ref{PartieFormeNormale}. From now on, we work with the Hamiltonian $\Hprime_\eps$.

\subsubsection*{Phase portrait for the normal form of degree 3} 

We first study the dynamics of the Hamiltonian $\Hprime_\eps$ truncated at degree 3
$$\Hprime_{3,\eps}(\xp)=\frac{1}{2}\left(\plp^2-\qlp^2\right)+2\sqrt{2} \ \qlp^3+\frac{\omega(\eps)}{2\eps^2} \ \Izp,$$
where $\Izp=\qzp^2+\pzp^2$. The associated differential system reads
\begin{equation}\nonumber
\left\{\begin{array}{rcl}
\qlp'(t)&=&\plp \\
\plp'(t)&=&\qlp-6\sqrt{2} \ \qlp^2 \\
\qzp'(t)&=&\frac{\om(\eps)}{2\eps^2} \ \pzp\\
\pzp'(t)&=&-\frac{\om(\eps)}{2\eps^2} \ \qzp.
\end{array}\right.
\end{equation}
We observe that in this system, the two couples of variables $(\qlp,\plp)$ and $(\qzp,\pzp)$ are uncoupled. The solutions of the half system in $(\qzp,\pzp)$ are
$$\left(\begin{matrix} \qzp \\ \pzp \end{matrix}\right)=\rot_{\frac{\om(\eps)}{2\eps^2}t}\left(\begin{matrix} \qzp(0) \\ \pzp(0) \end{matrix} \right), \qquad \text{ where } \quad \rot_{\theta}:=\left(\begin{matrix} \cos\theta & - \sin\theta \\ \sin\theta & \cos\theta \end{matrix}\right).$$
In particular $\Izp$ is constant. Let us denote by {\bf $\P^\a_0$  the periodic orbits satisfying $\Izp=\a$ and $(\qlp,\plp)=(0,0)$.} We get the phase portrait for the half system in $(\qlp,\plp)$ by drawing the energy level sets which read
\begin{equation}\label{IntroAlpha}
\left\{(\qlp,\plp), \frac{1}{2}\left(\plp^2-\qlp^2\right)+2\sqrt{2} \ \qlp^3=\alpha\right\}=\left\{(\qlp,\plp), \plp=\pm\sqrt{\qlp^2-4\sqrt{2} \ \qlp^3+\alpha}\right\}.
\end{equation}

\begin{figure}[!h]
\centering
\includegraphics[scale=0.4]{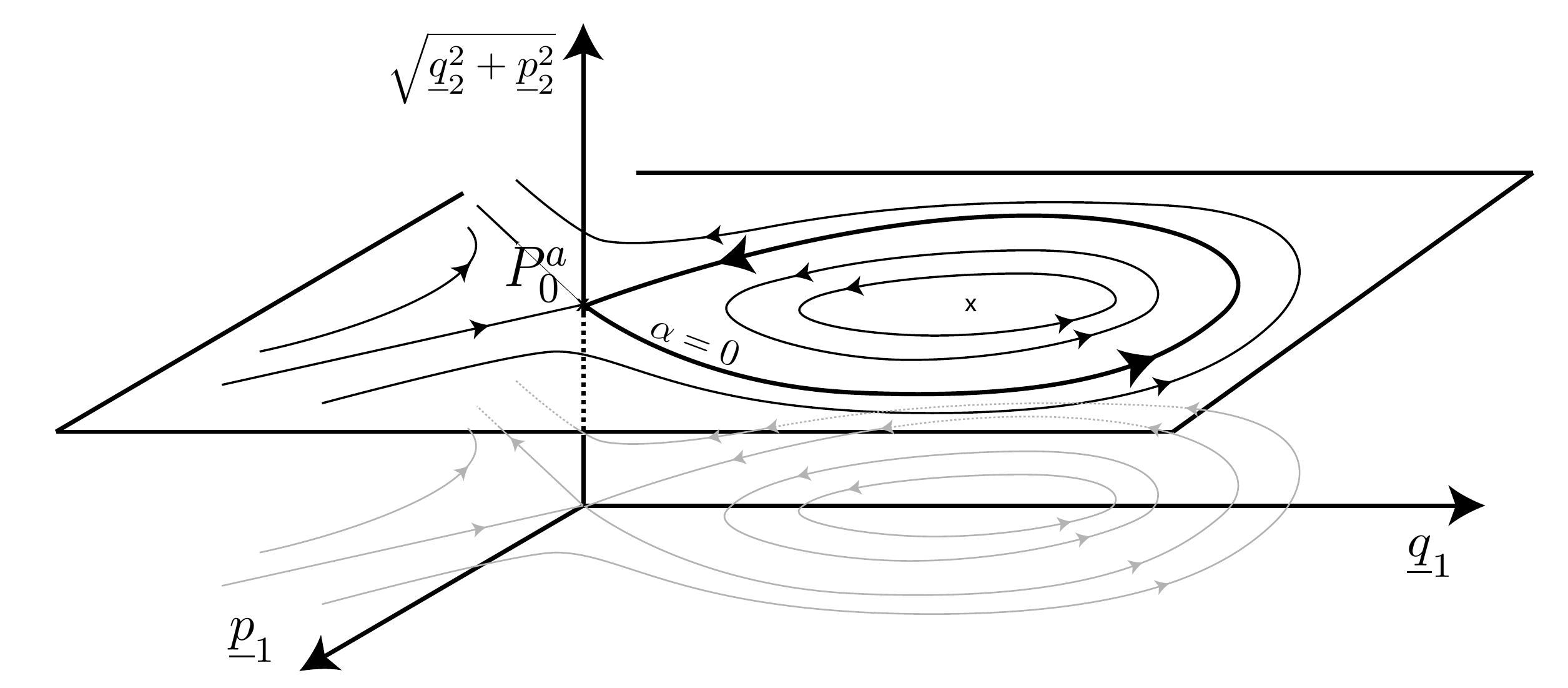}
\caption{Phase portrait for the normal form of degree 3.}
\label{DESSIN12}
\end{figure}

%
%
%
%
We then get the phase portrait of Figure \ref{DESSIN12}, in which there is an {\bf homoclinic orbit to $\P^\a_0$} when $\alpha=0$, whose explicit form is

$$\left(\begin{matrix} \qlp(t) \\ \plp(t) \end{matrix}\right)=\frac{1}{\sqrt{2\sqrt{2}}}\left(\begin{matrix} \Frac{1}{1+\cosh(t)} \\ \Frac{-\sinh(t)}{(1+\cosh(t))^2}\end{matrix}\right).$$ 

For the normal form of degree $n$, given that $\frac{d\Izp}{dt}=0,$ we also get the entire phase portrait, which is a deformation of the phase portrait (see Figure \ref{DESSIN12}) of the normal form of degree 3.

\subsubsection*{Linear change of coordinates and truncation at infinity}
From now on, it will be easier to work in the {\bf new coordinates $\qp$}, obtained by the following canonical linear change of coordinates in which the linearized hamiltonian system at the origin is in Jordan form :
\begin{equation}
\left(\begin{matrix} \qlp \\ \plp \end{matrix}\right)=\frac{1}{\sqrt{2}}\left(\begin{matrix} 1 & 1 \\ -1 & 1 \end{matrix}\right)\left(\begin{matrix} \ql \\ \pl \end{matrix} \right):= L \left(\begin{matrix} \ql \\ \pl \end{matrix} \right), \qquad  \left(\begin{matrix} \qzp \\ \pzp \end{matrix}\right)=\left(\begin{matrix} \qz \\ \pz \end{matrix}\right).\nonumber
\end{equation}
In these coordinates the Hamiltonian reads
\begin{equation}\label{hamiltonian}
-\ql\pl+(\ql+\pl)^3+\frac{\omega(\eps)}{2\eps^2}\Iz+\eps^2 \Q(\ql+\pl,\Iz,\eps)+ \eps^{4n-8} \Reste(x,\eps),
\end{equation}
where $\Iz=\qz^2+\pz^2$ \label{NotIz}, $\om$ is a $\C^1$ function of $\eps\in]-\eps_0,\eps_0[$ which reads $\om(\eps)=\omo+{\cal O}(\eps^4)$ ; the normal form $\Q$ is a real polynomial of degree less than $n$ in $\qp$ satisfying
$$\Q(\ql\+\pl,\Iz,\eps)={\cal O}\left( |\ql\+\pl||\Iz|+\eps^2(|\ql\+\pl|^2+|\Izp|)^2\right),$$
and the coefficients of which are $\C^1$ functions of $\eps\in]-\eps_0,\eps_0[$ ; the remainder $\Reste:]-\eps_0,\eps_0[\ \to \An(\B_{\R^4}(0,\rho_1);\R) : \eps \mapsto \Reste(\cdot,\eps)$ is a $\C^1$ one parameter family of analytic Hamiltonians such that $\Reste(\x,\eps)={\cal O}(|x|^{n+1})$.

\vspace{2ex}

We cut this Hamiltonian to get a bounded flow. Let $\mathcal{T}_{\roo}$ be a $\C^{\infty}$ map of $\R$ such that
$$\mathcal{T}_{\roo}(r)=\left\{\begin{array}{l}
 															1 \quad\text{for }r\in[-\frac{1}{4}\roo^2,\frac{1}{4}\roo^2]\\
 															\\
 															0 \quad\text{for }|r|\geq\roo^2.
															\end{array}\right.$$
We chose $\roo$  \label{Notroo} such that the homoclinic orbit obtained above for the truncated system is strictly included in $\B(0,\frac{1}{2}\roo)$. We finally consider the following Hamiltonian,
{\small $$\H(x,\eps):=\mathcal{T}_{\roo}(\ql^2)\mathcal{T}_{\roo}(\pl^2)\mathcal{T}_{\roo}(\Iz)\big(\moins\ql\pl\+(\ql\+\pl)^3+\frac{\omega(\eps)}{2\eps^2}\Iz+\eps^2 \Q(\ql\+\pl,\Iz,\eps)\+\eps^{4n-8} \Reste(x,\eps)\big),$$
}
which, in $\B(0,\frac{1}{2}\roo)$ is equal to the Hamiltonian (\ref{hamiltonian}) obtained above. This truncation is useful to work with a bounded flow, which cannot get out of $\B(0,\roo)$ : this will be useful to obtain uniform upper bounds. And then with the aid of these upper bounds we will get that, for $\eps$ sufficiently small, the solutions of interest stay in $\B(0,\frac{1}{2}\roo)$ : these solutions will also be solutions of the initial Hamiltonian (\ref{hamiltonian}).

{\it In the following, we always work with the Hamiltonian without making mention of the cutoff function $\mathcal{T}_{\roo}(r)$, given that we will always work in $\B(0,\frac{1}{2}\roo)$}.

\subsubsection*{New parameters $\ep, \nuu, \muu$}

We introduce new parameters for the Hamiltonian, $\ep:=(\eps,\nuu,\muu)$ \label{Notep} so that
\begin{equation}\label{3parametres}
\H(x,\ep)=-\ql\pl+(\ql+\pl)^3+\frac{\omega(\eps)}{2\eps^2}\Iz+\nuu \Q(\ql+\pl,\Iz,\eps)+ \muu\nuu\eps^{\No} \Reste(x,\eps).
\end{equation}
With these new parameters, we have
\begin{itemize}
\item for $(\eps,\nuu,\muu)=(\eps,0,0)$ the Hamiltonian is the normal form of degree 3 studied above in subsection \ref{subsubFN3};
\item for $(\eps,\nuu,\muu)=(\eps,\eps^2,0)$ we have the normal form of degree $n$ ;
\item for $(\eps,\nuu,\muu)=(\eps,\eps^2,\eps^{4n-8-(\No+2)})$ the complete system.
\end{itemize}
{\it The distance from the complete system to the normal form of degree 3 and the normal form of degree $n$ corresponds then to the smoothness in the parameters $\nuu$ (for degree 3) and $\muu$ (for degree $n$), uniformly in $\eps$}. We chose $\No$ and $n$ later in the proof.

\subsubsection*{Introducing heuristically our strategy of proof} 

At every order $n$, the normal form has homoclinic connections to the origin and to each periodic orbit of the family surrounding the origin. Moreover, if one deflects from the homoclinic trajectory "inward", one arrives in a region of space filled with trajectories periodic in $(\qlp,\plp)$. In the following, we consider the complete system as a perturbation of the normal form. The heuristic idea is then to show that if the homoclinic connection to a periodic orbit $\P^\a$ is perturbated, necessarily it deflects "inward", and then "follows" a trajectory periodic in $(\qlp,\plp)$, maybe making several loops and then finally joins the periodic orbit $\P^\a$ back. Such a trajectory for the complete system would then be an homoclinic orbit with several loops.

To give a mathematical sense to the idea of "doing several loops", a natural idea is then to introduce an appropriate Poincar\'e section intersecting transversally the homoclinic orbits of the normal form and to consider iterations of the first return map to this section. {\it This is our strategy in the following.}

\subsection{Construction of the first return map}\label{subsubFF}\label{subsubRet2}\label{subsubRet}

{\it Let us introduce the Poincar\'e section
\begin{equation}
\Sigl:=\left\{\qp\in\R^4/ \ql=\dll\right\},\label{NotSigl} \nonumber
\end{equation}}
where $\dll$ is fixed "small" : we will have several conditions of smallness on $\dll$ in the following, but none linked to the size of $\ep$.

To construct the first return map to $\Sigl$, denoted by $\Ret_{\ep}$, we proceed in two main steps, that we summarize here and detail below :
\begin{description}
\item[Global map.] We use the existence of an orbit homoclinic to the origin for the normal form of degree 3 and show by perturbation that there exists a return to the section for  $\a$ and $\nuu$ small. 
The perturbative method works only on a part of the homoclinic orbit which is covered in a finite period of time : we perform this strategy for {\bf a trajectory from a second section $\Sigz$ to the section $\Sigl$} (see Figure \ref{DESSIN10}).
\item[Local map.] We chose $\Sigz$ and $\Sigl$ close to the origin, where the homoclinic trajectory is covered in an infinite period of time. In order to construct {\bf a local map from $\Sigl$ to $\Sigz$} we have to build a local change of coordinates $\FFep$ in the neighborhood of $0$ which nearly linearizes the flow and then allows to show that the trajectories coming from $\Sigl$ do intersect $\Sigz$ (see Figure \ref{DESSIN11}).
\end{description}


\subsubsection*{Construction of a local change of coordinates, $\FF_{\protect\underline{\varepsilon}}$} 

\begin{prop} \label{propFF}
There exist $\epso$, $\roop<\frac{1}{2}\roo$ \label{Notroop} and a family of canonical analytic changes of coordinates
$$\FFep\xiet:=(\philep,\psilep,\phizep,\psizep)\xiet=\xiet+{\cal O}(|(\xi,\eta)|^2),\label{NotFFep}$$ 
defined for $|\ep=(\eps,\nuu,\muu)|\leq\epso$ such that the Hamiltonian $\H$ defined by \eqref{3parametres} \label{NotHLep} reads in the new coordinates $\xiet$
\begin{equation}\label{DefHL1}
\H\left(\FFep(\xiet),\ep\right)=\HLep(\xil\etl,\xiz^2+\etz^2)=-\xil\etl+\frac{\om(\eps)}{2\eps^2}(\xiz^2+\etz^2)+{\cal O}(|(\xil\etl,\xiz^2+\etz^2)|^2),
\end{equation}
and for all $\ep$, $\FFep$, $\FFep^{-1}$ $\in\An(\B_{\R^4}(0,\roop),\R^4)$. 

Moreover, $\FF_{(\eps,0,0)}:=\FFo$ does not depend on $\eps$, and there exists $\Mo\in\R$ such that
\begin{equation}\label{EstimFFep}
|\FFep\xiet-\FFo\xiet|\leq\nuu\Mo
 \end{equation}
for $|\ep|\leq\epso$. All the ${\cal O}$ correspond to upper bounds independent of $\ep$.
\end{prop}

The proof of this proposition and a more detailed statement are in Appendix \ref{constr} (Proposition \ref{propFFcomplete}).


\begin{rem}
This proposition is a fundamental step of the proof. Unfortunately, it requires very long and technical computations to obtain the estimate \eqref{EstimFFep}. The main interest (and main difficulty in the proof) of this proposition is to deal with the singularity in $\eps$ of the initial Hamiltonian $\H$ (see the explicit form \eqref{3parametres}) and to verify that despite of this singularity the estimate \eqref{EstimFFep} is uniform in term of $\eps$ small. See also the more detailed version (Proposition \ref{propFFcomplete}) of this proposition in Appendix \ref{constr} and the comments therein.
\end{rem}

\begin{figure}[!h]

\begin{minipage}[b]{60mm}
\includegraphics[scale=0.40]{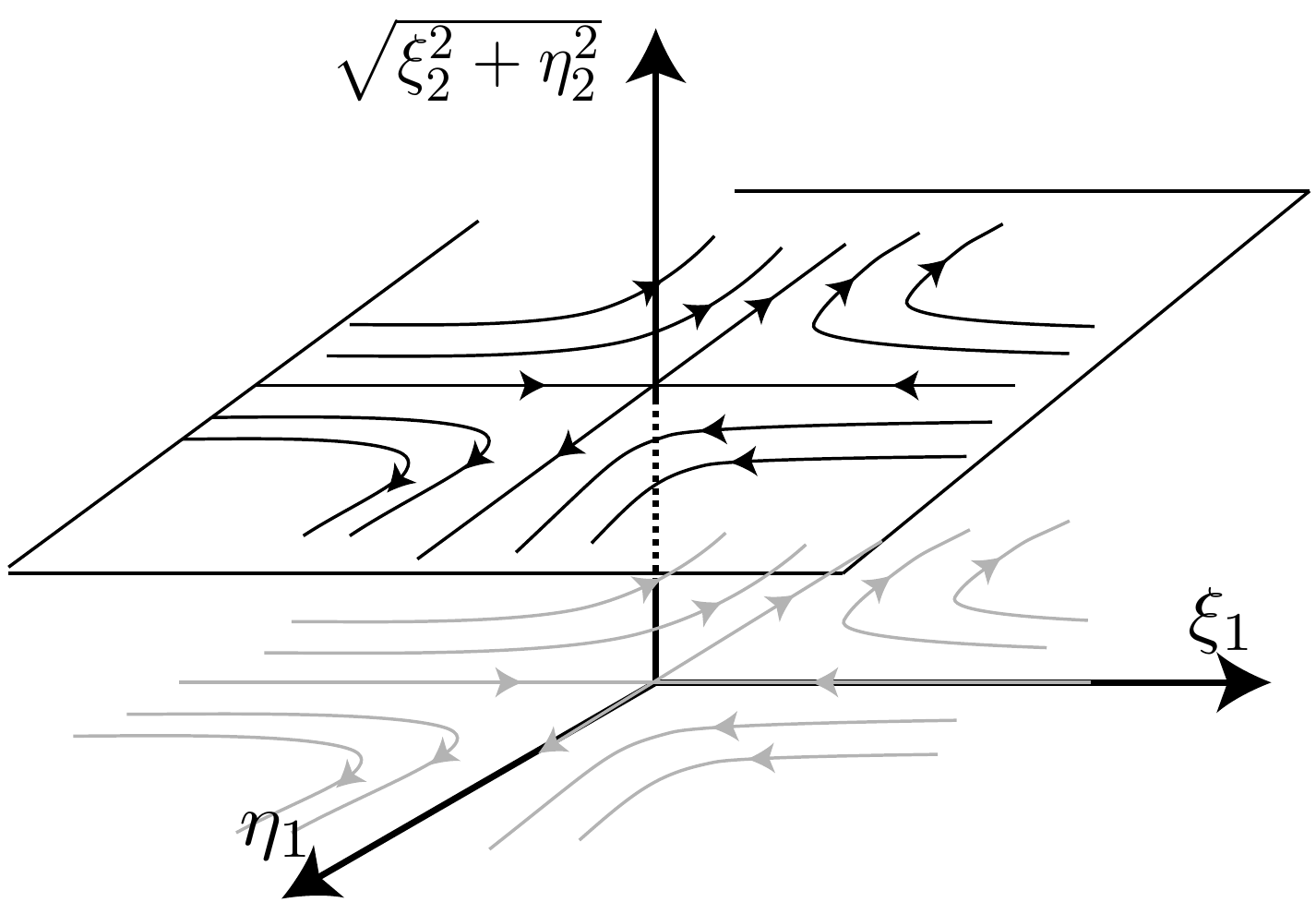}
\caption{Phase portrait in local coordinates $\xiet$.}
\label{DESSIN9}
\end{minipage}
\begin{minipage}{5mm}
\hfill
\end{minipage}
\begin{minipage}[b]{75mm}
\includegraphics[scale=0.39]{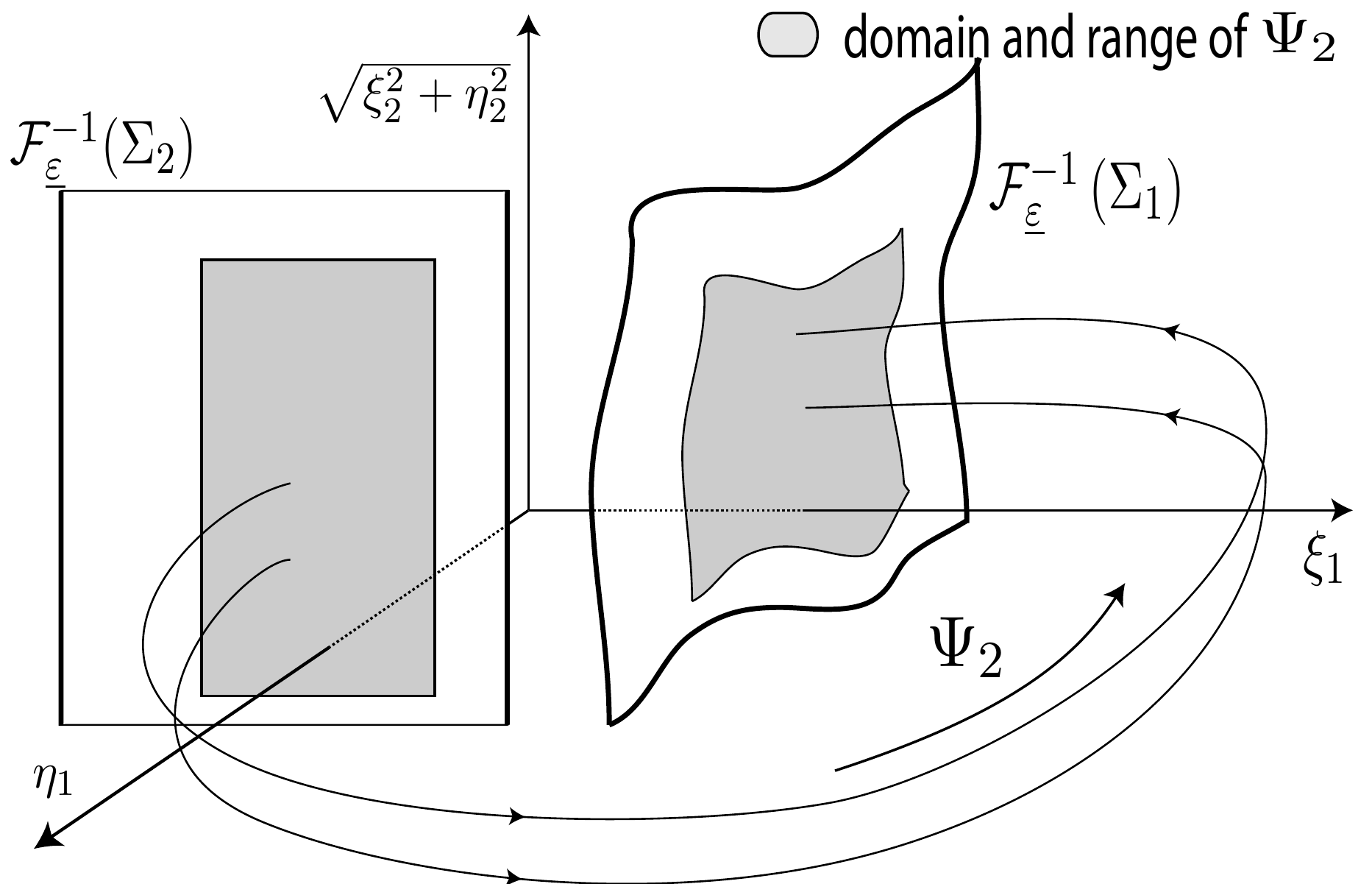}
\caption{Global return map $\Ret_2$.}
\label{DESSIN10}

\vspace{3ex}

\end{minipage}

\end{figure}

The form of the Hamiltonian obtained in the new coordinates allows to get the entire phase portrait. Indeed, the flow of the associated hamiltonian system satisfies
$$\frac{d(\xil\etl)}{dt}=0, \quad \frac{d(\xiz^2+\etz^2)}{dt}=0. $$
Moreover,
$$\frac{d\xil}{dt}=(\partial_1\HLep)(\xil\etl,\xiz^2+\etz^2)\cdot\xil$$
holds and, given that the ${\cal O}$ of equation (\ref{DefHL1}) is independent of $\ep$, we get
\begin{eqnarray}
\left|\HLep(\xil\etl,\xiz^2+\etz^2)-\left(-\xil\etl+\frac{\om(\eps)}{2\eps^2}(\xiz^2+\etz^2)\right)\right|\leq\Mo|(\xil\etl,\xiz^2+\etz^2)|^2.
\end{eqnarray}
Then, up to a reduction of the radius $\roop$ if necessary (independent of $\ep$), $\frac{d\xil}{dt} <0$ if $\xil>0$ and $\frac{d\xil}{dt}>0$ if $\xil<0$. So, the dynamics in the local coordinates $\xiet$ is as draw in the phase portrait of Figure \ref{DESSIN9}.

\subsubsection*{Global map, from a second section $\Sigz$ to $\Sigl$} 

{\it We define $\Sigz$ by defining its range in coordinates $\xiet$,
\begin{equation}
\Sigz:=\FFep\left(\{\xiet\in\B_{\R^4}(0,\roop)/ \etl=\dll\}\right). \label{NotSigz}\nonumber
\end{equation}}

By a perturbative method in the neighborhood of the homoclinic orbit to the origin of the normal form of degree 3, we show the existence of a Poincar\'e map $\Ret_2$ following the flow from $\Sigz$ to $\Sigl$. Precisely, denoting by $\flot(t,x,\ep)$ the flow of the hamiltonian system associated to$ H(\cdot,\ep)$, we show the following proposition, which gives moreover the upper bound (\ref{MajRet2}) useful later :

\begin{prop} \label{PropRet2}
For $\dll$ sufficiently small, there exist $T^-(\dll)<T^+(\dll)$ such that for all $\qp$ in
$$\Sigz\cap\FFep\left(\left\{\xiet/0\leq\xil\leq\frac{1}{16}\dll, \sqrt{\xiz^2+\etz^2}\leq\frac{1}{2}\dll\right\}\right),$$
there exists an unique $T_L(\qp,\ep)\in[T^-(\dll),T^+(\dll)]$ satisfying
$$\flot(T_L(\qp,\ep),\qp,\ep)\in\Sigl\cap\B(0,\dll).$$
Moreover, denoting by 
\begin{eqnarray}
\Ret_2(\qp,\ep)&:=&\flot(T_L(\qp,\ep),\qp,\ep)\nonumber
\end{eqnarray}
there exists $M_2$ such that 
\begin{equation}\label{MajRet2}
\left|\Ret_{2,\qz}(\qp,\ep)^2\+\Ret_{2,\pz}(\qp,\ep)^2\moins(\qz^2\+\pz^2)\right|\leq\muu\nuu\eps^{\No} M_2 T^+(\dll).
\end{equation}
\end{prop}

The proof of this proposition is in section \ref{SubPropDefRet}.

\subsubsection*{Existence of the first return map to $\Sigl$}

\begin{figure}[!h]
\centering
\includegraphics[scale=0.35]{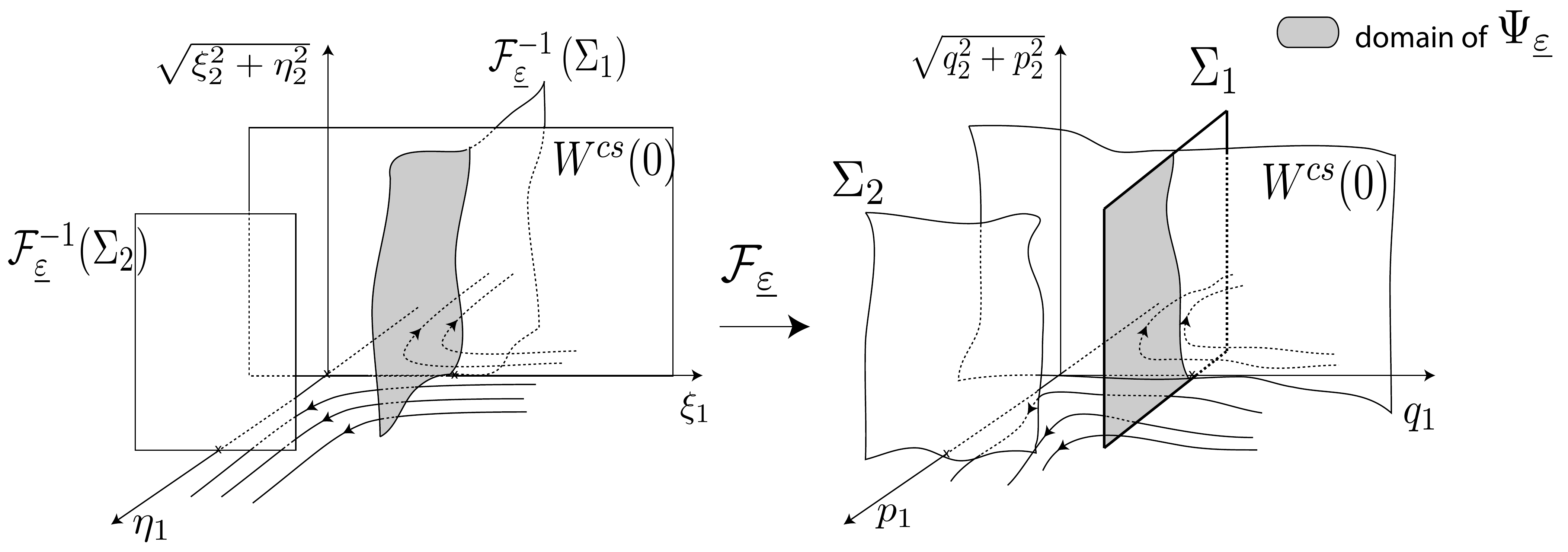}
\caption{Position of the domain of the Poincar\'e map and $\CS$.}
\label{DESSIN11}
\end{figure}

Observing the phase portrait in the neighborhood of the origin in the local coordinates $\xiet$ (Figure \ref{DESSIN9}), we see that the Poincar\'e map from $\Sigl$ to $\Sigz$ exists if and only if $\etl$ is positive. Given that in these coordinates, the center-stable manifold to the origin is the hyperplane $\{\etl=0\}$, we get that in coordinates $\qp$, the domain of existence corresponds to being "on the right side" of the center-stable manifold of the origin (see Figure \ref{DESSIN11}).

Precisely, we show the following
\begin{prop} \label{PropDefRetp}
It is possible to define a first return map
$$\Ret:\Sigl\cap\FFep\left(\{\xiet/0<\etl\leq\frac{1}{24}\dll,\sqrt{\xiz^2+\etz^2}\leq\dll\}\right)\times]-\epso,\epso[^3 \rightarrow \Sigl.$$
\end{prop}

A more detailed version of this proposition and its proof are given in part \ref{SubPropDefRet} (Proposition \ref{PropDefRet}).

\subsection{Construction of an invariant curve for the first return map using a KAM theorem}\label{SubStrategieKAM}\label{subsubinvariantcurve}

{\it From now on, we will not need to distinguish $\nuu$ from $\eps$, we work now with
$$\ep=(\eps,\nuu,\muu)=(\eps,\eps^2,\muu).$$}
In this section, we fix one periodic orbit $P$ and prove that the restriction of $\Ret$ to the energy level set of $P$ can be expressed as a diffeomorphism of an annulus of $\R^2$. Then we construct an invariant curve for this diffeomorphism with the aid of a KAM theorem. This curve will be useful later in part \ref{StrategieConclusion} to bound the iterations of the map $\Ret$, and then to conclude that $\WIP$ and $\WSP$ intersect each other.

\subsubsection*{The maps $\Reta$ and their expression as diffeomorphisms of an annulus of $\R^2$}

Thank to the canonical change of coordinates $\FFep$ of Proposition \ref{propFF}, we have a precise labelling of the periodic orbits in the neighborhood of 0. Indeed, in coordinates $\xiet$ we have the family of periodic orbits
$$\{(0,0,\xiz,\etz)/ \xiz^2+\etz^2=\a\},$$
labelled by their symplectic area $\a$. {\it Then, we denote
\begin{equation} \label{DefPaep}
\P^\a_{\ep}:=\FFep\left(\{(0,0,\xiz,\etz)/ \xiz^2+\etz^2=\a\}\right).
\end{equation}}
In particular, given that $\FFep$ is canonical, the symplectic area of $\P^\a_{\ep}$ is also $\a$.

\begin{figure}[!h]

\begin{minipage}[b]{67mm}
\includegraphics[scale=0.35]{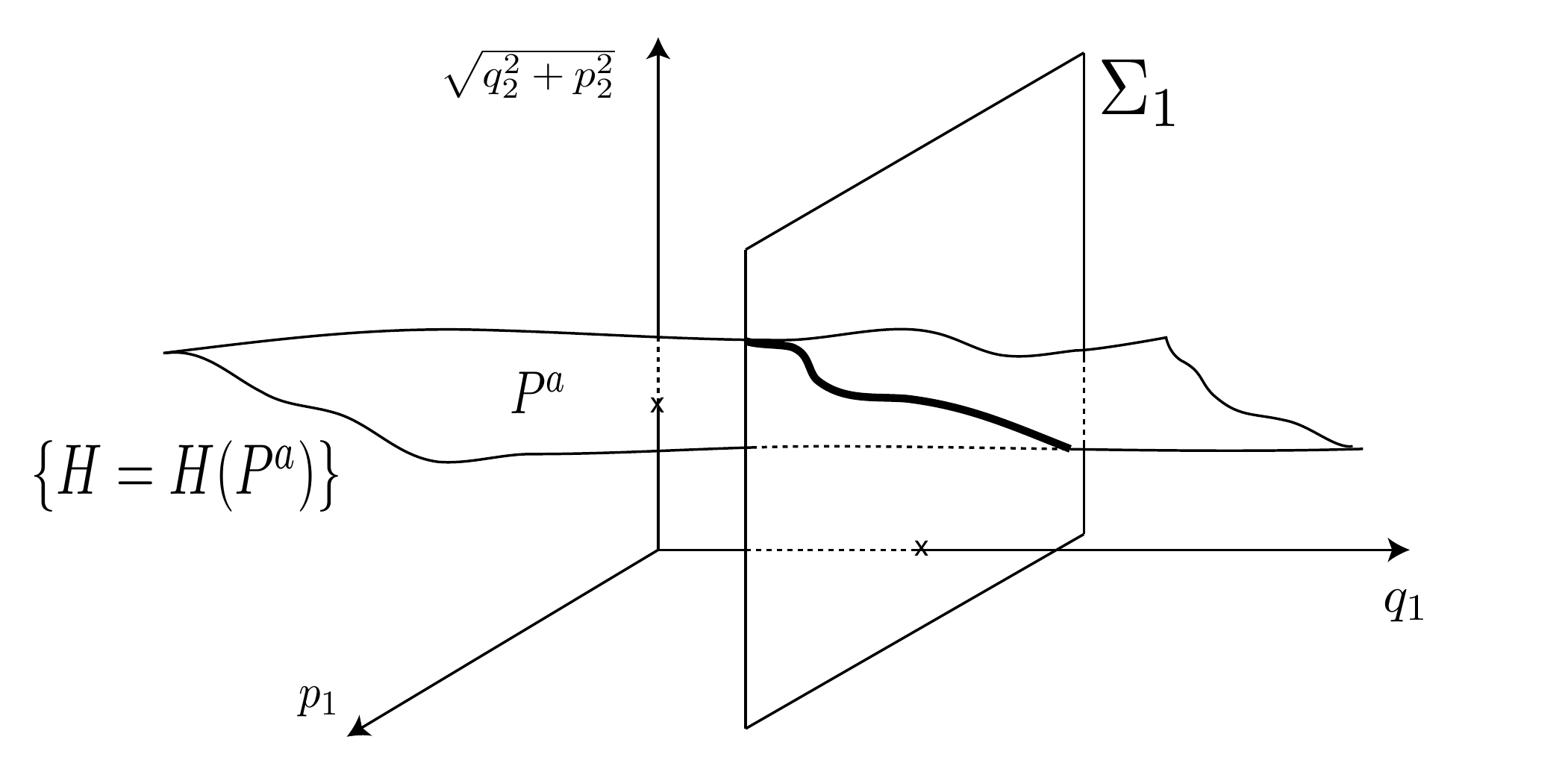}
\caption{Intersection of the energy level set $\{\H=\H(\P^\a)\}$ with the section $\Sigl$.}
\label{DESSIN5}\end{minipage}
\begin{minipage}{5mm}
\hfill
\end{minipage}
\begin{minipage}[b]{75mm}
\includegraphics[scale=0.3]{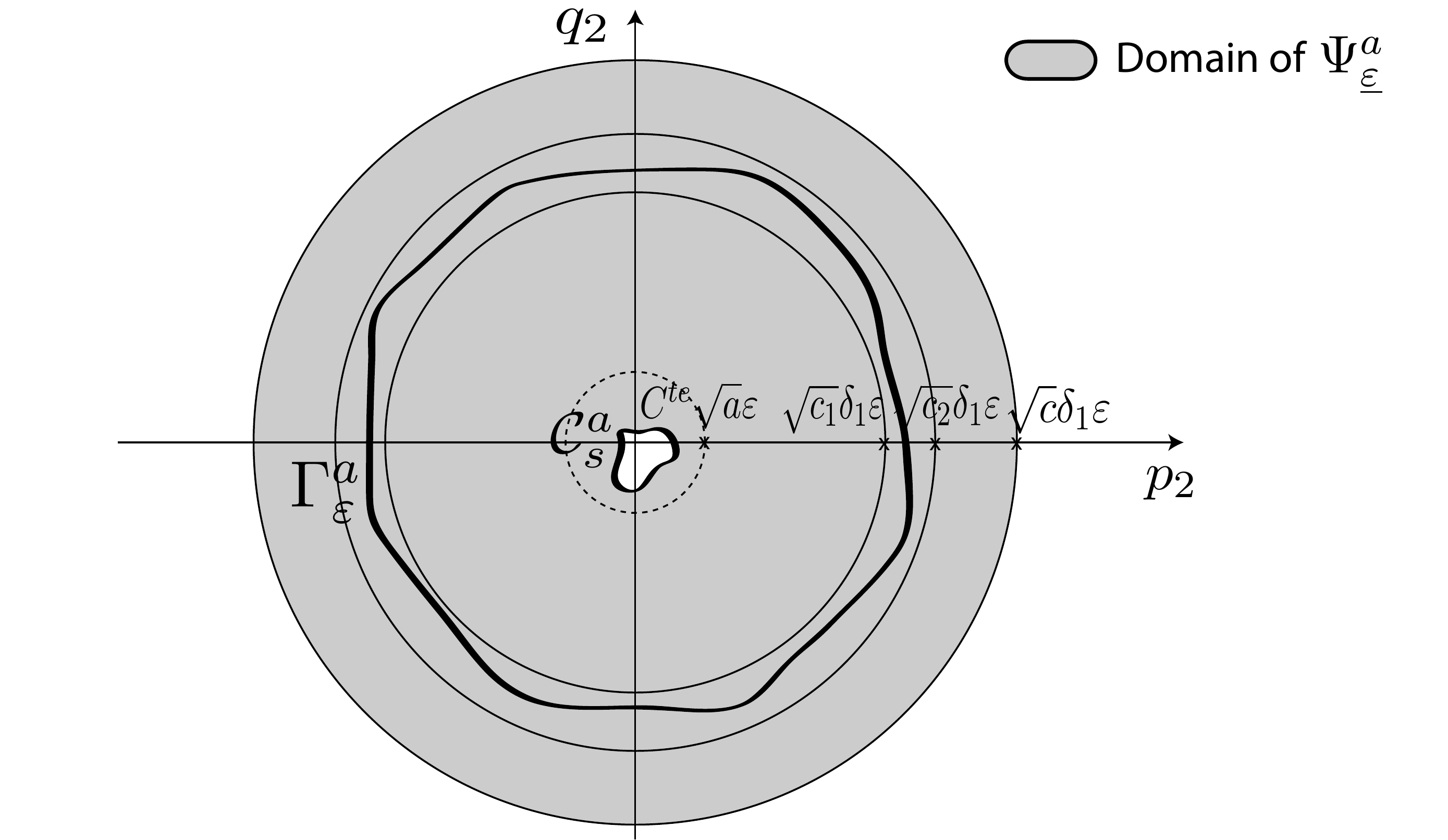}
\caption{Invariant curve $\Gamma^\a_\eps$ and size of the annulus used in the proof.}
\label{DESSIN14}

\end{minipage}

\end{figure}

Let us denote by $\Reta$ the restriction of $\Ret$ to the energy level set of $\P^\a$, $\ie$ to
$$\Sigl\cap\{\H=\H(\P^\a)\}=\{\qp/\ql=\dll, \H(\qp,\ep)=\H(\P^\a_{\ep},\ep)\}$$
(see Figure \ref{DESSIN5}). Proposition \ref{PropReta} below 
states that $\Reta$ can be considered as a diffeomorphism of a disc of $\R^2$.

\begin{prop}\label{PropReta}
%
%
%
Let us define the curve
$$\C^\a_s:=\Sigl\cap\CS\cap\{\H=\H(\P^\a)\}.$$
There exists $c,\co$ such that, for $\a\leq\co\dll^2\eps^2$, the restriction $\Ret^\a$ of $\Ret$ to the energy level set $\{\H=\H(\P^\a)\}$ reads as a diffeomorphism
$$\Ret^\a_{\ep} : \{(\qz,\pz)/ \qz^2+\pz^2\leq c\dll^2\eps^2, (\qz,\pz) \text{ outside of }\C^\a_s\} \to \{(\qz,\pz), \qz^2+\pz^2\leq \omo\dll^2\eps^2\}.$$
\end{prop}

The proof of this Proposition is in Part \ref{sub(iii)}.

\subsubsection*{Existence of an invariant curve}

We prove now the existence of invariant curves for each diffeomorphism $\Ret^\a_{\ep}$. This will be possible by an appropriate choice of some parameters : in this part, we chose and fix the order $n$ of the normal form and the power $\No$ of $\eps$ in the Hamiltonian expressed with the three parameters $\ep=(\eps,\nuu,\muu)$ (see (\ref{3parametres}) above), and also a value $\muu$ as a power of $\eps$. Here is the part of the proof where the normal form is the most fully used. Precisely, we show in this subsection the following

\begin{prop}\label{propGamma}
There exist $\lo$ and $\epso$, $\cl$, $\cz$ such that for $\eps\leq\epso, \a\leq\eps^{3\lo+1}, \muu\leq\eps^{4\lo+1}$ and $n\geq2(\lo+2), \No\geq4\lo+5$, the map $\Ret^\a_{\eps}:=\Ret^\a_{(\eps,\eps^2,\eps^{4\lo+1})}$ has an invariant curve $\Gamma^\a_\eps$ in the annulus of $\R^2$
$\{(\qz,\pz)/\Iz\in[\cl\dll^2\eps^2,\cz\dll^2\eps^2]\}$.
\end{prop}

\textbf{Proof.} The rest of this subsection is devoted to the proof of this proposition. For that purpose, we use the following KAM theorem stated by Moser \cite{Moser} :

\begin{thm}[KAM theorem]\label{KAMth}
Let 
$$\begin{array}{rcl}
\Phi : \R/2\pi\Z\times[a,b] & \rightarrow & \R/2\pi\Z \times\R \\
(\q,\ro) &\mapsto & (\q+\alpha(\ro)+F(\q,\ro),\ro+G(\q,\ro)).
\end{array}$$
We assume that the map $\Phi$ is exact and that there exists $m_0>0$ such that for all $\ro$, 
\begin{equation}\label{HypKAMalpha}
\frac{1}{m_0}\leq\frac{d\alpha}{d\ro}(\ro)\leq m_0.
\end{equation} 

Then there exist $\lo\in\NN$ and $\dlo(m_0)>0$ such that if
$$\NormeC{0}{F}+\NormeC{0}{G}\leq\dlo(m_0) \quad \text{and} \quad \NormeC{\lo}{\alpha}+\NormeC{\lo}{F}+\NormeC{\lo}{G}\leq m_0$$
are verified, then $\Phi$ admits an invariant curve of the form
\begin{equation} \label{FormeGamma}
\left\{(\q,\ro)=(\q'+f(\q'),\ro_0+g(\q')), \q'\in\R/2\pi\Z \right\}, 
\end{equation}
where $f$, $g$ are $\C^1$ functions.
\end{thm} 

The following Proposition \ref{prophypKAM} will involve (proof below) that, with an appropriate choice of the parameters, the maps $\Retap_{\ep}$ (which are the maps $\Reta_{\ep}$ after a change of coordinates, given in subsection \ref{SecEstimKAM}) satisfy the hypothesis of the KAM theorem \cite{Moser} applying it with
\begin{eqnarray}
&&\Retap_{(\eps,\eps^2,0)} : (\q,\ro)\mapsto (\q+\alpha^\a_\eps(\ro),\ro), \nonumber\\
&&\Retap_{(\eps,\eps^2,\muu)} : (\q,\ro) \mapsto (\q+\alpha^\a_\eps(\ro)+F^\a_{\ep}(\q,\ro),\ro+G^\a_{\ep}\q,\ro)). \nonumber
\end{eqnarray} 

%
\begin{prop} \label{prophypKAM} (see Figure \ref{DESSIN14})
There exist $\co,\cl,\cz$ and $m_0>0$ such that for $\eps$ sufficiently small, for all $\a\in[0,\co\dll^2\eps^2]$ and all $k\leq\entiere(\frac{\No-1}{4})-1$, in the annulus 
$$\{(\qz,\pz)/\Iz\in[\cl\dll^2\eps^2,\cz\dll^2\eps^2]\},$$
the $\Retap_{\ep}$ satisfy
\begin{enumerate}
\item the map $\Retap_{(\eps,\eps^2,\muu)}$ is exact;
\item $-m_0\leq\Frac{\partial\alpha^\a_\eps}{\partial\ro}\leq-\Frac{1}{m_0}$;
\item $\NormeC{k}{\alpha^\a_\eps}\leq m_k$;
\item $\NormeC{k}{F^\a_{\ep}}+\NormeC{k}{G^\a_{\ep}}\leq m_0 \cdot\left(\Frac{\a^2}{\eps^{6k}}+\Frac{\muu}{\eps^{4k}}\right).$
\end{enumerate}
\end{prop} 

\textbf{Proof :} in section \ref{SectionKAM}.

Recall that after the normalization until degree $n$, in Proposition \ref{PropNF}, we had $\eps^{4\n-8}$ in the expression of $\H$, that we rewrote $\eps^{4\n-8}=\muu\nuu\eps^{\No}$ when we introduced the parameters $\muu$ and $\nuu$. We recall also that we have already chosen the value of $\nuu$, $\nuu=\eps^2$, at the beginning of subsection \ref{SubStrategieKAM}. We now chose the values of $\n$ and $\No$ : we apply the normalization Proposition \ref{PropNF} with $\n:= 2\lo+3$ and chose $\No:=4\lo+1$.
Then necessarily $\muu=\eps^{4\lo+1}$, and for $\a\leq\eps^{3\lo}$ Proposition \ref{prophypKAM} implies that the hypotheses of the KAM theorem are satisfied in the annulus $\Iz\in[\cl\dll^2\eps^2,\cz\dll^2\eps^2]$ for all the maps $\Retap_\eps$. Then there exists a curve $\Gamma^\a_\eps$ of the form (\ref{FormeGamma}) in the annulus $\Iz\in[\cl\dll^2\eps^2,\cz\dll^2\eps^2]$, invariant by the map $\Reta_{\eps}=\Reta_{(\eps,\eps^2,\eps^{4\lo+1})}$.
\cqfd

\subsection{Proof of Theorem \ref{TH}} \label{StrategieConclusion}

This subsection is devoted to the proof of Theorem \ref{TH}, using the notations and results contained in the previous parts.
Let us fix $\a$ and $\eps$, and consider the unstable manifold $\WIPa$ of the periodic orbit $\P^\a$. 

\begin{figure}[!h]

\begin{minipage}[b]{85mm}
\includegraphics[scale=0.35]{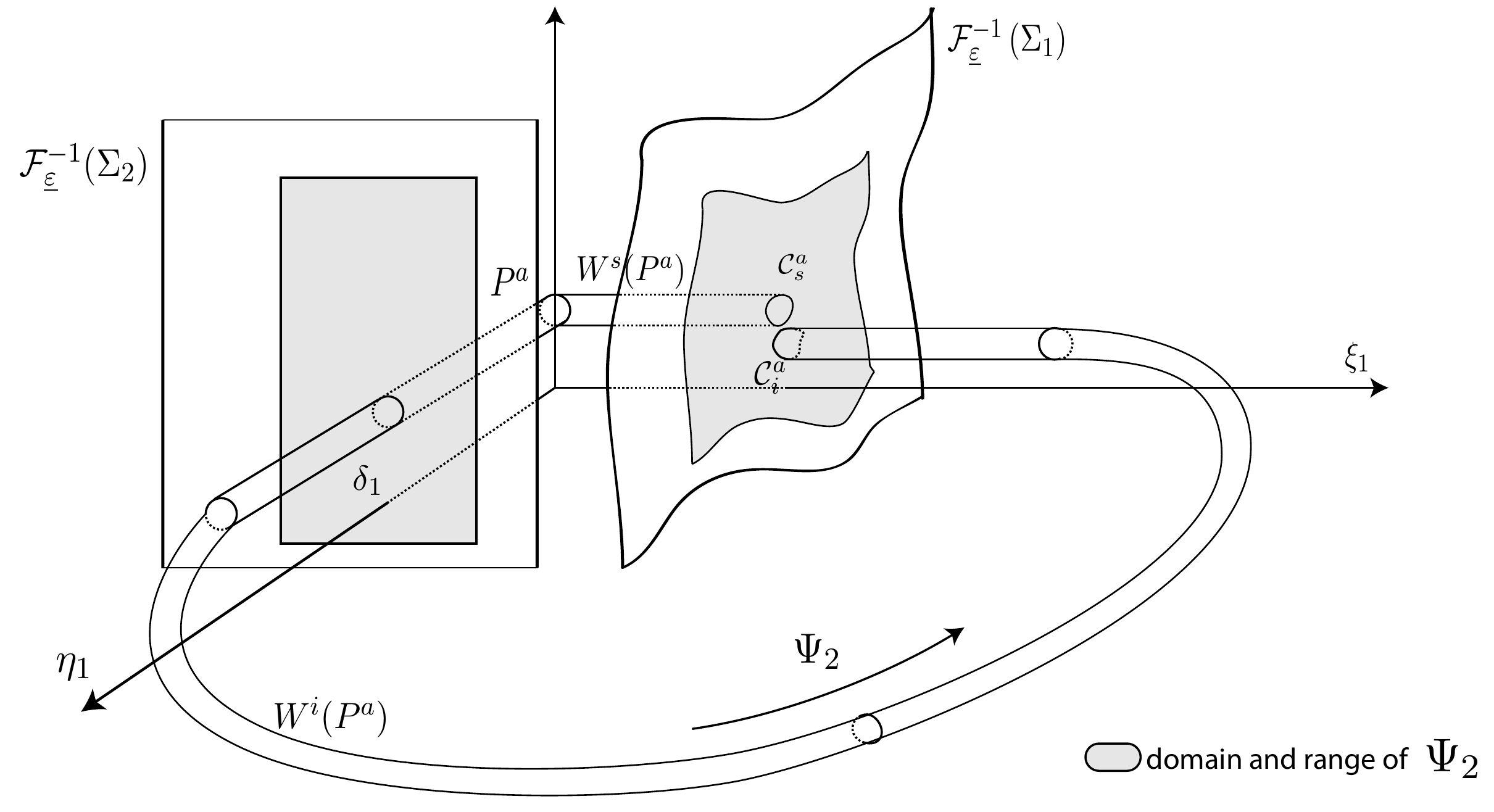}
\caption{$\WIPa$ and $\Sigl$ intersect.}
\label{DESSIN15}\end{minipage}
\begin{minipage}{5mm}
\hfill
\end{minipage}
\begin{minipage}[b]{55mm}
\includegraphics[scale=0.35]{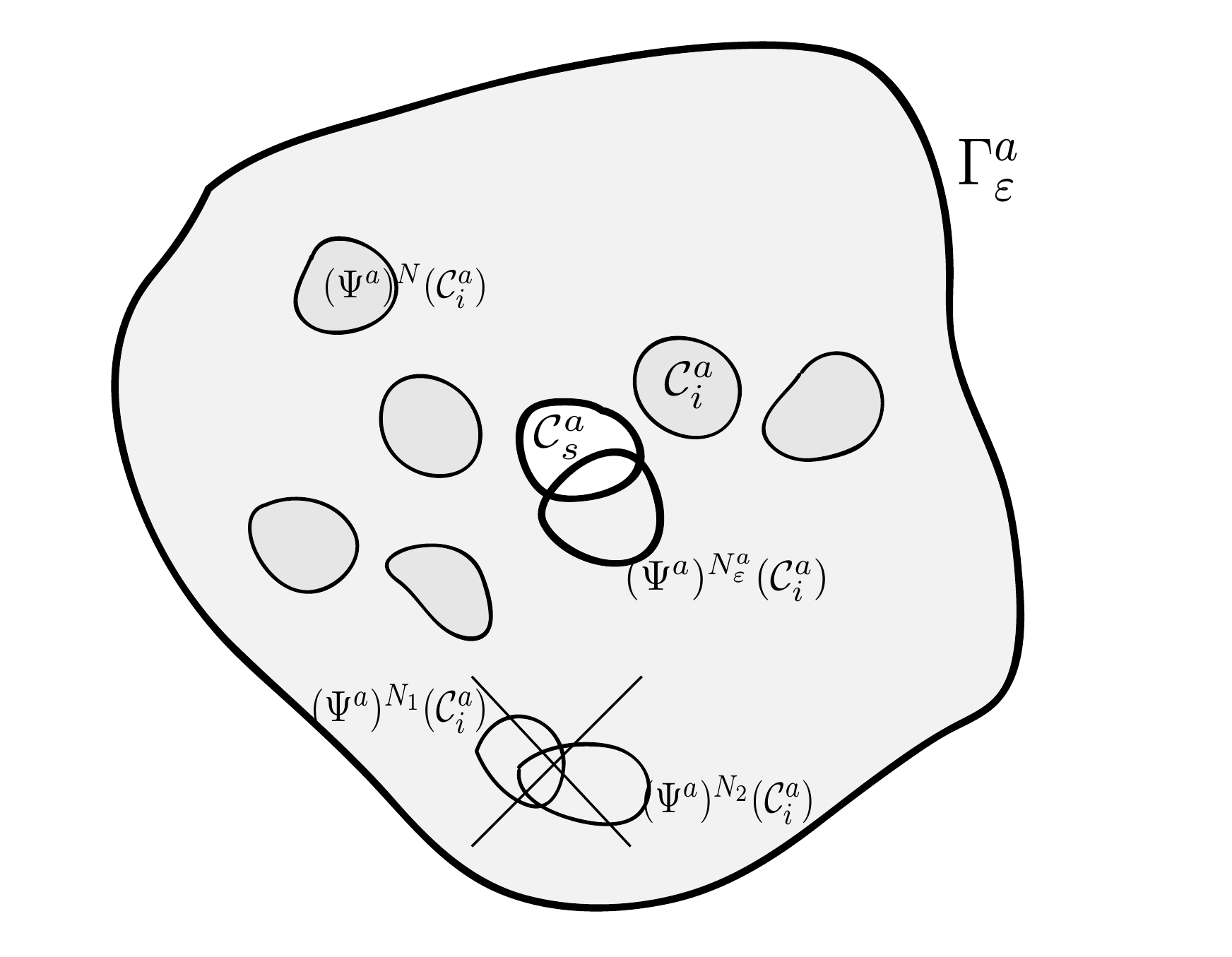}
\caption{Iterations of the map $\Reta$ on the curve $\C^\a_u$.}
\label{DESSIN18}

\end{minipage}

\end{figure}

Our first aim is to verify that $\WIPa$ does intersect $\Sigl$, and moreover that this intersection is in the neighborhood of $(\dll,0,0,0)$ in which $\Sigl \cap \{\H=\H(\P^\a)\}$ is a graph and also inside the invariant curve $\Gam^\a_\eps$. More precisely we prove below that for $\a$ and $\eps$ sufficiently small, $\WIPa$ hits the set
\begin{equation}\label{EnsDefReta}
\Sigl\cap\B(0,\dll)\cap\{\xiet/\xiz^2+\etz^2\leq\frac{1}{4}\cl\dll^2\eps^2\} ;
\end{equation}
To prove this, we use that in the coordinates $\xiet$, the unstable manifold of $\P^\a$ is the tube
$$\{\xiet\in\R^4/ \xil=0,\xiz^2+\etz^2=\a\},$$
whose intersection with the hyperplane $\{\xiet/\etl=\dll\}=\FFep^{-1}(\Sigz)$ is the circle $\{(0,\dll,\xiz,\etz)/\xiz^2+\etz^2=\a\}$. Then for $\a\leq\dll^2$ we can use the map $\Ret_2$ of Proposition \ref{PropRet2}, which maps $\Sigz$ onto $\Sigl$ : this way, we get that $\WIPa$ intersects $\Sigl\cap\B(0,\dll)$ (see Figure \ref{DESSIN15}). Using moreover the estimate (\ref{MajRet2}) of Proposition \ref{PropRet2}, we get that for $\a<\cl\dll^2\eps^2$ and for $\eps$ sufficiently small, $\WIPa$ intersects $\Sigl$ in the set (\ref{EnsDefReta}).

\vspace{1ex}

Denoting by $\C^\a_u$ \label{NotCai} the curve representing in $(\qz,\pz)$-coordinates the intersection $\WIPa\cap\Sigl$, we then have two possible alternatives :
\begin{itemize}
\item first option : $\C^\a_u$ and the curve $\C^\a_s$ of the stable manifold intersect {\bf in this case, the proof of the existence of an homoclinic connection to $\P^\a$ is completed},
\item second option : $\C^\a_u$ and $\C^\a_s$ do not intersect. 
\end{itemize} 
{\bf Let us consider the second option}. We know that both curves have the same symplectic area $\a$, so $\C^\a_u$ cannot be entirely contained inside $\C^\a_s$. Necessarily, $\C^\a_u$ is then outside of $\C^\a_s$. 


And we proved that $\Reta_\eps$ is defined on the set (\ref{EnsDefReta}) outside of the curve $\C^\a_s$, so we can consider $\Reta_\eps(\C^\a_i)$. Moreover, we know that for $\a$ and $\eps$ sufficiently small, the interior of the curve $\Gamma^\a_\eps$ of Proposition \ref{propGamma} above maps onto the interior  of $\Gamma^\a_\eps$ with the map $\Reta_\eps$. Then $\Reta_\eps(\C^\a_i)$ is also a curve whose symplectic area is $\a$ and contained inside $\Gamma^\a_\eps$. We again have two possible options :
\begin{itemize}
\item first option : $\Reta_\eps(\C^\a_i)$ and $\C^\a_s$ intersect, {\bf in this case, the proof of the existence of an homoclinic connection with two loops to $\P^\a$ is completed},
\item second option : $\Reta_\eps(\C^\a_i)$ is outside $\C^\a_s$, and then belongs to the domain of $\Reta_\eps$. 
\end{itemize}

\vspace{1ex}

Then we iterate this process. Since $\C^\a_u$ is defined as the first intersection of $\WIPa$ and $\Sigl$, for all $N$ $(\Reta_\eps)^N(\C^\a_i)$ and $\C^\a_u$ can not intersect. As $\Reta_\eps$ is a diffeomorphism (and then is invertible), then for all $N_1, N_2$, $(\Reta_\eps)^{N_1}(\C^\a_i)$ and $(\Reta_\eps)^{N_2}(\C^\a_i)$ do not intersect. Iterating this process $N_0$ times, the set of the $(\Reta_\eps)^N(\C^\a_i)$ for $N\leq N_0$ cover a surface whose area is $N_0\cdot\a$. And this surface is inside the curve $\Gam^\a_\eps$, whose area is finite. Then, the process must stop for one $N^\a_\eps$, $\ie$ necessarily there exists one $N=N^\a_\eps\in\NN$ for which $(\Reta_\eps)^{N}(\C^\a_i)$ and $\C^\a_s$ intersect (see Figure \ref{DESSIN18}). And this means that {\bf there exists an homoclinic connection with $N^\a_\eps$ loops to the periodic orbit $\P^\a_\eps$}.


\cqfd

\begin{rem}[about the importance of Hypothesis $(H6)$.]\label{RemH6}

\vspace{1ex}

The Hypothesis $(H6)$ is explicitly used below in the proof of Lemma \ref{LemGammaCroissante}, and we use this Lemma to prove the existence of $\Ret^\a_{\ep}$ \textbf{outside} of the curve $\C^\a_s$ in Proposition \ref{PropReta} above (if hypothesis $(H6)$ was not verified, the domain would be inside the curve). And the latter is crucial to allow the iterations of the first return map when the curves do not intersect (see Section \ref{StrategieConclusion}).

\vspace{1ex}

In the heuristic picture of the strategy outlined above in Section \ref{subsubStrategie}, the Hypothesis $(H6)$ is what allows to state that if the homoclinic trajectory deflects, it deviates \textbf{inwards} ("interior" in the $(\ql,\pl)$ coordinates).
\end{rem}


\section{Normal form and scaling : proof of Proposition \ref{PropNF}}
\label{PartieFormeNormale}
This section is devoted to the proof of Proposition \ref{PropNF}. We proceed in three main steps, constructing the three maps of $(i)$, $(ii)$ and $(iii)$ in the proposition.

\vspace{1ex}

\textbf{Step 1. Consequences of the Normal Form Theorem \ref{ThmNF}.} Under Hypotheses $(H1),\cdots,(H6)$, it is possible to find appropriate coordinates in $\R^4$ such that $\xbf$, $\Hlb$ and $\Om$ read 
$$\xbf=(\qbfun,\pbfun,\qbfd,\pbfd), \quad\Hlb(\xbf):=\frac{1}{2}\pbfun^2+\frac{1}{2}\omo(\qbfd^2+\pbfd^2)+\Ocal(\lb|\xbf|^2+|\xbf|^3), $$
and $\Om(\xbf,\ybf)=\scal{J\xbf}{\ybf}, $ where 
{\small$$  J= \left (\begin{array}{cccc}
           0 & 1 &0 & 0\\
           -1 & 0 &0 & 0\\
           0 & 0 &0 & 1\\
           0 &   0 &-1 & 0          
           \end{array}    
           \right).$$}
The matrix $J$ is not under the standard form $J_2$ (see the statement of Theorem \ref{ThmNF}), but this form can be obtained up to an isometric linear change of coordinates 
$$X=TX',\ T^*=T^{-1}, \quad
T=\left(\begin{array}{cccc}
1&0&0&0\\
0&0&1&0\\
0&1&0&0\\
0&0&0&1
\end{array}\right).$$

Hence, Theorem \ref{ThmNF} is still true in this case. So, applying it to the family of Hamiltonians $\Hlb(\xbf)$ we get the existence of a canonical transformation $\xbf=\phi_{n,\lb}(\xtilde)$ such that 
$$\Htildeulb(\xtilde) = \Hlb(\phi_{n,\lb}(\xtilde))=\frac{1}{2}\ptildeun^2+\frac{1}{2}\omo(\qtilded^2+\ptilded^2)+\Ncaltildeunlb(\xtilde)+ \Rcaltildeunlb(\xtilde),$$
where the rest $\Rcaltildeunlb(\xtilde)$ is a $\C^1$-one parameter family of analytic Hamiltonians satisfying $\Rcaltildeunlb(\xtilde)={\calO}\left(|\xtilde|^{n+1}\right)$ and where $\Ncaltildeunlb$ is a real polynomial of degree $\leq n$ satisfying $\Ncaltildeunlb(\xtilde)={\cal O}(\lb |\xtilde|^2+|\xtilde|^3)$ and
\begin{equation}\label{EqNsfmu}
\Ncaltildeunlb(\qtildeun,\ptildeun+t\qtildeun, \rot_{\omo t}(\qtilded,\ptilded))=\Ncaltildeunlb(\qtildeun,\ptildeun, \qtilded, \ptilded)\ \qquad \mbox{for all } t\in\R , 
\end{equation}
with
$$  \rot_{\omo t}= 
\left (\begin{array}{cr}
       \cos \omo t & -\sin \omo t\\
       \sin \omo t & \cos \omo t       
       \end{array}\right).$$

Setting $t=\frac{2\pi}{\omo} \ell$, $\ell\in \Z$ in (\ref{EqNsfmu}) and pushing $\ell$ to $\infty$, we get that necessarily $\Ncaltildeunlb$ does not depend on $\ptildeun$, $\ie$ 
$$\Ncaltildeunlb(\qtildeun,\ptildeun, \qtilded,\ptilded)=\Ncaltildeunlb^o(\qtildeun,\qtilded, \ptilded).$$
Then identifying $\R^2$ and $\CC$ via $(\qtilded,\ptilded)\mapsto(z_2=\qtilded+\I \ptilded,\overline{z}_2 =\qtilded-\I \ptilded)$ and defining 
$\Mcaltildeunlb^o(\qtildeun, z_2, \overline{z}_2 ) =\Ncaltildeunlb^o(\qtildeun,\qtilded,\ptilded ),$ identity (\ref{EqNsfmu}) reads
$$\Mcaltildeunlb^o(\qtildeun, \E^{\I\omo t}z_2, \E^{-\I\omo t}\overline{z}_2)=\Mcaltildeunlb^o(\qtildeun, z_2, \overline{z}_2) \qquad \mbox{pour tout } t\in\R.$$
Then,  setting $t=-\frac{{\rm arg} z_2}{\omo}$ and $t=-\frac{{\rm arg}z_2}{\omo}+\pi$, we obtain
$$\Mcaltildeunlb^o(\qtildeun, z_2, \overline{z}_2)=\Mcaltildeunlb^o(\qtildeun,|z_2|, |z_2|)=\Mcaltildeunlb^o(\qtildeun, -|z_2|, -|z_2|),$$
which ensures that $\Mcaltildeunlb^o(\qtildeun, z_2, \overline{z}_2 )=\Ptildeunlb(\qtildeun,|z_2|^2 )$ where $(\qtildeun,I_2)\mapsto \Ptildeunlb(\qtildeun, I_2 )$ is a real polynomial. Hence
$$ \Ncaltildeunlb(\qtildeun,\ptildeun, \qtilded,\ptilded )=\Ncaltildeunlb^o(\qtildeun,\qtilded, \ptilded )=\Ptildeunlb(\qtildeun, \qtilded^2+\ptilded^2 ). $$

\vspace{2ex}

{\it Hence we have proved the existence of a $\C^1$ one parameter family of canonical analytic transformations $\xbf=\phi_{n,\lb}(\xtilde )$ such that close to the origin $\Htildeulb(\xtilde ) = \Hlb(\phi_{n,\lb}(\xtilde))$ reads 
$$\Htildeulb(\xtilde ) = \frac{1}{2}\ptildeun^2+\frac{1}{2}\omo(\ptilded^2+\qtilded^2)+\Ptildeunlb(\qtildeun,\ptilded^2+\qtilded^2 ) + \Rcaltildeunlb(\xtilde),$$ 
where the rest $\Rcaltildeunlb(\xtilde)$ is a $\C^1$ one parameter family of analytic Hamiltonians satisfying $\Rcaltildeunlb(\xtilde)={\calO}\left(|\xtilde|^{n+1}\right)$ and where $\Ptildeunlb$ is a real polynomial with respect to $(\qtildeun,\ptilded^2+\qtilded^2)$ of degree less than $n$ with respect to $(\qtildeun,\ptilded,\qtilded)$, whose coefficients are $\C^1$ functions of $\lb$. Moreover, $\Ptildeunlb$ satisfies 
$$\Ptildeunlb(\qtildeun, \ptilded^2+\qtilded^2 )={\cal O}\Bigl(\lb(\qtildeun^2+\ptilded^2+\qtilded^2)+(|\qtildeun|+|\qtilded|+|\ptilded|)^3\Bigr).$$}

\vspace{1ex}

{\bf Step 2 : Change of parameter $\lb$.} Expanding $\Ptildeunlb$ we get that $\Htildeulb(\xtilde )$ reads
$$\begin{array}{ll}
      \Htildeulb(\xtilde ) =& \frac{1}{2}\ptildeun^2+\frac{1}{2}\underline{\widetilde{\om}}(\lb)(\qtilded^2+\ptilded^2) -\frac{1}{2}\uctilde_1(\lb)\qtildeun^2+\uctilde_2(\lb)\qtildeun^3+\uctilde_3(\lb)\qtildeun (\qtilded^2+\ptilded^2) \lba
      & +\Qtildeunlb(\qtildeun,\ptilded^2+\qtilded^2 ) + \Rcaltildeunlb(\xtilde),\end{array}$$
\nopagebreak
where
$$\Qtildeunlb(\qtildeun, \ptilded^2+\qtilded^2)={\cal O}\Bigl((|\qtildeun|+|\qtilded|+|\ptilded|)^4\Bigr),$$ 
and where $\uctilde_1$, $\uctilde_2,\uctilde_3$ and $\underline{\widetilde{\om}}$ are $\C^1$ functions of $\lb$ satisfying  $\underline{\widetilde{\om}}(0)=\omo$, $\uctilde_1(\lb)=c_{10}\lb+o(\lb)$,  $\uctilde_2(0)=c_{20}$. Moreover,  \textbf{because of hypotheses (H2), (H3) et (H4)} and  $\omo,\  c_{10}$ and $c_{20}$ are different from $0$.

\vspace{1ex}

Finally, since we only consider in this paper the "half bifurcation" corresponding to $c_{10}\lb>0$  (hypotheses (H3) and (H5)), the Implicit Function Theorem ensures that the identity  
$$\epstilde^2=  \uctilde_1(\lb)=c_{10}\lb+o(\lb)$$
can be inverted in a neighborhood of the origin, $\ie$ $\lb=\thetatilde(\epstilde^2)$ where $\thetatilde$ is a function of class $\C^1$.

\vspace{2ex}

{\bf Step 3 : scaling.} We perform a scaling in space and time suggested by the normal form of order 3.  Indeed, for any $n$, the normal form part of the Hamiltonian admits an homoclinic connection to $0$ which depends on $\epstilde$. Moreover, for $n=3$ this homoclinic connection  $\widetilde{h}_{\epstilde}$  can computed explicitly and has the form
$$\widetilde{h}_{\epstilde} (\ttilde)=(\epstilde^2\  \ql^h(\epstilde\ \ttilde),\epstilde^3\  \pl^h(\epstilde\  \ttilde), 0,0).$$   
To study the dynamics close to this homoclinic connection, it is more convenient to rescale the system so that the rescaled normal form of order 3 of the rescaled Hamiltonian admits an homoclinic connection which does not depend on $\epstilde$. So we perform the following scaling in space and time
\begin{equation}\label{Eqscalingtilde}
\begin{array}{lll}
\qtildeun =\TFrac{1}{2\sqrt{2} \ctilde_2(\epstilde)}\ \epstilde^2 \qlp, & \qquad\qtilded =\TFrac{1}{2\sqrt{2} \ctilde_2(\epstilde)}\ \epstilde^{\frac{5}{2}} \ \qzp, &\qquad t =\epstilde \ \widetilde{t}, \lba
\ptildeun =\TFrac{1}{2\sqrt{2} \ctilde_2(\epstilde)}\ \epstilde^3 \plp, & \qquad\ptilded =\TFrac{1}{2\sqrt{2} \ctilde_2(\epstilde)}\ \epstilde^{\frac{5}{2}} \ \pzp.& \qquad 
\end{array} 
\end{equation}
This scaling is a conformal mapping which is well defined for $\epstilde$ small since $\ctilde_2(\epstilde)=c_{20} +{\cal O}(\epstilde^2)$ and since by hypothesis (H4), $c_{20}\neq 0$ holds. Note that the change of coordinates on $(\qtildeun,\ptildeun,\qtilded,\ptilded)$ is not canonical. Nevertheless, together with the scaling in time, for $\epstilde\neq 0$ the rescaled differential system is an Hamiltonian system whose Hamiltonian reads 
$$\frac{8\ctilde_2(\epstilde)}{\epstilde^6}\Htilde_{\epstilde}(\TFrac{1}{2\sqrt{2} \ctilde_2(\epstilde)}\epstilde^2 \qlp,\TFrac{1}{2\sqrt{2} \ctilde_2(\epstilde)} \epstilde^3 \plp,\TFrac{1}{2\sqrt{2} \ctilde_2(\epstilde)} \epstilde^{\frac{5}{2}}\qzp,\TFrac{1}{2\sqrt{2} \ctilde_2(\epstilde)} \epstilde^{\frac{5}{2}}\pzp).$$ 

Moreover, to work with regular functions of the parameter, and  because of the square root in the scaling we also perform a last change of parameter
$$\eps^2=\epstilde.$$
For $\eps\in]-\eps_0, 0[\  \cup\  ]0,\eps_0[$ with $\eps_0=\sqrt{\epstilde_0}$, we get for $\xp=(\qlp,\pzp,\qzp,\pzp)$ the rescaled Hamiltonian
\begin{eqnarray} 
\Hprime_{\eps} (\xp) &=&\frac{8\ctilde_2(\epstilde)}{\epstilde^{12}} \Htilde_{\eps^2}(\sigma_\eps(\xp))\nonumber\\
&&\hspace{-4ex}=\Frac{1}{2}\plp^2-\Frac{1}{2}\qlp^2+2\sqrt{2}\qlp^3 +\Frac{\omega(\eps)}{2\eps^2}({\qzp}^2+{\pzp}^2)+\eps^2 \Nprime_{n,\eps}(\qlp,\qzp^2+\pzp^2)+\eps^{4n-8}\Rprime_{n,\eps}(\xp),\nonumber
\end{eqnarray}
where $\Nprime_{n,\eps}$ is a polynomial of degree less than $n$ with respect to $(\qlp,\qzp,\pzp)$ whose coefficients are $\C^1$ functions of $\eps$ and which satisfies
$$\begin{array}{ll}   
\Nprime_{n,\eps}(\qlp,\qzp^2+\pzp^2) &:=\frac{1}{\eps^2 \eps^{12}}\left[ \ctilde_3(\eps^2)\eps^{14} \qlp(\qzp^2+\pzp^2) +\Qtilde_{n,\eps^2}\Bigl(\eps^4 \qlp,\eps^{10}(\qzp^2+\pzp^2)\Bigr)\right],\lba
      &={\cal O}\left( |\qlp||\qzp^2+\pzp^2|+\eps^2(|\qlp|^2+|\qzp^2+\pzp^2|)^2\right ).
\end{array}$$
Recall finally that  $\Rtilde_{n,\epstilde}$ is a $\C^1$ one parameter family of analytic Hamiltonians satisfying $\Rtilde_{n,\epstilde}(\xtilde)={\cal O}\left(|\xtilde|^{n+1}\right)$ for all $\epstilde$ in $]-\epstilde_0,\epstilde_0[$. Thus the explicit formula
$$\Rprime_{n,\eps}(x): = \frac{1}{\eps^{12} \eps^{4n-8}}\Rtilde_{n,\eps^2}(\TFrac{1}{2\sqrt{2} \ctilde_2(\eps^2)}\eps^4 \qlp,\TFrac{1}{2\sqrt{2} \ctilde_2(\eps^2)} \eps^6 \plp,\TFrac{1}{2\sqrt{2} \ctilde_2(\eps^2)} \eps^5 \qzp,\TFrac{1}{2\sqrt{2} \ctilde_2(\eps^2)} \epstilde^5 \pzp)$$
ensures that $\Rprime_{n,\eps}(x)$ is a $\C^1$ family of the parameter $\eps$.
\cqfd

\section{Existence of the first return map on $\Sigl$ : proof of Propositions \ref{PropRet2} and \ref{PropDefRetp}}\label{Poincare}

This section is devoted to the proof of Propositions \ref{PropRet2} (given in subsections \ref{SubPropRet2}) and \ref{PropDefRetp} (given in Part \ref{SubPropDefRet}). The previous subsections \ref{SubLemRegPhiprime} is devoted to the proof of a lemma used to prove these propositions.

\subsection{Smoothness of the flow apart from the rotation}\label{SubLemRegPhiprime}

The following lemma will be the {\it key to prove the smoothness of the first return map}, and is {\it a consequence of the Normal Form Theorem applied up to degree $n$}. This lemma is first used in a weak way (the $\C^1$ smoothness is sufficient) in the proof of Proposition \ref{PropRet2} below.

\begin{lem} \label{LemFlotp}
Denote $$\Rot_\theta:=\left(\begin{matrix}
    								I & O\\
           					O & \rot_\theta
    					\end{matrix}\right)
    					:=\left(\begin{matrix}
    								1& 0 &  0 & 0\\
           					0 &1&0&0	 \\    
           					0&0&\cos\theta&-\sin\theta\\
           					0&0&\sin\theta&\cos\theta
    					\end{matrix}\right),$$
and let $\flot(t,\qp,\eps,\nuu,\muu)$ be the flow of the Hamiltonian
\begin{eqnarray}
\H(\qp,\eps,\nuu,\muu)&=&-\ql\pl+\frac{\om(\eps)}{2\eps^2}(\qz^2+\pz^2)+\frac{1}{2}(\ql+\pl)^3\nonumber\\
  &&+\nuu\Q(\ql+\pl,\qz^2+\pz^2,\eps)+\nuu\muu\eps^{\No}\Reste(\qp,\eps).\nonumber
\end{eqnarray}
Then
$$\flot(t,\qp,\eps,\nuu,\muu)=\Rotep\flotp(t,\qp,\eps,\nuu,\muu),$$
where $\flotp$ belongs to 
$$\C^1\left(]-\epso,\epso[^3,\C^{\ko}(\R^5)\right), \quad \text{where } \ko:=\entiere\left(\frac{\No-1}{4}\right),$$
meaning in particular that $\flotp$ is $\C^1$ with respect to $(t,\qp,\eps,\nuu,\muu)$ and $\C^{\ko}$ with respect to $(t,\qp,\nuu,\muu)$.   					
\end{lem}

\textbf{Proof.} Denoting 
$$\flotp(t,\qp,\eps,\nuu,\muu):=\Rotepm\flot(t,\qp,\eps,\nuu,\muu),$$
given that $(\ql,\pl)$ and $(\qz^2+\pz^2)$ are preserved by the rotation $\Rotep$, $\flotp$ is then the flow of the nonautonomous Hamiltonian
\begin{eqnarray}
\Hp(t,\qp,\eps,\nuu,\muu)&=&-\ql\pl+\frac{1}{2}(\ql+\pl)^3+\nuu\Q(\ql+\pl,\qz^2+\pz^2,\eps)\nonumber\\
  &&+\nuu\muu\eps^{\No}\Reste((\ql,\pl,\rotep(\qz,\pz),\eps).\nonumber
\end{eqnarray}
Thus $\flotp$ has the smoothness of $\nabla_X\Hp$, and in the definition of $\Hp$ all is $\C^{\infty}$ in terms of $(t,\qp,\nuu,\muu)$ and $\C^1$ with respect to $(t,\qp,\eps,\nuu,\muu)$, except $\rotep(\qz,\pz)$ (when $\eps=0$). But the derivatives of $\rotep$ read 
\begin{eqnarray}
&&D_{(\qz,\pz)}(\rotep(\qz,\pz))=\rotep(\qz,\pz),\quad \frac{\partial}{\partial t}(\rotep(\qz,\pz))=\frac{\om(\eps)}{2\eps^2}\Om\rotep(\qz,\pz),\nonumber\\
&&\frac{\partial}{\partial \eps}(\rotep(\qz,\pz))=\Frac{2\eps^2\om'(\eps)+4\eps\om(\eps)}{4\eps^4}\Om\rotep(\qz,\pz).\nonumber
\end{eqnarray}
We then get the estimate 
$$\left|\nabla_X\Reste(\ql,\pl,\rotep(\qz,\pz),\eps)\right|_{\C^k}=\grandO{\eps\rightarrow 0}\left(\frac{1}{\eps^{4k}}\right).$$
So $\nabla_X \big(\nuu\muu\eps^{\No}\Reste((\ql,\pl,\rotep(\qz,\pz),\eps)\big)$ is a $\C^k$-function in the neighborhood of $\eps=0$ as soon as $4k+1\leq\No$.
\cqfd

\subsection{Existence of the global map from $\Sigz$ to $\Sigl$ : proof of Proposition \ref{PropRet2}}\label{SubPropRet2}

Let us introduce
\begin{equation}
\Sigl^-:=\left\{\qp\in\R^4/ \ql=\dll-\frac{1}{2}\dll\right\}, \quad \Sigl^+:=\left\{\qp\in\R^4/ \ql=\dll+\frac{1}{2}\dll\right\}.\nonumber
\end{equation}
We proceed in several steps.

%
%
%
%

\vspace{1ex}

\textbf{Step 1 : case $\ep=(\eps,0,0)$. } Let us prove the existence of some $T^{\pm}(\xil,\dll)$ such that for all $\xiet$ in
\begin{equation}
\FFo^{-1}(\Sigz)\cap\left\{\xiet / 0\leq\xil\leq \frac{1}{16}\dll,\sqrt{\Iz}\leq\frac{1}{2}\dll \right\},\nonumber
\end{equation}
$\flot(T^{\pm}(\xil,\dll),\FFo(\xil,\dll,\xiz,\etz),(\eps,0,0))$ belongs to $\Sigl^{\pm}$.

To prove this, we first recall that when $\ep=(\eps,0,0)$ the flow and $\FFo$ are uncoupled and that $\flot_{\ql,\pl}(.,.,\eps,0,0)$ and $\FFo$ do not depend of $\eps$ (see Part \ref{subsubFN3} for the flow and Proposition \ref{propFF} for $\FFo$). From these facts, we get first that if the $T^{\pm}$ exist, they are independent of $(\xiz,\etz)$ and $\eps$. Let us prove their existence working with the \textbf{restriction of the flow to the $(\ql,\pl)$-plane}.

We then use the phase portrait drawn in Part \ref{subsubFN3} : let us study \textbf{which of the orbits hit $\Sigl^-$ and $\Sigl^+$}. For that purpose, we work first in the $(\qlp,\plp)$ coordinates of Part \ref{subsubFN3} and use the parameter $\alpha$ of \eqref{IntroAlpha}. For $\dll$ sufficiently small, if an orbit hits $\Sigl^-$, then it necessarily hits also $\Sigl^+$. Let us denote by $\alpha_\dll$ the parameter of the orbit passing through the point $(\qlp,\plp)=(\frac{\sqrt{2}}{2}\dll,0)$, which reads also $(\ql,\pl)=(\frac{1}{2}\dll,\frac{1}{2}\dll)$. Then the orbits labelled by any $\alpha\in[\alpha_\dll,0]$ hits $\Sigl^-$. A short computation gives
$$\alpha_\dll=-\frac{1}{2}\dll^2(1-4\dll).$$

It remains now to \textbf{get back to $(\ql,\pl)$-coordinates} and then to the local coordinates $(\xil,\etl)$. Firstly, from the study of the phase portraits, we get that the condition $\alpha\leq0$ means in local coordinates that $\xil\geq0$. Let us then study the condition $\alpha\geq\alpha_\dll$, by studying the orbit labelled by $\alpha_\dll$, in the neighborhood of $\Sigz$.

On $\Sigz$, $\etl=\dll$ is satisfied ; let us work in a domain a little larger in $(\ql,\pl)$ coordinates
\begin{equation}\label{domain} 
\{(\ql,\pl),\ql\geq0, \pl\in[\frac{1}{2}\dll,\frac{3}{2}\dll]\},
\end{equation}
and look for a condition on $\ql$ which ensures that $(\ql,\pl)$ belongs to an orbit satisfying $\alpha\geq\alpha_\dll$. We consider the orbit $\alpha=\alpha_\dll$ more precisely on the half part satisfying $\plp\geq0$, thus we obtain
$$\qlp\in[\sqrt{2}\frac{1}{2}\dll,\sqrt{2}\frac{3}{2}\dll].$$
From the equation $\plp=\sqrt{\qlp^2-4\sqrt{2}\qlp^3+\alpha_\dll}$ of the orbit labelled by $\alpha_\dll$, we compute a lower bound of $\ql$ on this orbit. Let us denote by $\ql^{\alpha_\dll}(\qlp)$ a graph description of the orbit ; for $\qlp$ in $[\sqrt{2}\frac{1}{2}\dll,\sqrt{2}\frac{3}{2}\dll]$, we get
\begin{equation}
\ql^{\alpha_\dll}(\qlp)=\frac{1}{\sqrt{2}}(\qlp-\sqrt{\qlp^2-4\sqrt{2}\qlp^3+\alpha_\dll}) = \frac{1}{\sqrt{2}} \frac{4\sqrt{2}\qlp^3-\alpha_\dll}{\qlp+\sqrt{\qlp^2-4\sqrt{2}\qlp^3+\alpha_\dll}}\geq\frac{1}{12}\dll.\nonumber
\end{equation}
Thus we obtain that if $\ql\in[0,\frac{1}{12}\dll]$ and $\pl\in [\frac{1}{2}\dll,\frac{3}{2}\dll]$, then $(\ql,\pl)$ belongs to an orbit labelled by $\alpha\geq\alpha_\dll$.

\textbf{Finally}, with the aid of $(vi)$ of Proposition \ref{propFFcomplete}, and given that $\Iz$ is preserved by the flow, we get that for $\dll$ sufficiently small, all the points $$\qp\in\Sigz\cap\FFo\left(\{\xiet\in\R^4,0\leq\xil\leq\frac{1}{16}\dll,\sqrt{\xiz^2+\etz^2}\leq\frac{1}{2}\dll\}\right)$$
belong to orbits hitting $\Sigl^-$ and $\Sigl^+$, which achieves the proof of existence of $T^{\pm}(\xil,\dll)$ as claimed above.

\vspace{2ex}

\textbf{Step 2 : upper and lower bounds for $T^{\pm}(\xil,\dll)$} (these bounds are useful to get back to the general case $\nuu,\muu\neq0$, see Step 3). Step 1 ensures the existence of the $T^{\pm}(\xil,\dll)$. $T^-$ and $T^+$ are also locally defined by the Implicit Equations
$$\flot_{\ql}(T^{\pm},\FFo(\xil,\dll,0,0),(\eps,0,0))=\dll\pm\frac{1}{2}\dll.$$
As in Step 2, we get the equivalent Implicit Equation 
\begin{equation} \label{impT}
\flotp_{\ql}(T^{\pm},\FFo(\xil,\dll,0,0),(\eps,0,0))=\dll\pm\frac{1}{2}\dll.
\end{equation}
Let us prove that (\ref{impT}) satisfies the hypotheses of the $\C^{\ko}$ implicit function theorem in the neighborhood of each $(T^{\pm}({\xil}_0,\dll),{\xil}_0)$. On one hand, the result of Lemma \ref{LemFlotp} ensures that this equation is $\C^{\ko}$. On the other hand,
\begin{eqnarray}
\partial_T\flotp_{\ql}(T,\qp,(\eps,0,0))&=&\partial_{\pl}\H(\flotp(T,\qp,(\eps,0,0)),(\eps,0,0)),\nonumber
\end{eqnarray}
where
\begin{eqnarray}
\partial_{\pl}\H(\qp,(\eps,0,0))&=&-\ql+\frac{3}{2}(\ql+\pl)^2\nonumber
\end{eqnarray}    
and 
\begin{equation}\nonumber
\flotp_{\ql}(T^{\pm}({\xil}_0,\FFo({\xil}_0,\dll,0,0),(\eps,0,0)))=\dll\pm\frac{1}{2}\dll \quad \left|\flotp_{\pl}(T^{\pm}({\xil}_0,\FFo({\xil}_0,\dll,0,0),(\eps,0,0)))\right|\leq\dll
\end{equation}
hold because of the definition $T^{\pm}$. Then the Implicit Functions Theorem applies for $\dll$ sufficiently small. This ensures that $T^{\pm}$ are continuous with respect to $\xil$, and thus bounded.

\vspace{2ex}

\textbf{Step 3 : existence of $T_L(X,\ep)$.} From Step 1, we know that
\begin{eqnarray}
 \flot_{\ql}(T^-(\xil,\dll),\FFo(\xil,\dll,\xiz,\etz),(\eps,0,0))&=&\frac{1}{2}\dll,\\
 \flot_{\ql}(T^+(\xil,\dll),\FFo(\xil,\dll,\xiz,\etz),(\eps,0,0))&=&\frac{3}{2}\dll.
\end{eqnarray}
Recall that $(\ql,\pl)$ are preserved by the $\Rot_\theta$, so that $\flot_{\ql}=\flotp_{\ql}$. And Step 2 ensures that there exists $T^-(\dll)$ and $T^+(\dll)$ such that for all $\xil$
$$T^-(\dll)\leq T^-(\xil,\dll),T^+(\xil,\dll)\leq T^+(\dll).$$
Then, from the $\C^1$-smoothness of $\flotp_{\ql}$ (Lemma \ref{LemFlotp} for $N_0\geq3$), together with $(iv)$ of Proposition \ref{propFFcomplete}, we obtain that for $\nuu$ and $\muu$ sufficiently small,
\begin{eqnarray}
 \flot_{\ql}(T^-(\xil,\dll),\FFo(\xil,\dll,\xiz,\etz),\ep)&\in&[\frac{1}{4}\dll,\frac{3}{4}\dll],\\
 \flot_{\ql}(T^+(\xil,\dll),\FFo(\xil,\dll,\xiz,\etz),\ep)&\in&[\frac{5}{4}\dll,\frac{7}{4}\dll].
\end{eqnarray}
Thus we get the existence of $T_L(X,\ep)$ thank to the Intermediate Value Theorem.

\vspace{2ex} 
 
\textbf{Step 4 : uniqueness of $T_L$.} For that purpose, it is sufficient to show that in the set
$$\left\{\qp / \frac{1}{4}\dll\leq\ql\leq\frac{7}{4}\dll,|\pl|\leq\dll,\sqrt{\Iz}\leq\dll\right\},$$
the flow satisfies $\frac{d\ql}{dt}<0$. And indeed
\begin{eqnarray}
\Frac{d\ql}{dt}&=&-\ql+\frac{3}{2}(\ql+\pl)^2+\nuu\partial_{\pl}(\Q(\ql+\pl,\qz^2+\pz^2,\eps)+\muu\eps^{\No}\Reste(\qp,\eps))\nonumber\\
     &\leq& -\frac{1}{32}\dll+\frac{3}{2}9\dll^2+M \nuu. \nonumber
\end{eqnarray}
Then, for $\nuu$ and $\dll$ sufficiently small, $\frac{d\ql}{dt}<0$, which ensures the uniqueness of $T_L$.

\vspace{1ex}

\textbf{Step 5 : upper bound (\ref{MajRet2}).} We use that the trajectories of the flow satisfy
\begin{eqnarray}
\left|\Frac{d (\qz^2+\pz^2)}{dt}\right|&=&\left|\muu\nuu\eps^{\No} (\pz\partial_{\pz}-\qz\partial_{\qz})\Reste(\qp,\eps)\right|\nonumber\\
  	&\leq& M \muu\nuu\eps^{\No}  \quad \text{ for } \qp\in\B(0,\roo).\nonumber
\end{eqnarray}
Given that $T_L(\FFep^{-1}\qp,\ep)\leq T^+(\dll)$, we get the upper bound (\ref{MajRet2}) claimed above.

\cqfd

\subsection{Existence and smoothness of the entire first return map : proof of Proposition \ref{PropDefRetp}}\label{SubPropDefRet}

The following Proposition is a more detailed version of Proposition \ref{PropDefRetp}.
\begin{prop} \label{PropDefRet}
There exists a first return map
\begin{eqnarray}
&\hspace{-10ex}\Ret:\Sigl\cap\FFep\left(\{\xiet/0<\etl\leq\frac{1}{24}\dll,\sqrt{\xiz^2+\etz^2}\leq\dll\}\right)\times]-\epso,\epso[^3 \rightarrow \Sigl\nonumber\\
&\hspace{22ex}(\pl,\qz,\pz,\ep)\mapsto\Rot_{\frac{\om(\eps)}{2\eps^2}T(\pl,\qz,\pz,\ep)}\flotp(T(\pl,\qz,\pz,\ep),(\dll,\pl,\qz,\pz),\ep),\nonumber
\end{eqnarray}
where $T$ belongs to $\C^1\left(]-\epso,\epso[^3,\C^{\ko}(\R^3)\right)$ with $\ko$ defined in Lemma \ref{LemFlotp}. 

\vspace{1ex}

Moreover, denoting $\Ret=(\Ret_{\ql},\cdots,\Ret_{\pz})$, there exists $M$ such that for $\ep$ and $\dll$ sufficiently small, 
\begin{equation}\label{MajRet}
\left|\Ret_{\qz}((\pl,\qz,\pz),\ep)^2+\Ret_{\pz}((\pl,\qz,\pz),\ep)^2-(\qz^2+\pz^2)\right|\leq\nuu M \dll^3
\end{equation}
holds on the domain of $\Ret$.
\end{prop}

\textbf{Proof of the existence.} As mentioned in part \ref{subsubRet}, observing the local phase portrait in coordinates $\xiet$ in $\B(0,\roop)$, the local map from $\Sigl$ to $\Sigz$ is very simple, following the level sets 
$$\{\xiet/\xil\etl=C^{te},\xiz^2+\etz^2=C^{te}\}.$$
This local map exists on the set $\{\xiet/0<\etl\leq\dlz\}$. We want to compose this local map with the global map $\Ret_2$ of Proposition \ref{PropRet2} whose domain is
$$\Sigz\cap\FFep\left(\{\xiet/0\leq\xil\leq\frac{1}{16}\dll, \sqrt{\xiz^2+\etz^2}\leq\dll\}\right).$$
Then we need to trim the domain of the local map so that its range is included in the domain of $\Ret_2$. Since $\xiz^2+\etz^2$ and $\xil\etl$ are conserved by the flow, it is sufficient to restrict the local map to trajectories for which $\xiz^2+\etz^2\leq\dll^2$ and $|\xil\etl|\leq\dlz\frac{1}{16}\dll$. We proved in Lemma \ref{LemGrapheSigl} that if $|\etl|\leq2\dll$ and $\sqrt{\qz^2+\pz^2}\leq2\dll$ then 
$$\xil=\gsig(\etl,\xiz,\etz)\leq\frac{3}{2}\dll.$$ 
So it is sufficient to trim the domain of the local map to the set
$$0<\etl\leq\dlz \text{ and }|\etl|\leq\dll, \sqrt{\qz^2+\pz^2}\leq\dll \text{ and } \etl\leq\frac{2}{3}\frac{1}{16}\dll=\frac{1}{24}\dll,$$
and get then the domain of $\Ret$ stated in the Proposition. 

\vspace{1ex}

\textbf{Proof of the smoothness.} We showed above the existence of $T(\pl,\qz,\pz,\ep)$. Recall that $T$ is locally defined by the implicit equation given by the intersection with $\Sigl$ :
\begin{equation}\label{EqImpT}
\flot_{\ql}(T,(\dll,\pl,\qz,\pz),\ep)=\dll.\nonumber
\end{equation}
Recall that $\flot_{\ql}=\flotp_{\ql}$, so that we get the equivalent implicit equation 
\begin{equation} \label{impT1}
\flotp_{\ql}(T,(\dll,\pl,\qz,\pz),\ep)=\dll.
\end{equation}
Let us prove that (\ref{impT1}) satisfies the hypotheses of the $\C^{\ko}$ implicit function theorem in the neighborhood of each $(T({\pl}_0,{\qz}_0,{\pz}_0,{\ep}_0),{\pl}_0,{\qz}_0,{\pz}_0,{\ep}_0)$. If they are fulfilled, we can construct a $\C^{\ko}$ map $T^*(\pl,\qz,\pz,\ep)$, and the uniqueness of the first return ensures that $T=T^*$ and finally that $T$ is $\C^{\ko}$.
The result of Lemma \ref{LemFlotp} ensures that this equation is $\C^{\ko}$. Let us prove that all $(\qp,\ep)$ in 
$$\Sigl\cap\FFep\left(\{\xiet/0<\etl\leq\frac{1}{24}\dll,\sqrt{\xiz^2+\etz^2}\leq\dll\}\right)\times]-\epso,\epso[^3,$$
satisfies $\partial_T\flotp_{\ql}(T,\qp,\ep)\neq 0$. We have
\begin{eqnarray}
\partial_T\flotp_{\ql}(T,\qp,\ep)&=&\partial_{\pl}\H(\flotp(T,\qp,\ep),\ep),\nonumber
\end{eqnarray}
where
\begin{eqnarray}
\partial_{\pl}\H(\qp,\ep)&=&-\ql+\frac{3}{2}(\ql+\pl)^2+\nuu\partial_{\pl}\Q(\ql+\pl,\qz^2+\pz^2,\eps)\nonumber\\
    &&+\muu\nuu\eps^{\No}\partial_{\pl}\Reste(\qp,\eps)).\nonumber
\end{eqnarray}    
Given that we work in $\B(0,\roo)$ (see the truncature in Part \ref{subsubchange}), $\partial_{\pl}\Q$ and $\partial_{\pl}\Reste$ are uniformly bounded. Moreover, 
\begin{equation}\nonumber
\flotp_{\ql}(T(\pl,\qz,\pz,\ep),(\dll,\pl,\qz,\pz),\ep)=\dll, \quad \quad \left|\flotp_{\pl}(T(\pl,\qz,\pz,\ep),(\dll,\pl,\qz,\pz),\ep)\right|\leq\dll
\end{equation}
hold because of the definition $T$ and the range of $\Ret_2$ (Proposition \ref{PropRet2}). Then
\begin{equation}\nonumber
\partial_T\flotp_{\ql}(T({\pl}_0,{\qz}_0,{\pz}_0,{\ep}_0),{\pl}_0,{\qz}_0,{\pz}_0,{\ep}_0)\leq-\dll+\frac{3}{2}(2\dll)^2+M\nuu<0
\end{equation}
holds for $\dll<\frac{1}{6}$ and $\nuu$ sufficiently small.

\vspace{1ex}

\textbf{Proof of the estimate (\ref{MajRet}).} Recall how $\Ret$ was constructed (proof of the existence above). We use two previous results :
on one hand, in local coordinates the flow preserves $\xiz^2+\etz^2$, and on the other hand the estimate (\ref{MajRet2}) gives an upper bound of the variation of $\qz^2+\pz^2$ by the map $\Ret_2$.

To complete the proof, we moreover compute estimates of the difference between $\xiz^2+\etz^2$ and $\qz^2+\pz^2$ on $\Sigl$ and $\Sigz$ by the changes of coordinates $\FFep$ and $\FFep^{-1}$. We detail the proof for $\FFep$ on $\Sigz$ ; we more precisely need estimates for $\xiet$ in the domain of $\Ret_2$. Recall that we denote $\FFep:=(\philep,\psilep,\phizep,\psizep)$ ; $(viii)$ and $(ix)$ of Proposition \ref{propFFcomplete} allows to get that, for all $\xiet\in\B(0,\roop)$,
\begin{eqnarray}
&&\hspace{-2ex}\left|\phizep\xiet^2+\psizep\xiet^2-(\xiz^2+\etz^2)\right|\nonumber\\
 &&\leq \Mo\nuu|\xiet|^2(2\Mo\nuu|\xiet|^2+|\xiz|+|\etz|).\nonumber
\end{eqnarray}
For $\xiet$ in the domain of $\Ret_2$, we then obtain
\begin{equation}\label{VarIz}
\left|\phizep\xiet^2+\psizep\xiet^2-(\xiz^2+\etz^2)\right|\leq \nuu\Mo'\dll^3.
\end{equation}
The same strategy for the change of coordinates $\FFep^{-1}$ in the domain of $\Ret$ allows to achieve the proof of (\ref{MajRet}).
\cqfd

\section{Restrictions of the first return map seen as diffeomorphisms of an annulus of $\R^2$ : proof of Proposition \ref{PropReta}}
\label{PartieReta}
This part is devoted to the proof of Proposition \ref{PropReta}. The proof itself is in Part \ref{sub(iii)}. The previous parts are devoted to the proof of lemmas useful for the proof : in Parts \ref{sub(i)} and \ref{SubLemGrapheSigl} we prove respectively that the center-stable manifold $\CS$ and $\Sigl$ locally read as a graphs. In Part \ref{GeomSPa} we prove two lemmas concerning the geometry of the different energy level sets on $\CS$. And in Part \ref{sub(ii)} we prove that the level sets $\{H=\H(\P^\a_{\ep})\}$ locally read as graphs and give a description of its position with respect to the graph of $\CS$, thank to the lemmas of the previous part. Finally, with the aid of all these graphs we can prove Proposition \ref{PropReta}.  

From now on, we only need two of the three parameters $\ep=(\eps,\nuu,\muu)$ : we only study the influence of the remainder on the dynamic (we explained it in Part \ref{subsubchange} when we introduced the parameters $\ep, \nuu,\muu$). Thus we introduce the notation \label{Notnup}
$$\nup:=(\eps,\eps^2,\muu).$$

\subsection{In $\Sigl$, the center-stable manifold $\CS$ reads as a graph}\label{sub(i)}

\begin{lem}\label{Lemgcsep}
For $\dll$ sufficiently small, there exists an analytic function $\gcsep$ such that
$$\CSep\cap\Sigl\cap\B(0,\dll)=\left\{(\ql,\pl,\qz,\pz)/\ql=\dll, \pl=\gcsep(\qz,\pz)\right\}.$$
Moreover, in $\B(0,\dll)$, $\gcsep$ satisfies
\begin{equation}\label{boncote}
\FFep^{-1}(\dll,\pl,\qz,\pz)\in\{\etl>0\} \Leftrightarrow \pl>\gcsep(q_2,p_2).
\end{equation}
\end{lem}

\textbf{Proof.} We proceed in three steps.

\textbf{Step 1 : existence.} Recall that we denote $\FFep^{-1}:=(\invphilep,\invpsilep,\invphizep,\invpsizep)$. For $\dll\leq\roop$, for any $(\dll,p_1,q_2,p_2)$ in $\B(0,\dll)$ the equivalence
\begin{eqnarray}
(\dll,p_1,q_2,p_2)\in\CSep &\Longleftrightarrow& \invpsilep(\dll,p_1,q_2,p_2)=0, \label{CSimp}
\end{eqnarray}
holds. Then, statement $(vii)$ of Proposition \ref{propFFcomplete} ensures the existence of $\Mo$ independent of $\ep$ such that
\begin{equation}
|\invpsilep(q_1,p_1,q_2,p_2)-p_1|\leq\Mo|(q_1,p_1,q_2,p_2)|^2 \quad \text{for all }(q_1,p_1,q_2,p_2)\in\B(0,\roop). \nonumber
\end{equation}
So for $\dll$ sufficiently small (independently of $\ep$), for all $(q_2,p_2)\in\B(0,\dll)$, we get
$$\invpsilep(\dll,-\dll,q_2,p_2)<0<\invpsilep(\dll,\dll,q_2,p_2).$$
Then the Intermediate value Theorem ensures the existence of a $\gcsep(q_2,p_2)\in[-\dll,\dll]$ such that
$$\invpsilep(\dll,\gcsep(q_2,p_2),q_2,p_2)=0.$$

\vspace{1ex}

\textbf{Step 2 : uniqueness and smoothness.} We use the implicit equation (\ref{CSimp}). Differentiating $(iii)$ of Proposition \ref{propFFcomplete} with respect to $p_1$ we get the existence of a convergent power series $\M_1$ (convergent on a ball $\B(0,\roop)$) independent of $\ep$ such that 
\begin{equation}
\partial_{p_1}(\invpsilep(q_1,p_1,q_2,p_2)-p_1)\prec (q_1\hspace{-0,5ex}+\hspace{-0,5ex}p_1\hspace{-0,5ex}+\hspace{-0,5ex}q_2\hspace{-0,5ex}+\hspace{-0,5ex}p_2)\M_1(q_1\hspace{-0,5ex}+\hspace{-0,5ex}p_1\hspace{-0,5ex}+\hspace{-0,5ex}q_2\hspace{-0,5ex}+\hspace{-0,5ex}p_2).\nonumber
\end{equation}
See Appendix \ref{tech} for definitions and properties of the relation $\prec$.

Then, if $\dll\leq\roop$, for all $(q_1,p_1,q_2,p_2)\in \B(0,\dll)$,
$$\big|\partial_{p_1}(\invpsilep(q_1,p_1,q_2,p_2)-p_1)\big|\leq 4\dll\M_1(4\dll).$$
So for $\dll$ sufficiently small $\partial_{p_1}\invpsilep>0$ in $\B(0,\dll)$. Thus, for any fixed $(q_2,p_2)$, the function $p_1\mapsto \invpsilep(\dll,p_1,q_2,p_2)$ is strictly increasing, and then we get the uniqueness of the $\gcsep(q_2,p_2)$ of Step 1. Moreover, the fact that $\partial_{p_1}\invpsilep$ is nonzero allows to apply the analytic Implicit Function Theorem to Equation (\ref{CSimp}) in the neighborhood of any fixed $(\dll,\gcsep(q_2,p_2),q_2,p_2)$ : we then obtain the analyticity of $\gcsep$.

\vspace{1ex}

\textbf{Step 3 : (\ref{boncote})} holds given that the function $\pl\mapsto \invpsilep(\dll,\pl,\ql,\pl)$ is increasing.
\cqfd

\subsection{$\Sigl$ as a graph in local coordinates $\xiet$}\label{SubLemGrapheSigl}

\begin{lem} \label{LemGrapheSigl}
For $\eps,\nuu,\muu<\epso$ and $\dll$ sufficiently small, there exists an analytic map $\gsigep$ defined on the domain $\B_{\R^3}(0,2\dll):=\{(\etl,\xiz,\etz)/|\etl|\leq2\dll,\sqrt{\xiz^2+\etz^2}\leq2\dll\}$ satisfying
$$\FFep^{-1}(\Sigl)\cap\B(0,2\dll)=\left\{\xiet/\xil=\gsigep(\etl,\xiz,\etz),(\etl,\xiz,\etz)\in\B_{\R^3}(0,2\dll)\right\}.$$
%
%
\end{lem} 



\vspace{1ex}

\textbf{Proof.} The proof is very similar to the proof of Lemma \ref{Lemgcsep}, so we only detail what is different.

\textbf{Step 1. Existence.} Recall that we denote $\FFep=(\philep,\psilep,\phizep,\psizep)$. As in the Step 1 of Lemma \ref{Lemgcsep}'s proof, thank to the result $(vi)$ of Proposition \ref{propFFcomplete}, we prove that for $\dll$, $\ep$ sufficiently small,
$$\philep(\frac{1}{2}\dll,\etl,\xiz,\etz)-\dll<0<\philep(\frac{3}{2}\dll,\etl,\xiz,\etz)-\dll$$
holds for any fixed $(\etl,\xiz,\etz)$ in $\B_{\R^3}(0,2\dll)$. We obtain the existence of $\gcsep(\etl,\xiz,\etz)$ in $]\frac{1}{2}\dll,\frac{3}{2}\dll[$ thank to the Intermediate Value Theorem.

\vspace{1ex}

\textbf{Step 2. Uniqueness and smoothness.} We proceed as in Step 2 of Lemma \ref{Lemgcsep}'s proof, showing with the aid of $(ii)$ of Proposition \ref{propFFcomplete} that for $\dll$ sufficiently small 
$$\partial_{\xil}\philep(\xil,\etl,\xiz,\etz)>0$$
holds for all $\xil$ in $]\frac{1}{2}\dll,\frac{3}{2}\dll[$ and all $(\etl,\xiz,\etz)$ in $\B_{\R^3}(0,2\dll)$.

\cqfd

\subsection{Positions of the $\SPa$ on the graph of $\CS$}\label{GeomSPa}

The following lemma ensures that the energy of $\P^\a$ is strictly increasing in term of $\a$.

\begin{lem} \label{LemGammaCroissante}
There exist $\a_0, \eps_1>0$ and a convergent power series $\M_1$ such that for all $\eps,\nuu,\muu<\eps_1$ and all $\a,\a'<\a_0$
\begin{enumerate}
\item $\big|\H(\P^{\a}_{\ep},\ep)-\frac{\om(\eps)}{\eps^2}\a\big|\leq \frac{\nuu}{\eps^2}\a^2\M_1(\a)$ holds
\item if $\a<\a'$ then $\H(\P^{\a}_{\ep},\ep)=\HLep(0,\a)<\H(\P^{\a'}_{\ep},\ep)=\HLep(0,\a').$
\end{enumerate}
\end{lem}

\textbf{Proof.} Recall that $\HLep$ was introduced in (\ref{DefHL1}), and $\P^\a_{\ep}$ was defined by (\ref{DefPaep}), so that for any fixed $(\xi_{2,\a},\eta_{2,\a})$ such that $\xi_{2,\a}^2+\eta_{2,\a}^2=\a$
$$\H(\P^\a_{\ep},\ep)=\H(\FFep(0,0,\xi_{2,\a},\eta_{2,\a}),\ep)=\HLep(0,\a).$$
Denoting
$$\FFep-\FFo:=\nuu(\philep',\psilep',\phizep',\psizep')$$
and using the particular form of $\FFo$ stated in Proposition \ref{propFFcomplete}, we get
\begin{equation}
\FFep(0,0,\xi_{2,\a},\eta_{2,\a})=(0,0,\xi_{2,\a},\eta_{2,\a})+\nuu(\philep',\psilep',\phizep',\psizep')(0,0,\xi_{2,\a},\eta_{2,\a}).
\end{equation}
Let us denote
$$\H(q_1,p_1,q_2,p_2,\ep)=-q_1p_1+\frac{\om(\eps)}{2\eps^2}(q_2^2+p_2^2)+\H_3(q_1,p_1,q_2,p_2,\ep).$$
$\H_3$ is of order 3 in term of $\ql$ and is $\C^1$ in term of $\ep=(\eps,\nuu,\muu)$. Thus $\H_3$ admits a convergent upper bound for $\prec$ uniformly in $\ep$. So do $\philep',\psilep',\phizep',\psizep'$ as stated in $(i)$ of Proposition \ref{propFFcomplete}. We finally obtain the existence of a convergent power series $\mathcal{M}_1$ independent of $\ep$ such that
\begin{eqnarray}
\H(\FFep(0,0,\xiz,\etz),\ep)-\frac{\om(\eps)}{2\eps^2}(\xiz^2+\etz^2)&\prec&\frac{\nuu}{\eps^2}\mathcal{M}_3(\xiz,\etz).
\end{eqnarray}
Moreover, thank to the particular form of $\HLep=\H(\FFep(\cdot),\ep)$ stated in (\ref{DefHL}), we get that $\mathcal{M}_3$ can be chosen of order 4 and as a power series of $\xiz^2+\etz^2$. Thus there exists a convergent power series $\M_1$ such that
\begin{equation}\label{premreslem}
\HL_{\ep}(0,\a)-\frac{\om(\eps)}{2\eps^2}\a\prec\frac{\nuu}{\eps^2}\a^2\mathcal{M}_1(\a),
\end{equation}
which achieves the proof of $(i)$ of the lemma. As a consequence of (\ref{premreslem}), there exist $\a_0$ and a real $M_1$ such that 
\begin{equation}\label{elliptic1}
\partial_\a\HLep(0,\a)_{|\a=0}=\frac{\om(\eps)}{2\eps^2}, \quad\quad \text{for all }  \a\leq\a_0, \quad \partial_\a^2\HLep(0,\a)_{|\a}\leq\frac{\nuu}{\eps^2}M_1.
\end{equation}
So, \textit{given that $\om(0)>0$ (Hypothesis (H6))} for $\nuu$ sufficiently small (this small size is independent of $\eps$) and $\a<\a'\leq\a_0$,
\begin{equation}\label{elliptic2}
\H(\P^{\a'}_{\ep},\ep)-\H(\P^{\a}_{\ep},\ep)=\HL_{\ep}(0,\a')-\HL_{\ep}(0,\a)>0.
\end{equation}
\cqfd
\begin{rem}\label{RemH6b}
Hypothesis $(H6)$ only appears in this lemma ! But in an essential way (to derive (\ref{elliptic2}) from (\ref{elliptic1})).
\end{rem}


%

We introduce in the following lemma the diffeomorphism $\g$ which is the restriction of $\FFep$ to $\Sigl\cap\CS$, seen as a map from $\R^2$ onto $\R^2$ thank to the graph form of $\CS$ stated in Lemma \ref{Lemgcsep}. The following result gives a hint on how the stable manifolds $\SPa$ of the $\P^\a$ intersects with $\Sigl$ (recall that they are in $\CS$). Observe that in the $(\qz,\pz)$ coordinates, these intersections are the images through $\g$ of the circles of area $\a$. The following Lemma describes some properties of the ranges of the circles through the map $\g$.

\begin{lem}\label{concentriques}
Let $\Ca:=\{(\xiz,\etz)/ \xiz^2+\etz^2=a\}$ circle of area $\a$ in $\R^2$, and let us introduce
$$\g: \begin{array}[t]{rl}         
			\B_{\R^2}(0,\frac{1}{2}\dll) &\rightarrow \R^2\\
			(\xiz,\etz) &\mapsto (\phizep,\psizep)(\gsigep(0,\xiz,\etz),0,\xiz,\etz),
			\end{array}$$
where we recall the notation $(\philep,\psilep,\phizep,\psilep)=\FFep$. Then for $\dll$ and $\nuu$ sufficiently small,
\begin{description}
\item[$(a)$] $\g(\Ca)$ is a Jordan curve,
\item[$(b)$] if $\a<\a'$ then $\g(\Ca)$ is inside $\g(\Cap)$.
\end{description}
\end{lem}

\textbf{Proof of $(a)$.} Given that $\Ca$ is a Jordan curve it is sufficient to show that $\g$ is an homeomorphism from $\B_{\R^2}(0,\frac{1}{2}\dll)$ to $\g(\B_{\R^2}(0,\frac{1}{2}\dll))$. For this purpose, let us prove that the map 
$$\g^-:(\qz,\pz)\mapsto (\invphizep,\invpsizep)(\dll,\gcsep(\qz,\pz),\qz,\pz)$$
is the inverse function of $\g$, where $(\invphilep,\invpsilep,\invphizep,\invpsizep)=\FFep^{-1}$. 

\vspace{1ex}

We first verify that $\g^-$ is well defined on $\g(\B_{\R^2}(0,\frac{1}{2}\dll))$. On one hand, Lemma \ref{Lemgcsep} and Proposition \ref{propFFcomplete} ensures that $\g^-$ is well defined if $\sqrt{\qz^2+\pz^2}\leq\dll$ and $(\dll,\gcsep(\qz,\pz),\qz,\pz)$ is in $\B(0,\roop)$, and that $|\gcsep(\qz,\pz)|\leq\dll$. On the other hand, thank to Lemma \ref{LemGrapheSigl} and $(vi)$ of Proposition \ref{propFFcomplete}, we get that for $\nuu$ sufficiently small $\g(\B_{\R^2}(0,\frac{1}{2}\dll))\subset\B(0,\dll)$. Then $\g^-\circ\g$ is well-defined for small values of $\nuu$ and $2\dll\leq\roop$.

\vspace{1ex}

Let us prove now that $\g^-\circ\g$ is equal to identity. From the definitions of $\gsigep$ and $\g$ and given that 
$$\FFep(\gsigep(0,\xiz,\etz),0,\xiz,\etz)\in\CS$$
(because $\etl=0$), we get that
$$(\dll,\gcsep(\g(\xiz,\etz)),\g(\xiz,\etz))=(\philep,\psilep,\phizep,\psizep)(\gsigep(0,\xiz,\etz),0,\xiz,\etz).$$
Thus
\begin{equation}
\g^-\circ\g(\xiz,\etz)=(\invphizep,\invpsizep)(\FFep(\gsigep(0,\xiz,\etz),0,\xiz,\etz))=(\xiz,\etz).\nonumber
\end{equation}
This achieves the proof of $(a)$.

\vspace{2ex}

\textbf{Proof of $(b)$.} Let us first show that $\g$ preserves the areas. Indeed, let $\C$ be a curve of $\B_{\R^2}(0,\frac{1}{2}\dll)$. We denote by $\mathcal{A}(\C)$ the symplectic area of $\C$ in $\R^2$ endowed with the restriction of $\Om$ to $\R^2$. Given that the curves of $\R^4$
$$\{\xiet/\xil=\gsigep(0,\qz,\pz),\etl=0,(\xiz,\etz)\in \g(\C)\}$$
$$\{\qp/\ql=\dll,\pl=\gcsep(\qz,\pz),(\qz,\pz)\in\C\},$$
are respectively subsets of $\{\xiet/\etl=0\}$ and of $\{\qp/\ql=\dll\}$, their areas are respectively $\mathcal{A}(\g(\C))$ and $\mathcal{A}(\C)$. So, given that $\FFep$ is symplectic and thus area-preserving, we get
\begin{eqnarray}
\mathcal{A}(\g(\C))&=&\mathcal{A}(\{\xiet/\xil=\gsigep(0,\qz,\pz),\etl=0,(\xiz,\etz)\in \g(\C)\})\nonumber\\
 &=&\mathcal{A}(\FFep(\{\xiet/\xil=\gsigep(0,\qz,\pz),\etl=0,(\xiz,\etz)\in \g(\C))\})\nonumber\\
  &=&\mathcal{A}(\{\qp/\ql=\dll,\pl=\gcsep(\qz,\pz),(\qz,\pz)\in\C\})=\mathcal{A}(\C).\nonumber
\end{eqnarray}
So $\g$ is area-preserving.

The result $(a)$ ensures that $\g(\C_{\a'})$ divides $\g(\B_{\R^2}(0,\frac{1}{2}\dll))$ into two connected subsets. Given that $\g$ is an homeomorphism (see the proof of $(a)$), $\g(\C_{\frac{1}{2}\dll})$ belongs to one of these subsets and $\g(\C_\a)$ to the other. Area preservation ensures that $\g(\C_{\frac{1}{2}\dll})$'s area is greater than $\g(\C_{\a'})$'s area, so $\g(\C_\a)$ is inside $\g(\C_{\a'})$ and $\g(\C_{\frac{1}{2}\dll})$ is outside.
\cqfd

\subsection{In $\Sigl$, the energy level set of $\P^\a$ reads as a graph ; position of this graph with respect to the graph of $\CS$}\label{sub(ii)}


\begin{lem}\label{LemGrapheHPa}
Let us define
 $$\hep(\qz^2+\pz^2,\a):=\frac{\om(\eps)}{2\eps^2}(\qz^2+\pz^2)-\H(\P^\a_{\ep},\ep).$$
For $\dll$ sufficiently small and $\a\leq\frac{1}{2}\dll^2$, for $\ep$ sufficiently small, there exists $\ghpaep$ such that
\begin{eqnarray}
&&\hspace{-6ex}\left\{\qp/\H(\qp,\ep)=\H(\P_{\ep}^\a,\ep)\right\}\cap\Sigl\cap\B(0,\dll) \nonumber\\
&&\hspace{6ex}\cap\big\{\qp/|\hep(\qz^2+\pz^2,\a)|\leq\dll^2\big\}\nonumber\\
&&=\left\{\qp/\ql=\dll,\pl=\ghpaep(\qz,\pz,\a), (\qz,\pz)\in\B_{\R^2}(0,\dll), |\hep|\leq\dll^2\right\}\nonumber
\end{eqnarray}
Moreover, $\ghpaep$ reads
\begin{equation}\label{Defghpaep} 
\ghpaep(\qz,\pz,\a):=\ghpah(\qz,\pz,\hep(\qz^2+\pz^2,\a),\ep),
\end{equation}
with $\ghpah$ analytic with respect to $(\qz,\pz,\h)$ and $\C^1$ with respect to $\eps$ ; more precisely $\ghpah$ belongs to $\C^1(]-\epso,\epso[^3,\An)$ (see Definition \ref{DefAn} of Part \ref{PartStatement}).
\end{lem}

\textbf{Proof.} We proceed in two main steps.

\textbf{Step 1. Existence.} Let us denote $\H_3:=\H-\H_2$ where $\H_2$ is the quadratic part of $\H$. We get
\begin{equation}
\H(\dll,\pl,\qz,\pz,\ep)=\H(\P^\a_{\ep},\ep)\Longleftrightarrow \hep(\qz^2+\pz^2,\a)=\dll\pl-\H_3(\dll,\pl,\qz,\pz,\ep).\label{defghpa}
\end{equation}
We first consider $h$ as an independent variable, $\ie$ we consider the equation
\begin{equation}\label{imph}
\h=\dll\pl-\H_3(\dll,\pl,\qz,\pz,\ep),
\end{equation}
From the explicit expression of $\H_3$ for $\nuu=0$ and given that $\H_3$ is $\C^1$ with respect to $\nuu$ and all its variables, for $\dll$ and $\nuu$ sufficiently small, we get that 
$$\dll\cdot2\dll-\H_3(\dll,2\dll,\qz,\pz,\eps,\nuu,\muu)\geq\dll^2,\quad \dll\cdot(-2\dll)-\H_3(\dll,-2\dll,\qz,\pz,\eps,0,\muu)\leq-\dll^2.$$
Thus the Intermediate Value Theorem ensures that the equation (\ref{imph}) admits a solution $\ghpah(\qz,\pz,\h,\ep)\in[-\dll,\dll]$ when $|\h|\leq\dll^2$.

\vspace{1ex}

\textbf{Step 2. Uniqueness and smoothness.} From the explicit form of $\H_3$ when $\nuu=0$ and given that $\H_3$ is $\C^1$ with respect to $\nuu$ and all its variables, we get that for $\dll$ and $\nuu$ sufficiently small, for $|\pl|\leq\dll$,
$$\partial_{\pl}(\dll\pl-\H_3(\dll,\pl,\qz,\pz,\ep))\geq\frac{1}{2}\dll>0.$$
Considering the implicit equation \eqref{imph} and proceeding like in the Step 2 of Lemma \ref{Lemgcsep}'s proof, we then get the uniqueness and the analyticity of $\ghpah$.
\cqfd

\vspace{2ex}

The following lemma is a refinement of Lemma \ref{LemGrapheHPa} : it gives an expression of the domain of $\ghpa$ in terms of $(\qz,\pz,\a)$ and a description of the position of the graph $\ghpa$ with respect to $\gcsep$.

\begin{lem}
There exists $\co$ such that for $\dll$ sufficiently small and for $\a\leq\co\dll^2\eps^2$, the function $\ghpa$ satisfies
\begin{eqnarray}
&&\hspace{-6ex}\{\H\hspace{-1ex}=\hspace{-0.5ex}\H(\P^\a)\}\cap\Sigl\cap\B(0,\dll)\cap\{\qp/\Iz\leq\frac{1}{2}\omo\dll^2\eps^2\}\nonumber\\
&&\hspace{20ex}=\left\{(\ql,\pl,\qz,\pz)/\ql=\dll, \pl\hspace{-1ex}=\ghpa(\qz,\pz,\a), \Iz\leq\frac{1}{2}\omo\dll^2\eps^2 \right\}\nonumber
\end{eqnarray}
Moreover, for all $(\qz,\pz)\in\g(\B_{\R^2}(0,\dll))$, $\ghpaep(\qz,\pz,\a)>\gcsep(\qz,\pz)$ holds if and only if
\begin{equation}\label{ext}
(\qz,\pz) \text{ is outside of }\C^\a_s:=\{(\qz,\pz)/\ghpa(\qz,\pz,\a)=\gcsep(\qz,\pz)\}.
\end{equation}
\end{lem}

\begin{figure}[!h]
\centering
\includegraphics[scale=0.45]{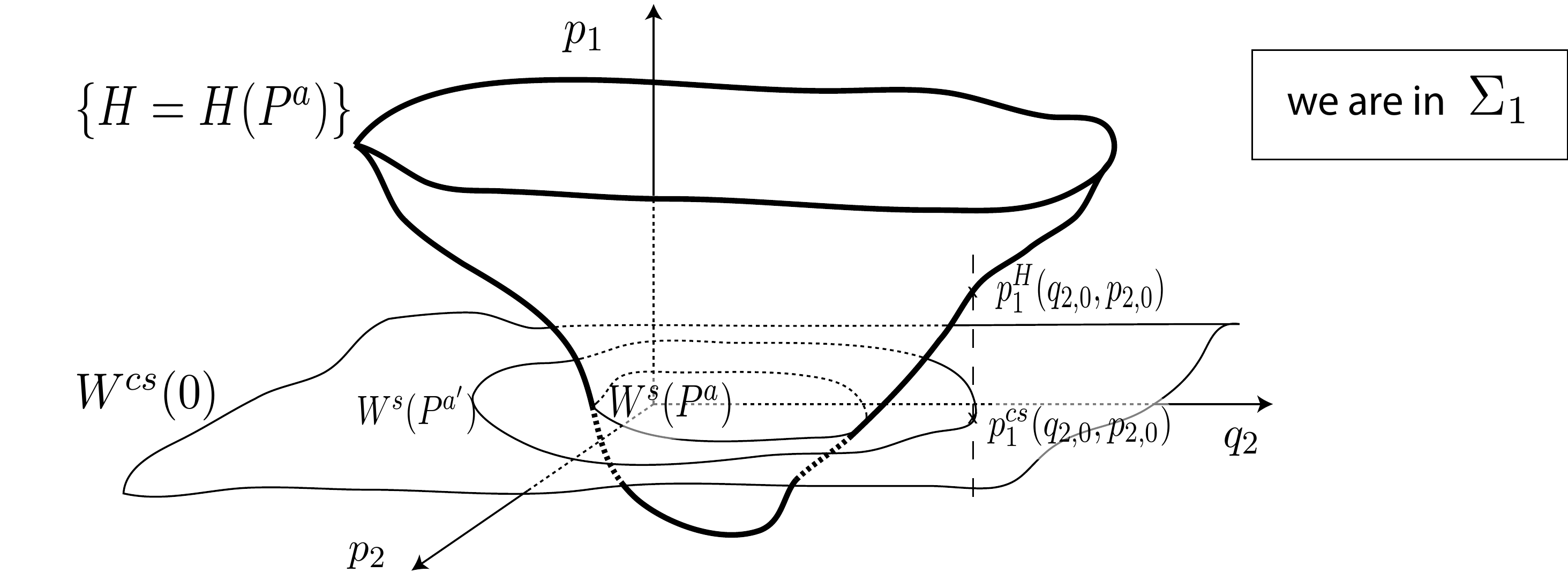}
\caption{In $\Sigl$, positions of the energy level set $\{\H=\H(\P^\a)\}$ and of the center-stable manifold $\CS$.}
\label{DESSIN6b}
\end{figure}

\textbf{Proof of the domain of $\ghpa$.} Recall that
$$|\hep(\qz^2+\pz^2,\a)|=|\frac{\om(\eps)}{2\eps^2}(\qz^2+\pz^2)-\H(\P^\a_{\ep},\ep)|.$$
On one hand, given that $\om$ is continuous, for $\eps$ sufficiently small we get
$$\Iz\leq\frac{1}{2}\omo\dll^2\eps^2 \quad \Longrightarrow \quad |\frac{\om(\eps)}{2\eps^2}(\qz^2+\pz^2)|\leq\frac{1}{2}\dll^2.$$
On the other hand, from $(i)$ of Lemma \ref{LemGammaCroissante}, we get
\begin{equation}
|\H(\P^\a_{\ep},\ep)|\leq\big|\H(\P^{\a}_{\ep},\ep)-\frac{\om(\eps)}{2\eps^2}\a\big| + \big|\frac{\om(\eps)}{2\eps^2}\a\big|\leq \a^2\M_1(\a)+\big|\frac{\om(\eps)}{2\eps^2}\a\big|.\nonumber
\end{equation}
So, there exists $\co$ such that if $\a\leq \co\eps^2\dll^2$ then $|\H(\P^\a_{\ep},\ep)|\leq\frac{1}{2}\dll^2$.
This proves the expression of the domain of $\ghpa$ claimed in the lemma, given that if $\Iz\leq\frac{1}{2}\omo\dll^2\eps^2$ and $\a\leq \co\eps^2\dll^2$ then $|\hep(\qz^2+\pz^2,\a)|\leq\dll^2$.
\cqfd

\vspace{1ex}

\textbf{Proof of the equivalence (\ref{ext}).} Recall that $\g$ was defined in Lemma \ref{concentriques}. First, let us prove that $\C^\a_s=\g(\C_\a)$. For that purpose, we work on the intersection of the domains of $\ghpaep$ and $\gcsep$. We consider then $\qz,\pz,\a$ such that $\Iz\leq\omo\dll^2\eps^2$ and $\a\leq\co\dll^2\eps^2$. Using Lemma \ref{LemGammaCroissante}, for $\eps$ sufficiently small we get that 
\begin{eqnarray}
&&\hspace{-6ex}(\dll,\gcsep(\qz,\pz),\qz,\pz)\in\{\qp/\H(\qp,\ep)=\H(\P^\a_{\ep},\ep)\}\nonumber\\
&\Longleftrightarrow& \HLep(0,(\g^{-1}(\qz,\pz))_{\qz}^2+(\g^{-1}(\qz,\pz))_{\pz}^2)=\HLep(0,\a).\nonumber\\
\hspace{10ex}&\Longleftrightarrow& \g^{-1}(\qz,\pz)\in \C_\a, \nonumber
\end{eqnarray}
which proves that $\C^\a_s=\g(\C_\a)$.

\vspace{1ex}

Let us consider a fixed value $({\qz}_0,{\pz}_0)$ outside of $\g(\C_\a)$ and suppose that moreover $({\qz}_0,{\pz}_0)$ belongs to $\g(\B_{\R^2}(0,\dll))$. Then Lemma \ref{concentriques} ensures that there exists $\a'>\a$ such that $({\qz}_0,{\pz}_0)\in\g(\C_\a')$. Thus given that $\C^\a_s=\g(\C_\a)$, the proof will be achieved if we show 
\begin{equation} \label{ghpaepcroit}
\ghpaep(\qz,\pz,\a')<\ghpaep(\qz,\pz,\a).
\end{equation}
Recall that $\ghpaep$ is defined by (\ref{defghpa}), $\ie$
\begin{equation}\label{defghpa2}
\frac{\om(\eps)}{2\eps^2}(\qz^2+\pz^2)-\H(\P^\a_{\ep},\ep)=\dll\ghpaep(\qz,\pz,\a)-\H_3(\dll,\ghpaep(\qz,\pz,\a),\qz,\pz,\ep).
\end{equation}
Recall also that we proved in the Step 2 of Lemma \ref{LemGrapheHPa}'s proof that 
$$\pl\mapsto\dll\pl-\H_3(\dll,\pl,\qz,\pz,\ep)$$
is strictly increasing for $|\pl|\leq\dll$. Finally, Lemma \ref{LemGammaCroissante} ensures that $\H(\P^\a_{\ep},\ep)<\H(\P^{\a'}_{\ep},\ep)$, so that (\ref{ghpaepcroit}) is a consequence of (\ref{defghpa2}). This achieves the proof of (\ref{ext}).
\cqfd

\subsection{Proof of Proposition \ref{PropReta}}\label{sub(iii)}

\textit{Let us define\label{NotReta}}
\begin{eqnarray}
\Reta_{\ep}:\{(\qz,\pz)/\qz^2+\pz^2\leq c\dll^2\eps^2, (\qz,\pz) \text{ outside of } \C^\a_s\}\hspace{-1ex}&\rightarrow&\hspace{-1ex} \{(\qz,\pz)/\qz^2+\pz^2\leq\omo\dll^2\eps^2\}\nonumber\\
(\qz,\pz)\hspace{-1ex}&\mapsto&\hspace{-1ex}(\Ret_{\ep})_{(\qz,\pz)}(\ghpa(\qz,\pz,\a),\qz,\pz).\nonumber
\end{eqnarray}
We proceed in several steps.

\vspace{1ex}

\textbf{Step 1. } Let us prove first that for any $c$, for $\eps$ sufficiently small, the set
\begin{equation}\label{Subset1}
\left\{\qp/\ql=\dll, 0<\pl-\gcsep(\qz,\pz)\leq \frac{1}{24\Mo}\dll,\qz^2+\pz^2\leq c\dll^2\eps^2\right\}
\end{equation}
is a subset of $\Ret$'s domain. 

Firstly, the equivalence (\ref{boncote}) of Lemma \ref{Lemgcsep} ensures that
$$\gcsep(\qz,\pz)<\pl \Rightarrow \etl>0.$$
Secondly, let us prove that
\begin{equation}\label{impliespo}
|\pl-\gcsep(\qz,\pz)|\leq \frac{1}{24\Mo}\dll \Rightarrow |\psilep^-(\dll,\pl,\qz,\pz)|\leq \frac{1}{24}\dll ;
\end{equation}
where we recall that $\FFep^{-1}=(\philep^-,\psilep^-,\phizep^-,\psizep^-)$. On one hand, the definition of $\gcsep$ ensures that $\psilep(\dll,\gcsep(\qz,\pz),\qz,\pz)=0$, and on the other hand from Lemma \ref{LemPrecTAF} we get that
{\small $$\psilep^-(\dll,\pl,\qz,\pz)-\psilep^-(\dll,\gcsep(\qz,\pz),\qz,\pz)\prec|\partial_{\pl}\psilep^-|(\dll,|\gcsep|(\qz,\pz)+\pl,\qz,\pz)|\pl-\gcsep(\qz,\pz)|.$$}
Then, from Lemma \ref{ProprPrec} and $(iii)$ of Proposition \ref{propFFcomplete}, we obtain
$$|\psilep^-(\dll,\pl,\qz,\pz)-\psilep^-(\dll,\gcsep(\qz,\pz),\qz,\pz)|\leq\Mo|\pl-\gcsep(\qz,\pz)|.$$
This achieves the proof of \eqref{impliespo}.

Finally, for $c$ and $\dll$ sufficiently small,
$$\qz^2+\pz^2\leq c\dll^2\nuu \Rightarrow \sqrt{\xiz^2+\etz^2}\leq\dll.$$
Indeed, this result is a consequence of $(vii)$ of Proposition \ref{propFFcomplete}.

\vspace{1ex}

\textbf{Step 2 : domain of $\Reta_{\ep}$.} Let us prove now that for $c$ sufficiently small, the set
$$\Sigl\cap\{\H=\H(\P^\a_{\eps})\}\cap\left\{\qp / (\qz,\pz)\text{ is outside of }\C^\a_s, \qz^2+\pz^2\leq c\dll^2\eps^2\right\}$$
is a subset of the set (\ref{Subset1}), and then also of $\Ret$'s domain.

Firstly, observe that (\ref{ext}) ensures that
$$(\qz,\pz)\text{ outside of }\C^\a_s \Rightarrow \ghpa(\qz,\pz,\a)>\gcsep(\qz,\pz).$$
Secondly, let us show that for $c$ and $\eps$ sufficiently small,
$$\qz^2+\pz^2\leq c\dll^2\eps^2 \Rightarrow \ghpa(\qz,\pz,\a)-\gcsep(\qz,\pz)\leq \frac{1}{24\Mo}\dll.$$
On one hand, we know that on the curve $\C^\a_s$, $\ghpa(\qz,\pz,\a)-\gcsep(\qz,\pz)=0$. On the other hand, from the definition of $\gcsep$, we get that
$$D\gcsep(\qz,\pz)=-\frac{D_{(\qz,\pz)}\psilep^-}{\partial_{\pl}\psilep^-}(\dll,\gcsep(\qz,\pz),\qz,\pz).$$
and so $(iii)$ of Proposition \ref{propFFcomplete} allows to obtain that
$$|D\gcsep(\qz,\pz)|\leq 2$$ 
for $\dll$ sufficiently small. From the definition (\ref{defghpa}) of $\ghpa$ we get in a similar way that
$$|D_{\qz,\pz}\ghpa(\qz,\pz,\a)|\leq \frac{C}{\eps}.$$
Thus we obtain that
$$\qz^2+\pz^2\leq c\dll^2\eps^2 \Rightarrow |\ghpa(\qz,\pz,\a)-\gcsep(\qz,\pz)|\leq \left(2+\frac{C}{\eps}\right)\sqrt{c}\dll\eps.$$
So, for $c$ sufficiently small, the result claimed in the summary above holds.

\vspace{1ex}

\textbf{Step 3 : range of $\Reta_{\ep}$.} Let us prove that the image of (\ref{Subset1}) through $\Ret$ is a subset of the set where the energy level set $\{\H=\H(\P^\a_{\ep})\}$ reads as the graph of $\ghpa$. 

From Proposition \ref{PropRet2}, we already know that the range of $\Ret$ is a subset of $\B(0,\dll)$. Let us show that the image of (\ref{Subset1}) is a subset of
$$\{\qp,\Iz\leq\frac{1}{2}\omo\dll^2\eps^2\}.$$
From (\ref{MajRet}) with $\nuu=\eps^2$, we get that if $\Iz\leq c\dll^2\eps^2$ then
$$|\frac{\om(\eps)}{2\eps^2}(\Ret_{\qz}(\qp,\ep)^2+\Ret_{\pz}(\qp,\ep)^2)|\leq M (c+\dll) \dll^2.$$
So the result claimed in the summary of the main steps above holds for $c$ and $\dll$ sufficiently small.
\cqfd

\def\nuu{\eps}
\section{Construction of an invariant curve for the restrictions of the first return map with the aid of a KAM theorem : proof of Proposition \ref{prophypKAM}}\label{SectionKAM}
\def\nuu{\eps}

This section is entirely devoted to the proof of Proposition \ref{prophypKAM}. 

In this part we consider $\a\leq\co\dll^2\eps^2$, and work in annulus of the form
$$\{(\qz,\pz)/\Iz=\qz^2+\pz^2\in[\cl(\dll)\dll\eps^2,\cz(\dll)\dll\eps^2]\},$$
where $\cl,\cz$ are lower than the $c$ of $\Reta_{\ep}$'s domain (see Proposition \ref{PropReta}) : the choice of $\cl,\cz$ is made in Lemma \ref{LemChoixci}. Here is an \textbf{outline of the proof of Proposition \ref{prophypKAM}} :

\begin{itemize}
\item \textbf{Parts \ref{SubHypReta0} and \ref{SubHypRetaEps}} are devoted to the proof of estimates concerning the $\Reta_{\ep}$. The results $(ii)$ and $(iii)$ will be the consequences of the estimates on the $\Reta_{(\eps,\eps^2,0)}$ showed in Part \ref{SubHypReta0}, and $(iv)$ a consequence of the upper bounds of $\Reta_{(\eps,\eps^2,\muu)}-\Reta_{(\eps,\eps^2,0)}$ computed in Part \ref{SubHypRetaEps} . 
\item In \textbf{Part \ref{SecEstimKAM},} we introduce the changes of coordinates which transforms $\Reta_{\ep}$ into $\Retap_{\ep}$. The proof of $(ii)$, $(iii)$ and $(iv)$ follows directly from the results of Parts \ref{SubHypReta0} and \ref{SubHypRetaEps}.  
\item\textbf{Part \ref{SubExacte}} is devoted to the proof of $(i)$.  
\end{itemize}

\subsection{Estimates for the first return map $\Reta_{(\eps,\eps^2,0)}$ of the normal form }\label{SubHypReta0}

In this part, we give first in Lemma \ref{LemRetFN} an explicit form of the $\Reta_{(\eps,\eps^2,0)}$ in polar coordinates. We then use this form to compute upper and lower bounds. Up to the change of coordinates that we will perform in Part \ref{SecEstimKAM}, Lemma \ref{LemMajkFN} states the upper bound $(iii)$ of Proposition \ref{prophypKAM} and Lemma \ref{LemChoixci} is the result $(ii)$. The latter lemma requires the proof of the two preliminary results of Lemmas \ref{LemTret0} and \ref{LemEncadrementGhpa}.

\begin{lem} \label{LemRetFN}
In polar coordinates, $\Reta_{(\eps,\eps^2,0)}$ reads
$$(\Reta_{\teta,(\eps,\eps^2,0)},\Reta_{\r,(\eps,\eps^2,0)})(\teta,\r)=(\teta+\Teta(\r,\a,\eps),\r),$$
with
\begin{eqnarray}
\Teta(\r,\a,\eps)&=&\frac{\om(\nuu)}{2\nuu^2}\Tret({\ghpa}_{\hspace{-1ex}(\eps,\eps^2,0)}(\r,0,\a),\r,0,(\eps,\eps^2,0))\nonumber\\
&&\hspace{-7ex}+\nuu\int^{\Tret({\ghpa}_{\hspace{-1ex}(\eps,\eps^2,0)}(\r,0,\a),\r,0,(\eps,\eps^2,0))}_{0}\hspace{-1ex}\partial_2\Q((\flot_{\ql}\hspace{-0.5ex}\+\flot_{\pl})(s,(\dll,{\ghpa}_{\hspace{-1ex}(\eps,\eps^2,0)}(\r,0,\a),\r,0),\eps),\r,\eps)ds,\nonumber
\end{eqnarray}
where ${\ghpa}_{\hspace{-0.5ex}(\eps,\eps^2,0)}(\r\cos\teta,\r\sin\teta,\a)$ and $\Tret(\pl,\r\cos\teta,\r\sin\teta,(\eps,\eps^2,0))$ (recall that $\Tret$ was defined in Proposition \ref{PropDefRet}) are independent of $\teta$.
\end{lem}

\textbf{Proof.} Recall that for $\muu=0$, the Hamiltonian system reads
\begin{equation}\label{SystHamFN}
\left\{\begin{array}{rcl}
\ql'(t)&=&-\pl+\frac{3}{2}(\ql+\pl)^2+\eps^2\partial_1\Q(\ql+\pl,\qz^2+\pz^2,\eps),\\
\pl'(t)&=&\ql-\frac{3}{2}(\ql+\pl)^2-\eps^2\partial_1\Q(\ql+\pl,\qz^2+\pz^2,\eps),\\
\qz'(t)&=&\frac{\om(\eps)}{2\eps^2}\pz+\eps^2\partial_2\Q(\ql+\pl,\qz^2+\pz^2,\eps)\pz,\\
\pz'(t)&=&-\frac{\om(\eps)}{2\eps^2}\qz-\eps^2\partial_2\Q(\ql+\pl,\qz^2+\pz^2,\eps)\qz.
\end{array}\right.
\end{equation}
This system satisfies
$$\frac{d}{dt}(\qz^2+\pz^2)=\frac{d}{dt}\Iz=0.$$
Then, the $(\ql,\pl)$ component of the flow reads $\flot_{(\ql,\pl)}(t,(\ql,\pl,{\Iz}_0),\nuu)$ for a fixed value $\Iz={\Iz}_0$. So for ${\qz}_0,
{\pz}_0$ such that ${\qz}_0^2+{\pz}_0^2={\Iz}_0$, 
$$\left(\begin{matrix}
\qz \\
\pz
\end{matrix}\right)(t)=\rot_{\teta(t,{\Iz}_0,\nuu)}\left(\begin{matrix}
{\qz}_0\\
{\pz}_0
\end{matrix}\right),$$
where
$$\teta(t,{\Iz}_0,\nuu):=\frac{\om(\nuu)}{2\nuu^2}t+\int^{t}_{0}\nuu^2\partial_2\Q((\flot_{\ql}+\flot_{\pl})(s,(\ql,\pl,{\Iz}_0),\nuu),\Iz,\nuu)ds.$$
The result claimed by Lemma \ref{LemRetFN} follows, up to the proof that ${\ghpa}_{\hspace{-0.5ex}(\nuu,\nuu^2,0)}(\r\cos\teta,\r\sin\teta,\a)$ and $\Tret(\pl,\r\cos\teta,\r\sin\teta,(\nuu,\nuu^2,0))$ are independent of $\teta$. 

Indeed, ${\ghpa}_{\hspace{-0.5ex}(\nuu,\nuu^2,0)}$ is defined as the unique $\ghpa$ such that
\begin{eqnarray}
  &&\H((\dll,\ghpa,\qz,\qz),(\nuu,\nuu^2,0))=\H(\P^\a_{(\nuu,\nuu^2,0)},(\nuu,\nuu^2,0))\nonumber\\
  &\Leftrightarrow& \dll\ghpa-\frac{1}{2}(\dll+\ghpa)^2-\nuu^2\Q(\dll+\ghpa,\r^2,\nuu)=\frac{\om(\nuu)}{2\nuu^2}\r^2-\H(\P^\a_{(\nuu,\nuu^2,0)},(\nuu,\nuu^2,0)).\nonumber
\end{eqnarray}
where the latter equation is independent of $\teta$. Similarly, $\Tret$ is defined by
$$\flotp_{\ql}(\Tret,(\dll,\pl,\qz,\pz),(\nuu,\nuu^2,0))=\dll,$$
where from the form (\ref{SystHamFN}) of the hamiltonian system we know that the $(\ql,\pl)$ component of the flow reads $\flot_{(\ql,\pl)}(t,(\ql,\pl,{\Iz}_0),\nuu)$. Then, the definitions of $\ghpa$ and $\Tret$ are independent of $\theta$ and we get the result. 
\cqfd

\begin{lem}\label{LemMajkFN}
There exists $M$ such that for $0\leq\a\leq\co$ and $\r^2\in[\cl\dll^2\nuu^2,\cz\dll^2\nuu^2]$, for all $k\leq\ko=\entiere(\frac{\No-1}{4})$, $\Teta$ satisfies
$$\left|\partial_r^k\Teta(\r,\a,\nuu)\right|\leq \frac{M}{\nuu^{k+2}}.$$
\end{lem}

\textbf{Proof.} On one hand, from the result of Lemma \ref{LemFlotp} the explicit formula of $\Teta$ in Lemma \ref{LemRetFN} reads
\begin{equation}\label{EcritTeta}
\begin{array}{rcl}
\Teta(\r,\a,\nuu)&=&\frac{\om(\nuu)}{2\nuu^2}\Tret({\ghpa}_{\hspace{-1ex}(\nuu,\nuu^2,0)}(\r,0,\a),\r,0,(\nuu,\nuu^2,0))\\
&&+\nuu^2\mathcal{F}\left(\Tret({\ghpa}_{\hspace{-1ex}(\nuu,\nuu^2,0)}(\r,0,\a),\r,0,(\nuu,\nuu^2,0)),{\ghpa}_{\hspace{-1ex}(\nuu,\nuu^2,0)}(\r,0,\a),\r,\nuu\right),
\end{array}
\end{equation}
with $\mathcal{F}$ and $\Tret$ in $\C^1(]-\epso,\epso[,\C^{\ko})$. On the other hand, recall that 
\begin{equation}\label{EcritGhpa}
{\ghpa}_{,\nup}(\qz,\pz,\a)=\ghpah(\qz,\pz,\hnup(\qz^2+\pz^2,\a),\nup),
\end{equation}
where $\ghpah$ belongs to $\C^1(]-\epso,\epso[^3,\An)$. Thus the only irregularity is in 
$$\hnup(\r^2,\a)=\frac{\om(\nuu)}{2\nuu^2}\r^2-\H(\P^\a_{\nup},\nup)$$
for $\nuu=0$. And we check that there exists $M_1$ such that for $k\leq\ko$ and $\r^2\in[\cl\dll^2\nuu^2,\cz\dll^2\nuu^2]$,
$$\left|\partial^k_\r\hnup(\r^2,\a)\right|\leq \frac{M_1}{\nuu^k}$$
holds. Finally, we obtain the existence of $M$ such that for $k\leq\ko$ and $\r^2\in[\cl\dll^2\nuu^2,\cz\dll^2\nuu^2]$,
$$\left|\partial^k_\r\Teta(\r,\a,\nuu)\right|\leq \frac{M}{\nuu^{k+2}}.$$

\cqfd

\def\nuu{\nu}
To prove $(ii)$ of Proposition \ref{prophypKAM}, we need a lower bound of $\partial_\r \Teta$. In view of the explicit form of $\Teta$ obtained in Lemma \ref{LemRetFN}, we first compute estimates of $\Tret$ in Lemma \ref{LemTret0} below and then of $\ghpa$ in the following Lemma \ref{LemEncadrementGhpa}.

\begin{lem} \label{LemTret0}
For the truncated Hamiltonian
$$\H(\qp,(\eps,0,0))=-\ql\pl+\frac{1}{2}(\ql+\pl)^3+\frac{\om(\eps)}{2\eps^2}(\qz^2+\pz^2),$$
the time $\Tret$ of first return to the section satisfies
$$\partial_{\pl}\Tret(\qp,(\eps,0,0))=\partial_{\pl}\Tret^0(\pl)<0. $$
\end{lem}

\textbf{Proof.} In the coordinates $(\qlp,\plp)$ (see Part \ref{SubNormalisation}), the truncated Hamiltonian reads
$$\frac{1}{2}(\plp^2-\qlp^2)+\frac{1}{2}\qlp^3+\frac{\om(\eps)}{2\eps^2}(\qz^2+\pz^2).$$
Recall that $\qz^2+\pz^2$ is constant, so that the flow follows the level sets
$$\frac{1}{2}(\plp^2-\qlp^2)+\frac{1}{2}\qlp^3=-\alpha \Leftrightarrow \plp=\pm\sqrt{\qlp^2-\qlp^3-\alpha}.$$
We get periodic orbits when $0<\alpha<\frac{4}{27}$ and the homoclinic orbit for $\alpha=0$. Let us denote by $(\qlp(t),\plp(t),\qz(t),\pz(t))$ the periodic solution associated with $\alpha$ and by $\tau(\alpha)$ its period. Then $\plp(\frac{1}{2}\tau(\alpha))=0$, and $\qlp(0)$ and $\qlp(\frac{1}{2}\tau(\alpha))$ are the two positive roots of $z^2-z^3-\alpha$. And $\frac{1}{2}\tau(\alpha)$ satisfies
$$\frac{1}{2}\tau(\alpha)=\int^{\frac{1}{2}\tau(\alpha)}_{0}\frac{\plp(t)}{\sqrt{\qlp(t)^2-\qlp(t)^3-\alpha}}dt=\int^{\qlp(\frac{1}{2}\tau(\alpha)):=z_2(\alpha)}_{\qlp(0):=z_1(\alpha)}\frac{1}{\sqrt{z^2-z^3+\alpha}}dz$$
given that $\qlp'(t)=\plp(t)$. By studying the map $z\mapsto z^2-z^3$, we obtain that $\partial_{\alpha}z_1(\alpha)>0$ and $\partial_{\alpha}z_2(\alpha)<0$ for all $\alpha\in]0,\frac{4}{27}[$. This proves that $\partial_{\alpha}\tau(\alpha)<0$, and thus that $\partial_{\pl}\Tret^0(\pl)<0$.


\cqfd

\begin{lem}\label{LemEncadrementGhpa}
Let $\ucl,\ucz$ be two fixed positive reals satisfying $\ucz-\ucl>4\dll$. 

Let us introduce $\ucl'',\ucz''$ such that $\ucl<\ucl''<\ucz''<\ucz-4\dll$ and define
$$\cl:=\frac{4\ucl''}{\omo},\quad \cz:=\frac{4\ucz''}{\omo}.$$
Then there exist $\co$ such that for all $\eps$ sufficiently small, for all 
$\a\in[0,\co\dll^2\eps^2]$ and $\r^2\in[\cl\dll^2\eps^2,\cz\dll^2\eps^2]$,
$${\ghpa}_{\hspace{-1ex}(\eps,\eps^2,0)}(\r,0,\a)\in[\ucl\dll,\ucz\dll], \quad \partial_h\ghpah\geq\frac{1}{2}\dll>0.$$ 
Moreover, on this domain, $\hnup$ satisfies $\hnup(\r^2,\a)\in[\ucl\dll^2,(\ucz-4\dll)\dll^2]$.
\end{lem}

\textbf{Proof.} Recall that Lemma \ref{LemGrapheHPa} defines $\ghpah$ and asserts that it is $\C^1$. 

\textbf{Step 1. Case $\ep=(\eps,\nuu=0,\muu=0)$.} Let us introduce $\ucl', \ucz'$ such that 
$$\ucl<\ucl'<\ucl''<\ucz''<\ucz'-4\dll<\ucz-4\dll.$$
In the case considered in this step, $\ghpah$ is defined by
$$\dll\ghpah-\frac{1}{2}(\dll+\ghpah)^3=\h.$$ 
Recall also that on its domain, $|\ghpah|\leq\dll$ holds. From these two results we get that, on one hand
\begin{equation}\label{Casnuuzero}
\partial_\h\ghpah(\r,0,\h,(\eps,0,0))\geq\frac{1}{\dll}>0
\end{equation}
holds on the domain of $\ghpah$. And on the other hand, for all $\h$ in $[\ucl'\dll^2,(\ucz'-4\dll)\dll^2]$,
\begin{equation}\label{casnuuzero}
\ghpah(\r,0,\h,\a,\nuu=0)\in [\ucl'\dll,\ucz'\dll].
\end{equation}

\vspace{1ex}

\textbf{Step 2. Case $\ep=(\eps,\eps^2,0)$.} Given that $\ghpah$ is $\C^1$ in term of $\nuu$ and $\eps$, from (\ref{casnuuzero}) we get that for $\eps$ sufficiently small, for all $\h$ in $[\ucl'\dll^2,(\ucz'-4\dll)\dll^2]$,
$$\ghpah(\r,0,\h,\a,(\eps,\eps^2,0))\in [\ucl\dll,\ucz\dll].$$
and from (\ref{Casnuuzero}) we obtain
$$\partial_\h\ghpah(\r,0,\h,(\eps,\eps^2,0))\geq\frac{1}{2\dll}>0.$$
Recall that
$$\hnup(\r^2,\a)=\frac{\om(\eps)}{2\eps^2}\r^2-\H(\P^\a_{\nup},\nup).$$
The choice of $\cl,\cz$ ensures that for $\eps$ sufficiently small if $\r^2$ belongs to $[\cl\dll^2\eps^2,\cz\dll^2\eps^2]$ then $\frac{\om(\eps)}{2\eps^2}\r^2$ is in $[\ucl''\dll^2,\ucz''\dll^2]$. And Lemma \ref{LemGammaCroissante} with $\nuu=\eps^2$ ensures that
$$\left|\H(\P^\a_{\nup},\nup)\right|\leq \frac{\om(\eps)}{2\eps^2}\a+\frac{\eps^2}{\eps^2}\Mo\a^2.$$
So, for $\co$ sufficiently small, if $ \a\leq\co\dll^2\eps^2$, then 
$$\left|\H(\P^\a_{\nup},\nup)\right|\leq \min(\ucl''-\ucl',(\ucz'-4\dll)-\ucz'')\dll^2.$$
\cqfd

\def\nuu{\eps}

\begin{lem} \label{LemChoixci}
If $\No\geq 5$, there exist $\co,\cl,\cz$ and $m,M>0$ such that for $\nuu$ sufficiently small, for all $\a\in[0,\co\dll^2\nuu^2]$ and all $\r^2\in[\cl\dll^2\nuu^2,\cz\dll^2\nuu^2]$,
$$-\frac{M}{\nuu^3}\leq\partial_\r\Teta(\r,\a,\nuu)\leq-\frac{m}{\nuu^3}.$$
Moreover, we can chose $\cl,\cz$ satisfying
$$M\dll<\cl<\cz\leq c,$$ 
with the constant $c$ introduced in Proposition \ref{PropReta} and $M$ is the constant of the upper bound \eqref{MajRet}.

And there exist $\ucl,\ucz$ such that on this domain, $\hnup$ and $\ghpa$ satisfies 
\begin{equation}\label{Boundghpa}
\hnup(\r^2,\a)\in[\ucl\dll^2,(\ucz-4\dll)\dll^2], \quad {\ghpa}_{(\nuu,\nuu^2,0)}(\r,\teta,\a)\in[\ucl\dll,\ucz\dll].
\end{equation}

%
\end{lem}


\textbf{Proof.} The existence of $M$ (without conditions on $a$ and $\r$) is a direct consequence of Lemma \ref{LemMajkFN} for $k=1$, given that we suppose $\No\geq 5$.

In order to prove the existence of $m$, we use the form (\ref{EcritTeta}) of $\Teta$ together with the form (\ref{EcritGhpa}) of $\ghpa$ : $\Teta$ reads
\begin{eqnarray}
\Teta(\r,\a,\nuu)&=&\frac{\om(\nuu)}{2\nuu^2}\Tret({\ghpa}_{\hspace{-1ex}(\nuu,\nuu^2,0)}(\r,0,\a),\r,0,(\nuu,\nuu^2,0))+\nuu^2\tilde{\mathcal{F}}\left(\hnup(\r^2,\a),\r,\nuu\right),\nonumber
\end{eqnarray}
where $\Tret$ and $\tilde{\mathcal{F}}$ belong to $\C^1(]-\epso,\epso[,\C^{\ko})$. Differentiating with respect to $\r$ we get
{\small 
\begin{eqnarray}
&&\hspace{-3ex}\partial_\r\Teta(\r,\a,\nuu)\nonumber\\
&\hspace{-2ex}=&\hspace{-2ex}\frac{\om(\nuu)}{2\nuu^2}\partial_{\pl}\Tret({\ghpa}_{\hspace{-1ex}(\nuu,\nuu^2\hspace{-0.5ex},0)}(\r,0,\a),\r,0,(\nuu,\nuu^2\hspace{-0.5ex},0))\hspace{-0.5ex}\left(\hspace{-0.5ex}\partial_\h\widetilde{\pl}^{\hspace{-0.5ex}\H}\hspace{-0.5ex}(\r,0,\a,(\nuu,\nuu^2\hspace{-0.5ex},0))\partial_\r\hnup(\r^2\hspace{-0.5ex},\a)\+\partial_\r\widetilde{\pl}^{\hspace{-0.5ex}\H}\hspace{-0.5ex}(\r,0,\a,(\nuu,\nuu^2\hspace{-0.5ex},0))\hspace{-0.5ex}\right)\nonumber\\
&&\hspace{-4ex}+\frac{\om(\nuu)}{2\nuu^2}\partial_\r\Tret({\ghpa}_{\hspace{-1ex}(\nuu,\nuu^2\hspace{-0.5ex},0)}(\r,0,\a),\r,0,(\nuu,\nuu^2\hspace{-0.5ex},0))\+\nuu^2\hspace{-0.5ex}\left(\hspace{-0.5ex}\partial_\h\tilde{\mathcal{F}}\left(\hnup(\r^2\hspace{-0.5ex},\a),\r,\nuu\right)\hspace{-0.5ex}\big(\frac{\om(\nuu)}{2\nuu^2}\r\big)\+\partial_\r\tilde{\mathcal{F}}\left(\hnup(\r^2\hspace{-0.5ex},\a),\r,\nuu\right)\hspace{-1ex}\right),\nonumber\\
 &\hspace{-2ex}=&\hspace{-2ex}\frac{\om(\nuu)}{2\nuu^2}\partial_{\pl}\Tret({\ghpa}_{\hspace{-1ex}(\nuu,\nuu^2,0)}(\r,0,\a),\r,0,(\nuu,\nuu^2,0))\partial_\h{\ghpah}(\r,0,\a,(\nuu,\nuu^2,0))\partial_\r\hnup(\r^2,\a) + \grandO{}\left(\frac{1}{\nuu^2}\right),\label{GrandO}
\end{eqnarray}}

where (\ref{GrandO}) holds because $\partial_{\pl}\Tret, \partial_\r\ghpah, \partial_\r\Tret, \partial_\h\tilde{\mathcal{F}}, \partial_\r\tilde{\mathcal{F}}$ are continuous and for any fixed choice of $\co,\cl,\cz$, $\hnup$ is bounded and Lemma \ref{LemEncadrementGhpa} ensures that $\ghpa$ is also. 

Let us show that the principal part in (\ref{GrandO}) admits an upper bound of the form $-\frac{m}{\nuu^3}$ for an appropriate choice of $\co,\cl$ and $\cz$. For that purpose, we use Lemmas \ref{LemTret0} and \ref{LemEncadrementGhpa}. 

On one hand, Lemma \ref{LemTret0} ensures that $\partial_{\pl}\Tret^0(\pl)<0$ when $\pl>\pl(\dll)$. Then, given that $\Tret$ is $\C^1$, for any set $\{\pl\in[\ucl\dll,\ucz\dll]\}$, there exists $m$ such that $\partial_{\pl}\Tret^0\leq -2m$ on this set. For $\nuu,\muu$ sufficiently small, we get
$$\partial_{\pl}\Tret(\pl,\qz,\pz,\ep)\leq-m<0 \text{ for all }\pl\in[\ucl\dll,\ucz\dll].$$
Let us chose $\ucl,\ucz$ satisfying moreover
$$4\frac{\ucz-4\dll}{\om_0}< c, \quad \ucl<\ucz+4\dll, \quad M\dll<\frac{4\ucl}{\om_0}, $$
where $M$ is the constant of the upper bound \eqref{MajRet}.
On the other hand, with this choice of $\ucl,\ucz$, Lemma \ref{LemEncadrementGhpa} ensures that there exists $\co,\cl$ and $\cz$ satisfying $\cl,\cz\leq c$ and $\cl>M\dll$ such that for all $\a\in[0,\co\dll^2\nuu^2]$ and $ \r^2\in[\cl\dll^2\nuu^2,\cz\dll^2\nuu^2]$,
$${\ghpa}_{\hspace{-1ex}(\nuu,\nuu^2,0)}(\r,0,\a)\in[\ucl\dll,\ucz\dll] \quad \quad \partial_\h{\ghpah}(\r,0,\a,(\nuu,\nuu^2,0))\geq\frac{1}{2\dll}>0.$$
Moreover, for $\r^2\in[\cl\dll^2\nuu^2,\cz\dll^2\nuu^2]$
$$\partial_\r\hnup(\r^2,\a)=2\frac{\om(\nuu)}{2\nuu^2}\r \geq \sqrt{\cl}\frac{\om(\nuu)}{\nuu}\dll>0.$$
We finally obtain that for $\a\in[0,\co\dll^2\nuu^2]$ and $\r^2\in[\cl\dll^2\nuu^2,\cz\dll^2\nuu^2]$,
$$\frac{\om(\nuu)}{2\nuu^2}\frac{\partial\Tret}{\partial\pl}({\ghpa}_{\hspace{-1ex}(\nuu,\nuu^2,0)}(\r,0,\a),\r,0,(\nuu,\nuu^2,0))\frac{\partial{\ghpah}}{\partial\h}(\r,0,\a,(\nuu,\nuu^2,0))\frac{\partial\hnup}{\partial\r}(\r^2,\a)\leq -\frac{m\om(\nuu)^2\sqrt{\cl}}{\nuu^3},$$
which, together with (\ref{GrandO}) achieves the proof of the Lemma.
\cqfd

\subsection{Upper bound $\C^{k}$ of $\Reta_{(\nuu,\nuu^2,\muu)}-\Reta_{(\nuu,\nuu^2,0)}$}\label{SubHypRetaEps}

Given that the KAM theorem of \cite{Moser} (see Theorem \ref{KAMth} in part \ref{subsubinvariantcurve}) is stated in polar coordinates, we need the following lemma, which gives upper bounds of the polar form $(\Reta_{\teta}(\teta,\r),\Reta_{\r}(\teta,\r))$ of $\Reta$ in terms of $(\Reta_{\qz}(\qz,\pz),\Reta_{\pz}(\qz,\pz))$. This computation relies on the fact that we work on the domain $$\{(\qz,\pz)/\qz^2+\pz^2\in[\cl\dll^2\nuu^2,\cz\dll^2\nuu^2]\},$$ which is away from $0$.

\begin{lem}\label{LemPolaire}
If $(\pl,\qz,\pz,\nup)$ belongs to the domain of $\Ret$ and $\cl\dll^2\nuu^2\leq\qz^2+\pz^2\leq\cz\dll^2\nuu^2$, then
$$\Ret_{(\qz,\pz)}(\pl,\qz,\pz,\nup)=\Ret_\r(\pl,\qz,\pz,\nup)\Rot_{\Ret_\teta(\pl,\qz,\pz,\nup)},$$
where $\Ret_\r$ is $\C^{\ko}$ and in $\R/2\pi\Z$,
$$\Ret_{\teta}(\pl,\qz,\pz,\nup)\equiv \frac{\om(\nuu)}{\nuu^2}\Tret(\pl,\qz,\pz,\nup)+\flotp_{\teta}(\Tret(\pl,\qz,\pz,\nup),(\dll,\pl,\qz,\pz),\nup),$$
holds with $\Tret$ and $\flotp_{\teta}$ $\C^{\ko}$ in $\C^1(]-\epso,\epso[^3,\C^{\ko})$ and there exists $M$ such that for every variable $x$ of $\flotp$,
$$\left|\partial^j_x\flotp_{\teta}\right|\leq\frac{M}{\nuu^j}\left(\sum_{k=0}^j|\partial^j_x\flotp_{\qz}|+|\partial^j_x\flotp_{\pz}|\right).$$
\end{lem}

\textbf{Proof.} Thank to the upper bound (\ref{MajRet}) of Proposition \ref{PropDefRet} (with $\nuu=\eps^2$), there exists $M$ such that if $(\pl,\qz,\pz,\nup)$ is in the domain of $\Ret$ and $(\qz,\pz)$ satisfies $\{(\qz,\pz)/\qz^2+\pz^2\in[\cl\dll^2\nuu^2,\cz\dll^2\nuu^2]\}$ then 
\begin{equation}\label{anneau}
\Ret(\pl,\qz,\pz,\nup)\in\left\{\qp/\qz^2+\pz^2\in[(\cl\dll^2-M\dll^3)\nuu^2,(\cz\dll^2+M\dll^3)\nuu^2]\right\}.
\end{equation}
On one hand, Lemma \ref{LemChoixci} ensures that $\cl>M\dll$, thus, denoting
$$\Ret_{\r}(\pl,\qz,\pz,\nup):=\sqrt{\Ret_{\qz}(\pl,\qz,\pz,\nup)^2+\Ret_{\pz}(\pl,\qz,\pz,\nup)^2},$$
$\Ret_{\r}$ is $\C^{\ko}$ as $\Ret$.
On the other hand, recall that
$$\Ret(\pl,\qz,\pz,\nup)=\Rot_{\frac{\om(\nuu)}{\nuu^2}\Tret(\pl,\qz,\pz,\nup)}\flotp(\Tret(\pl,\qz,\pz,\nup),(\dll,\pl,\qz,\pz),\nup),$$
so  
$$\Ret_{\teta}\equiv \frac{\om(\nuu)}{\nuu^2}\Tret(\pl,\qz,\pz,\nup)+\flotp_{\teta}(\Tret(\pl,\qz,\pz,\nup),(\dll,\pl,\qz,\pz),\nup),$$
where $\flotp_{\teta}$ is defined by
\begin{equation}
\Frac{\flotp_{\qz}+\I\flotp_{\pz}}{\flotp_{\r}}=\E^{\I\flotp_{\teta}}.\label{derivee}
\end{equation}
And (\ref{anneau}) ensures that $|\flotp_{\r}|^2\geq(\cl-M\dll)\dll^2\nuu^2$, then if $x$ is any variable of $\flotp$, we get
$$\left|\partial_x\flotp_{\teta}\right|\leq\frac{1}{\sqrt{(\cl-M\dll)}\dll\nuu}\left(|\partial_x\flotp_{\qz}|+|\partial_x\flotp_{\pz}|\right).$$
Differentiating (\ref{derivee}) many times, we obtain the upper bounds claimed above for the higher order derivatives of $\flotp_{\teta}$.

\cqfd

\def\nup{(\nuu,\nuu^2,\muu)}

\vspace{1ex}

\begin{lem} \label{LemMajTeta}
For $\a$,$\nuu$,$\muu$ sufficiently small and $0\leq k\leq\ko$,
\begin{eqnarray} 
 &&\NormeC{k}{\Reta_{\r}(\r,\teta,(\nuu,\nuu^2,\muu))-\Reta_{\r}(\r,\teta,(\nuu,\nuu^2,0))}\leq \left(\a^2+\frac{\muu}{\nuu^{k}}\right)M,\nonumber\\
&&\NormeC{k}{\Reta_{\teta}(\r,\teta,(\nuu,\nuu^2,\muu))-\Reta_{\teta}(\r,\teta,(\nuu,\nuu^2,0))}\leq \left(\frac{\a^2}{\nuu^{k+2}}+\frac{\muu}{\nuu^{k+1}}\right)M.\nonumber
\end{eqnarray}
\end{lem}

\textbf{Proof.} To compute an upper bound of $\NormeC{k}{\Reta_{(\nuu,\nuu^2,\muu)}-\Reta_{(\nuu,\nuu^2,0)}}$, we firstly wanted to use an upper bound of $\partial_\muu\Reta_{(\nuu,\nuu^2,\muu)}$ together with the Mean Value Theorem. Unfortunately, $\Reta_{\nup}$ is not smooth with respect to $\muu$ because of
$$\hnup(\r^2,\a)=\frac{\om(\nuu)}{\nuu^2}\r^2-\H(\P^\a_{\nup},\nup),$$
Indeed, $\P^\a_{\ep}$ is defined as $\P^\a_{\ep}:=\FFep(\C^\a)$ and $\FFep=\FF_{(\eps,\nuu,\muu)}$ is not smooth with respect to $\muu$. So in the following we are going to use $\partial_\muu$ together with the Mean Value Theorem as soon as it is possible, and we will complete by a use of $\partial_\h$ together with an upper bound of $|\hnup(\r^2,\a)-\tilde{\h}_{(\nuu,\nuu^2,0)}(\r^2,\a)|$. 

From Lemma \ref{LemGammaCroissante}, for $\a<\a_0$ and $\nuu,\muu$ sufficiently small we get the upper bound 
\begin{equation}
|\hnup(\r^2,\a)-\tilde{\h}_{(\nuu,\nuu^2,0)}(\r^2,\a)|=\left|\H(\P^\a_{(\nuu,\nuu^2,0)},(\nuu,\nuu^2,0))-\H(\P^\a_{\nup},\nup)\right|\leq \Mo \a^2.\nonumber 
\end{equation}

On one hand, recall that from the definitions of $\Reta$ (p.\pageref{NotReta}) and from the form \eqref{Defghpaep} of $\ghpaep$, $\Reta_\r$ reads 
\begin{eqnarray}
\Reta_\r(\r,\teta,\nup)\hspace{-2.5ex}&=&\hspace{-2.5ex}\Ret_\r\hspace{-0.5ex}\left(\ghpah(\r\cos\teta,\r\sin\teta,\hnup(\r^2,\a),\nup),\r\cos\teta,\r\sin\teta,\nup\hspace{-0.5ex}\right)\nonumber\\
 &:=& \Retah_\r(\r,\teta,\hnup(\r^2,\a),\nup)\nonumber
\end{eqnarray}
where $\Ret_\r$ and $\ghpah$ belong to $\C^1(]-\epso,\epso[^3,\C^{\ko})$ and thus so is the new function $\Retah_\r$ for $\r^2\in[\cl\dll^2\nuu^2,\cz\dll^2\nuu^2]$ and $\a\leq\co\dll^2\nuu^2$. Recall also that Lemma \ref{LemEncadrementGhpa} ensures that $\hnup$ is bounded on this domain. Then we get that for any $j+\ell\leq \ko$, 
 \begin{eqnarray}
&& \hspace{-3ex}\left|\partial_\r^j\partial_{\teta}^\ell\Retah_\r(\r,\teta,\hnup(\r^2,\a),\nup)-\partial_\r^j\partial_{\teta}^\ell\Retah_\r(\r,\teta,\hnup(\r^2,\a),(\nuu,\nuu^2,0))\right|\nonumber\\
&\leq&\hspace{-3ex} \supl_{\dindice{\r^2\in[\cl\dll^2\nuu^2,\cz\dll^2\nuu^2]}{\h\in[\ucl\dll^2,\ucz\dll^2]}}\left|\partial_\h\left(\partial_\r^j\partial_{\teta}^\ell\Retah_\r(\r,\teta,\h,\nup)\right)\right||\hnup(\r^2\hspace{-0.5ex},\a)\hspace{-0.5ex}-\hspace{-0.5ex}\tilde{\h}_{(\nuu,\nuu^2,0)}(\r^2\hspace{-0.5ex},\a)|\nonumber\\
&&+\supl_{\dindice{\r^2\in[\cl\dll^2\nuu^2,\cz\dll^2\nuu^2]}{\h\in[\ucl\dll^2,\ucz\dll^2]}}\left|\partial_\muu\left(\partial_\r^j\partial_{\teta}^\ell\Retah_\r(\r,\teta,\tilde{\h}_{(\nuu,\nuu^2,0)}(\r^2,\a),\nup)\right)\right|\cdot \muu. \nonumber
\end{eqnarray}
And given that $|\partial_\r^i\hnup(\r^2,\a)|\leq \frac{M}{\eps^i}$ (see the proof of Lemma \ref{LemMajkFN}), we obtain
$$\NormeC{k}{\Reta_{\r}(\r,\teta,(\nuu,\nuu^2,\muu))-\Reta_{\r}(\r,\teta,(\nuu,\nuu^2,0))}\leq \left(\a^2+\frac{\muu}{\nuu^{k}}\right)M.$$

\vspace{1ex}

On the other hand, recall that $\Reta_\teta(\r,\teta,\nup)$ reads
\begin{eqnarray}
&&\frac{\om(\nuu)}{2\nuu^2}\Tret(\ghpah(\r\cos\teta,\r\sin\teta,\hnup(\r^2,\a),\nup),\qz,\pz,\nup)\nonumber\\
&\hspace{-2ex}+&\hspace{-3ex}\flotp_{\teta}(\Tret(\ghpah(\r\cos\teta,\r\sin\teta,\hnup(\r^2\hspace{-0.5ex},\a),\nup),\qz,\pz,\nup),(\dll,\pl,\qz,\pz),\nup) \nonumber
\end{eqnarray}
where $\Tret$ and $\ghpah$ belong to $\C^1(]-\epso,\epso[^3,\C^{\ko})$. From Lemma \ref{LemPolaire} and given that $\flotp_{\qz,\pz}$ belongs $\C^{\ko}$, with the same computations as that of the upper bounds of $\Reta_\r(\r,\teta,\nup)$, we get
\begin{eqnarray}
\NormeC{k}{\Reta_{\teta}(\r,\teta,(\nuu,\nuu^2,\muu))-\Reta_{\teta}(\r,\teta,(\nuu,\nuu^2,0))}&\leq&\frac{1}{\nuu}M\cdot\left(\a^2+\frac{\muu}{\nuu^{k}}\right)+\frac{1}{\sqrt{\nu}^{k}}M\cdot\left(\a^2+\frac{\muu}{\nuu^{k}}\right),\nonumber\\
&\leq&   M'\cdot\left(\frac{\a^2}{\nuu^{k+2}}+\frac{\muu}{\nuu^{2k+2}}\right).\nonumber
\end{eqnarray}
\cqfd

\subsection{Change of coordinates, proof of $(ii)$,$(iii)$ and $(iv)$ of Proposition \ref{prophypKAM}} \label{SecEstimKAM}

\def\nup{\underline{\nuu}}

\textbf{Change of coordinates.} Denoting by $\left[\frac{1}{\nuu^2}\right]$ the integer part of $\frac{1}{\nuu^2}$, let us define 
\begin{equation} \label{defnuup}
\nuup:=\frac{1}{\left[\frac{1}{\nuu^2}\right]},
\end{equation}
and perform the following change of coordinates.
$$(\teta,\r):=\left(\frac{1}{\nuup}\q,\sqrt{\nuup}\rho\right), \quad \q\in\R/2\pi\Z, \ro^2\in\frac{\nuu^2}{\nuup}[\cl\dll^2,\cl\dll^2].$$
Observe that $1-\nuu^2\leq\frac{\nuu^2}{\nuup}\leq1$. Let us denote by $\Retap_{\nup}$ the map $\Reta_{\nup}$ expressed in those new coordinates. Thus $\Retap_{\nup}$ is defined for $\ro^2\in[\cl\dll^2,(1-\nuu^2)\cz\dll^2]$, $\ie$ for $\nuu$ sufficiently small $\Retap_{\nup}$ is defined on a set $\ro^2\in[\cl\dll^2,d_1\dll^2]$ independent of $\nuu$. $\Retap_{\nup}$ reads

\def\nup{(\nuu,\nuu^2,\muu)}

\begin{eqnarray}
\Retap_{\nup}(\q,\ro)&=&\left(\q+\nuup\Teta(\sqrt{\nuup}\ro,\a,\nuu)+\nuup(\Reta_{\nup}-\Reta_{(\nuu,\nuu^2,0)})_{\teta}(\frac{1}{\nuup}\q,\sqrt{\nuup}\ro)\right., \nonumber\\
 & & \left. \ro+ \frac{1}{\sqrt{\nuup}}(\Reta_{\nup}-\Reta_{(\nuu,\nuu^2,0)})_{\teta}(\frac{1}{\nuup}\q,\sqrt{\nuup}\ro)\right), \nonumber\\
 &:=&\left(\q+\alpha^\a_\nuu(\q,\ro)+F^\a_{\nup}(\q,\ro) , \r+G^\a_{\nup}(\q,\ro)\right).
\end{eqnarray}
%
%


\vspace{2ex}

\textbf{Proof of $(ii)$ of Proposition \ref{prophypKAM}.} $$\partial_\ro\alpha^\a_{\nuu}=\nuup\sqrt{\nuup}\partial_\r\Teta(\sqrt{\nuup}\ro,\a,\nuu).$$
From Lemma \ref{LemChoixci} and the definition of $d_1$, we get that for $\ro^2\in[\cl\dll^2,d_1\dll^2]$,
$$-\frac{M}{\nuu^3}\nuup\sqrt{\nuup}\leq\frac{\partial\alpha^\a_\nuu}{\partial\ro}\leq-\frac{m}{\nuu^3}\nuup\sqrt{\nuup}.$$
Observe that from the definition (\ref{defnuup}) of $\nuup$, we get that $1\leq\frac{\nuup}{\nuu^2}\leq 2$ if $\nuu^2\leq\frac{1}{2}.$ We then obtain
$$-2\sqrt{2}M\leq\frac{\partial\alpha^\a_\nuu}{\partial\ro}\leq-m ;$$
so $(ii)$ holds with $m_0=\max(2\sqrt{2}M,\frac{1}{m})$.

\vspace{1ex}

\textbf{Proof of $(iii)$ of Proposition \ref{prophypKAM}.} Given that if $\ro^2\in[\cl\dll^2,d_1\dll^2]$ then $\r^2\in[\cl\dll^2\nuu^2,\cz\dll^2\nuu^2]$, Lemma \ref{LemMajkFN} ensures that for $k\leq\ko$,
$$\left|\frac{\partial^k\alpha^\a_\nuu}{\partial\ro^k}(\q,\ro)\right|=\left|\nuup\sqrt{\nuup}^k\frac{\partial\Teta}{\partial\r}(\sqrt{\nuup}\ro,\a,\nuu)\right|\leq \nuup\sqrt{\nuup}^k\frac{M}{\nuu^{k+2}} \leq \sqrt{2}^k M \quad\text{for }\nuu\leq\frac{1}{2}.$$

\textbf{Proof of $(iv)$ of Proposition \ref{prophypKAM}.} Lemma \ref{LemMajTeta} ensures that for $i+j\leq k\leq\ko$,
\begin{eqnarray}
\left|\partial_\q^i\partial_\ro^j F^\a_{\nup}(\q,\ro)\right|&=&\nuup\frac{\sqrt{\nuup}^j}{\nuup^i}\left|\partial_\teta^i\partial_\r^j(\Reta_{\nup}-\Reta_{(\nuu,\nuu^2,0)})_{\teta}(\frac{1}{\nuup}\q,\sqrt{\nuup}\ro)\right|,\nonumber\\
 &\leq& \nuup\frac{\sqrt{\nuup}^j}{\nuup^i} M\cdot\left(\frac{\a^2}{\nuu^{i+j+2}}+\frac{\muu}{\nuu^{2(i+j+1)}}\right)\leq M\cdot\left(\frac{\a^2}{\nuu^{3k}}+\frac{\muu}{\nuu^{6k}}\right);\nonumber
\end{eqnarray}
and
\begin{eqnarray}
\left| \partial_\q^i\partial_\ro^j G^\a_{\nup}(\q,\ro)\right|&=&\nuup\frac{\sqrt{\nuup}^j}{\nuup^i}\left|\partial_\teta^i\partial_\r^j (\Reta_{\nup}-\Reta_{(\nuu,\nuu^2,0)})_{\r}(\frac{1}{\nuup}\q,\sqrt{\nuup}\ro)\right|,\nonumber\\
 &\leq& \nuup\frac{\sqrt{\nuup}^j}{\nuup^i}M \cdot \left(\a^2+\frac{\muu}{\nuu^{i+j}}\right)\leq M\cdot\left(\frac{\a^2}{\nuu^{2k-2}}+\frac{\muu}{\nuu^{6k-4}}\right).\nonumber
\end{eqnarray}
So finally, for $0\leq k\leq\ko-1$,
$$\NormeC{k}{F^\a_{\nup}}+\NormeC{k}{G^\a_{\nup}}\leq M\cdot\left(\frac{\a^2}{\nuu^{3k}}+\frac{\muu}{\nuu^{6k}}\right).$$

\cqfd
\def\nup{\underline{\nuu}}

\subsection{Proof of $(i)$ of Proposition \ref{prophypKAM} : the $\Retap$ are exact maps}\label{SubExacte}

We first prove that the $\Retap$ are area-preserving maps (Lemma \ref{LemAire}), and then that for every $\Retap$ there exist some Jordan curves intersecting their range through $\Retap$ (Lemma \ref{LemExacte}). These two results together ensures that $(i)$ of Proposition \ref{prophypKAM} holds.

\begin{lem} \label{LemAire}
$\Retap_{\nup}$ is an area-preserving map.
\end{lem}

\textbf{Proof.} $\Reta_{\nup}$ is symplectic given that it is a first return map associated to a Hamiltonian flow. Moreover, we verify that the change of coordinates  
$$\varphi:(\ro,\q)\mapsto(\sqrt{\nuup}\ro,\frac{1}{\nuup}\q)$$
is symplectic. Then $\Retap_{\nup}$ is symplectic and thus area-preserving.
\cqfd

\vspace{1ex}

\begin{lem}\label{LemExacte}
We fix $\Delta>0$.

For $\nup$ sufficiently small with respect to $\Delta$, every Jordan curve $\C$ in the set $$\{(\q,\ro)/\ro^2\in[\cl'\dll^2,d_1\dll^2]\}$$ and of the form
$$\C=\{(\q,\ro), \ro=f(\q)\} \quad\text{ with for all } \q\in\R/2\pi\Z, |f'(\q)|\leq\Delta, $$
intersects its range $\Retap_{\nup}(\C)$.
\end{lem}



\textbf{Proof of an upper bound for the map $\Reta_{\nup}$.}
Firstly, we prove that there exists a constant $M''$ such that the upper bound 
\begin{equation}
\left|(\Retap_{\nup})_{\ro}(\q,\ro)-\ro\right|\leq M''\muu\nuu^{\No}\label{MajRetap}
\end{equation}
holds in the domain $\{(\q,\ro)/\ro^2\in[\cl\dll^2,d_1\dll^2]\}$. For that purpose, we use two results : 
\begin{enumerate}
\item in the coordinates $(\qz,\pz)$, there exists $M$ such that on the domain $\B(0,\roo)$ the flow of $\H(.,(\nuu,\nuu^2,\muu))$ satisfies
$$\left|\frac{d(\qz^2+\pz^2)}{dt}\right|\leq M \muu\nuu^{\No+2} ;$$
\item on the domain $\{(\teta,\r)/\r^2\in[\cl\dll^2\nuu^2,\cz\dll^2\nuu^2]\}$ the time $\Tret$ of first return to $\Sigl$ is bounded : this is the consequence of the $C^k$-smoothness of $\Tret(\pl,\ql,\pz,\nup)$ together with the fact that $\ghpa$ is bounded on this domain thank to \eqref{Boundghpa} of Lemma \ref{LemChoixci}.
\end{enumerate}
From $(i)$ and $(ii)$ we derive that
$$\left|(\Reta_{\nup})_{\qz}(\qz,\pz,\a)^2+(\Reta_{\nup})_{\pz}(\qz,\pz,\a)^2- (\qz^2+\pz^2)\right|\leq M \muu\nuu^{\No+2} (Sup  |\Tret|)\leq M'\muu\nuu^{\No+1}.$$ 
Thus in polar coordinates we get
\begin{eqnarray}
&&\left|(\Reta_{\nup})_\r(\r,\teta,\a)^2-\r^2\right|\leq M'\muu\nuu^{\No+2}. \nonumber\\
&\Leftrightarrow&\left|(\Reta_{\nup})_\r(\r,\teta,\a)-\r\right|\left|(\Reta_{\nup})_\r(\r,\teta,\a)+\r\right|\leq M'\muu\nuu^{\No+2}. \nonumber\\
& \Rightarrow&\left|(\Reta_{\nup})_\r(\r,\teta,\a)-\r\right|\leq M''\muu\nuu^{\No+1} \quad \text{given that }\r\geq\sqrt{\cl}\dll\nuu.\nonumber
\end{eqnarray}
Which, in coordinates $(\q,\ro)$ reads
$$\left|(\Reta_{\nup})_\r(\r,\teta,\a)-\r\right|=\left|(\Reta_{\nup})_{\ro}(\frac{1}{\nuup}\q,\sqrt{\nuup}\ro)-\sqrt{\nuup}\ro\right|=\left|\sqrt{\nuup}(\Retap_{\nup})_{\r}(\q,\ro)-\sqrt{\nuup}\ro\right|.$$
This completes the proof of \eqref{MajRetap}.

\vspace{1ex}

\begin{figure}[!h]
\centering
\includegraphics[scale=0.68]{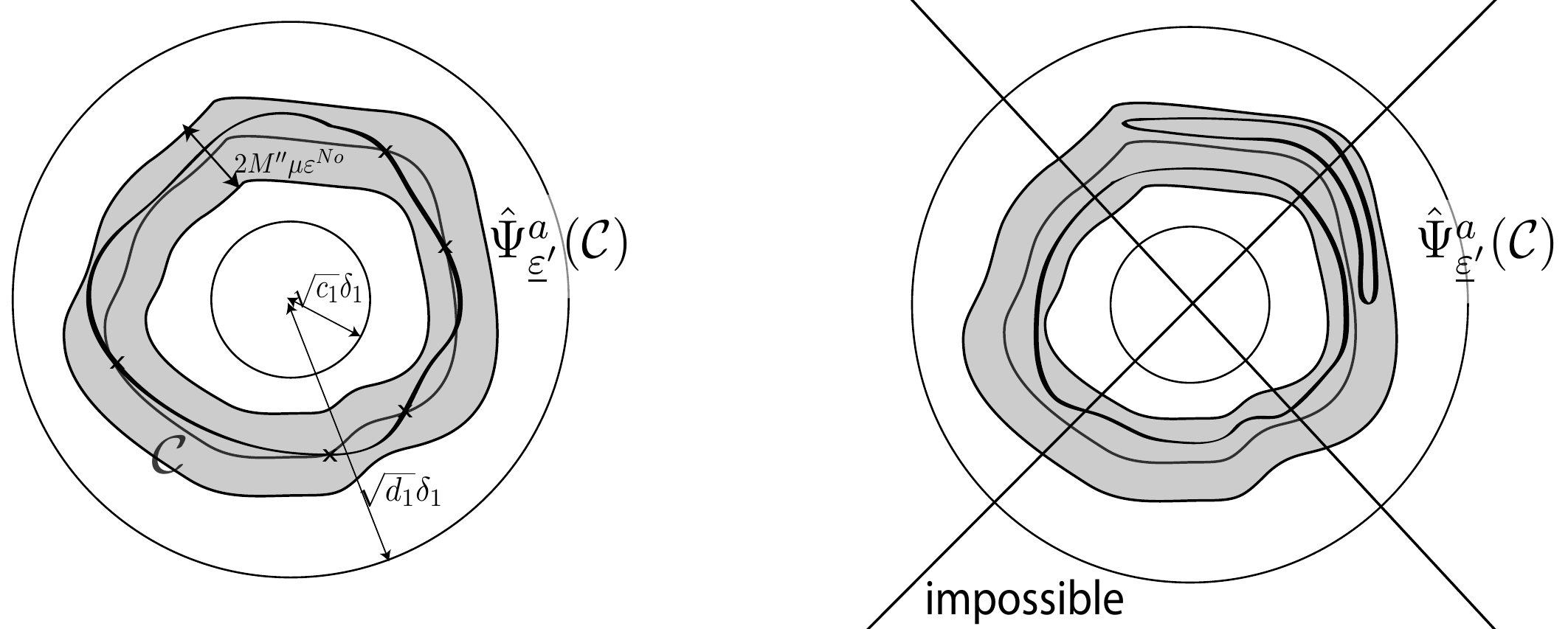}
\caption{Range of $\C$ through $\Retap_{\protect\underline{\varepsilon}}$.}
\label{DESSIN10bis}
\end{figure}

\textbf{Proof of Lemma \ref{LemExacte}.} Consider a curve $\C$ as described in the statement of the Lemma. From (\ref{MajRetap}) we get that $\Retap_{\nup}(\C)$ is in the tube
$$\{(\q,\ro),\ro\in[f(\q)-M''\muu\nuu^{No},f(\q)+M''\muu\nuu^{No}]\}$$
(see Figure \ref{DESSIN10bis}) whose area admits the upper bound $2 M''\muu\nuu^{No}\cdot$(length of the curve $\C$). Where the length of $\C$ reads
$$\int^{2\pi}_{0}\sqrt{f'(\q)^2+f(\q)^2}d\q\leq 2\pi(\Delta+\sqrt{d_1}).$$

Suppose that $\Retap_{\nup}(\C)\cap\C=\emptyset$ holds. Given that $\C$ is a Jordan curve, so is $\Retap_{\nup}(\C)$, and then their inside and outside sets are well-defined. From Lemma \ref{LemAire} we get that their inside sets have the same area, and thus necessarily $\Retap_{\nup}(\C)$ is in the outside set of $\C$. Then $\Retap_{\nup}(\C)$ is in the tube $\{(\q,\ro),\ro\in[f(\q),f(\q)+M''\muu\nuu^{No}]\}$, and so is its inside set whose area is greater than $\pi\cl$. 

So for $\muu,\nuu$ sufficiently small with respect to $\cl,d_1$ and $\Delta$ the area of the tube is lower than the area of the curve $\Retap_{\nup}(\C)$, so that necessarily $\Retap_{\nup}(\C)\cap\C\neq\emptyset$ is impossible.

\cqfd

%

\def\nuu{\nu}

\setcounter{section}{0}
\renewcommand{\thesection}{\Alph{section}}

\section{Appendix. Proof of Theorem \ref{ThmNF}}\label{AppendixA}
This appendix is devoted to the statement and proof of the following general Hamiltonian Normal Form Theorem \ref{ThmNF} used in Part \ref{PartieFormeNormale}. We thank G\'erard Iooss who suggested us this version of the Normal Form Theorem and another proof based on generatrix functions.

The main novelty of this result is that it does not require the linear part to be semi-simple. Moreover, unlike the usual Birkhoff Normal Form results, this theorem holds for elliptic but also non elliptic fixed points. 

\begin{thm}[Normal form theorem]\label{ThmNF}
Let $\R^{2m}$ be endowed with the symplectic form $\Omega_m$ given by
\begin{equation}\label{EqMatrixOmegam}
\Omega_m(x,y) = \scal{J_mx}{y} \qquad \textrm{where}\qquad  J_m= \left (\begin{array}{cc}           0 & I_m \\           -I_m & 0                      \end{array}    \right).
\end{equation}
Let ${\cal U}_\Lambda$ be an open set of a Banach space $\Lambda$ such that $0\in{\cal U}_\Lambda$. Let $\Hcal_\lb$
with $\lb\in{\cal U}_\Lambda$, be a $C^1$ one parameter family of $\C^k$ (\textit{resp. analytic and in $\An\left(\B_{\R^4}(0,\rho),\R\right)$}) Hamiltonian such that $D_x\Hcal_\lb(0)=0$. We denote by $\Hcal_{2,\lb}(x)=\frac{1}{2} D_{x,x}^2\Hcal_\lb(0).[x,x]$ the quadratic part of $\Hcal_\lb$ and by $L_0$ the linear part at $\lb=0$ of the associated hamiltonian vector field, $\ie$ $L_0x=J_m\nabla_{\hsp[-.5]x}\Hcal_{2,0}(x)$.

\vspace{1ex}

\noindent Then, for all $n$ and $k$ such that $k-1\geq n\geq 3$ (\textit{resp. for all $n\geq3$}), there exists a $\C^1$ one parameter family $\phi_{n,\lb}$ of canonical analytic mappings in $\An\left(\B_{\R^4}(0,\rho'),\R\right)$, close to identity, defined in a neighborhood of $0$ in $\Lambda$, such that in the neighborhood of $0$ in $\R^{2m}$, ${\cal H}_\lb$ satisfy
$$\begin{array}{ll}  
\Hcaltilde_{\lb}(\xtilde) = {\cal H}_\lb(\phi_{n,\lb}(\xtilde)) = \Hcal_{2,0}(\xtilde)+\Ncal_{n,\lb} (\xtilde)+\Rcal_{n+1,\lb} (\xtilde)
\end{array},$$
where  $\Ncal_{n,\lb}$ is a real polynomial of degree less than $n$, and $\Rcal_{n,\lb}$ is $\C^k$ (\textit{resp. belongs to $\An\left(\B_{\R^4}(0,\rho'),\R\right)$}) and 
\begin{eqnarray}
&&\hspace{-7ex}\Ncal_{n,\lb}(\xtilde)={\cal O}(\lb |\xtilde|^2+|\xtilde|^3), \qquad \Rcal_{n+1,\lb}(\xtilde)={\cal O}\left(|\xtilde|^{n+1}\right),\nonumber\\
&&\hspace{-7ex}\Ncal_{n,\lb}(\E^{t L_0^*} \xtilde)=\Ncal_{n,\lb}(\xtilde) \ \ \mbox{for all } t\in\R, \quad \text{or equivalently} \quad \{\Hcal_{2,0}(-J_m),\Ncal_{n,\lb}\}=0 \label{NFCrit}
\end{eqnarray}
Moreover, the coefficients of $\Ncal_{n,\lb}$ are $C^{1}$ functions of $\lb$.
\end{thm}

\begin{rem}
Observe that in the case of the Birkhoff normal form, the normal form $\Ncal$ belongs to the kernel of $\{\Hcal_{2,0},\cdot\}$ while here it lies in the kernel of $\{\Hcal_{2,0}\circ(-J),.\}$. This is due to the fact that in the Birkhoff normal form theorems the fixed point is elliptic and hence in that case $\Hcal_{2,0}\circ(-J)=\Hcal_{2,0}$.
\end{rem}

\textbf{Proof.} We begin with a summary of the strategy of proof.

\noindent{\bf Strategy of proof.} We perform the proof by induction on the statement of the theorem at degree $n=\ell\geq 2$. We define $\phi_{1,\lb}=Id$ and at each order $\ell$ we construct the change of coordinates $\phi_{\ell,\lb}=\phi_{\ell-1,\lb}\circ\varphi_{\ell,\lb}$. We look for $\varphi_{\ell,\lb}$ as the Lie transform of a homogeneous polynomial $S_{\ell,\lb}$ of degree $n$. Our aim is to construct $S_{\ell,\lb}$ (and hence $\varphi_{\ell,\lb}$) such that the polynomial part of degree $\ell$ of the new Hamiltonian is as simple as possible, $\ie$ with as few monomials as possible. Let us denote by $\Hcaltilde_{\ell-1,\lb}$ the Hamiltonian obtained after the $\ell^{th}$ step of the induction. It reads
\begin{equation}
\Hcaltilde_{\ell-1,\lb}=\Hcal_\lb\circ\phi_{\ell-1,\lb}=\Hcal_{2,0}+\Ncal_{\ell-1,\lb}+P_{\ell,\lb}+R_{\ell+1,\lb},\label{BeginInd}
\end{equation}
where $\Ncal_{\ell-1,\lb}$ is of the form \eqref{NFCrit}, $P_{\ell,\lb}$ is an homogeneous polynomial of degree $\ell$ whose coefficients are $C^1$ in term of $\lb$, and $R_{\ell+1,\lb}={\cal O}(|x|^{\ell+1})$ is $C^1$ in term of $\lb$. Our aim is then to construct $\varphi_{\ell,\lb}$ such that 
\begin{equation}
\Hcaltilde_{\ell,\lb}=\Hcal_\lb\circ\phi_{\ell-1,\lb}\circ\varphi_{\ell,\lb}=\Hcal_{2,0}+\Ncal_{\ell-1,\lb}+N_{\ell,\lb}+\Rcal_{\ell+1,\lb},\label{EndInd}
\end{equation}
where $N_{\ell,\lb}$ is an homogenous polynomial of degree $\ell$ satisfying the criterium \eqref{NFCrit} and $\Rcal_{\ell+1,\lb}={\cal O}(|x|^{\ell+1})$.

We will first (see {\bf Step 1} below) compute the equation satisfied by $S_{\ell,\lb}$ and $N_{\ell,\lb}$ in the set of homogeneous polynomial of degree $\ell$. We prove ({\bf Step 2}) that the family $\phi_{\ell,\lb}$ inherits the regularity of the family $S_{\ell,\lb}$ in term of $\lb$. Then, we construct ({\bf Step 3}) a scalar product on the space of homogeneous polynomial of degree $\ell$ which will allow to get uniqueness in the construction of $S_{\ell,\lb}$ and $N_{\ell,\lb}$. Finally, thanks to the uniqueness of the decomposition we prove ({\bf Step 4}) by an Implicit Function Theorem the existence of smooth families $S_{\ell,\lb}$ and $\N_{\ell,\lb}$ satisfying the equation of the Step 1.
 
\vspace{1ex} 

\noindent{\bf Step 1.1 Normalization at order 2.}  The initial Hamiltonian $\Hcal_\lb$ reads
$$\Hcal_\lb=\Hcal_{2,\lb}+R_{3,\lb}=\Hcal_{2,0}+P_{2,\lb}+R_{3,\lb},$$
where $P_{2,\lb}$ is an homogenous polynomial of degree 2 with $C^1$ coefficients in term of $\lb$ verifying $P_{2,0}\equiv 0$, and $R_{3,\lb}={\cal O}(|x|^3)$.

We recall that if we denote by $\varphi_{2,\lb}^t$ the Hamiltonian flow of the function $S_{2,\lb}$ then the Lie transform $\varphi_{2,\lb}$ of $S_{2,\lb}$ is defined by $\varphi_{2,\lb}:=\varphi_{2,\lb}^1$. We moreover recall that then for any regular function $G$,
\begin{equation}\label{ResLie}
\frac{d}{dt}(G\circ\varphi_{2,\lb}^t)=\{G,S_{2,\lb}\}\circ\varphi_{2,\lb}^t
\end{equation}
holds. Thus, we get
$$\Hcal_{2,\lb}\circ\varphi_{2,\lb}=\Hcal_{2,\lb}+\{\Hcal_{2,\lb},S_{2,\lb}\}+\int^{1}_{0}(1-s)\{\{\Hcal_{2,\lb},S_{2,\lb}\},S_{2,\lb}\}\circ\varphi_{2,\lb}^s ds +{\cal O}(|x|^3).$$
So that $N_{2,\lb}$ and $S_{2,\lb}$ satisfy the following equation in the space of homogeneous polynomials of degree 2 :
\begin{equation}\label{Homol2}
N_{2,\lb}-\{\Hcal_{2,\lb},S_{2,\lb}\}-\int^{1}_{0}(1-s)\{\{\Hcal_{2,\lb},S_{2,\lb}\},S_{2,\lb}\}(\E^{sJ_m\nabla S_{2,\lb}})ds=P_{2,\lb}.
\end{equation}
 
\vspace{1ex} 
 
\noindent{\bf Step 1.2. Equation at order $\ell\geq 3$}. If $\Hcal_{\ell,\lb}$ is of the form \eqref{BeginInd} after the $(\ell-1)^{th}$ step of the induction, then we get
\begin{eqnarray}
\Hcal_{\ell,\lb}&=&\Hcal_{\ell-1,\lb}+\{\Hcal_{\ell-1,\lb},S_{\ell,\lb}\}+\int^{1}_{0}(1-s)\{\{\Hcal_{\ell-1,\lb},S_{\ell,\lb}\},S_{\ell,\lb}\}\circ\varphi_{\ell,\lb}^sds,\nonumber\\
 &=& \Hcal_{2,0}+\Ncal_{\ell-1,\lb}+P_{\ell,\lb}+\{\Hcal_{2,\lb},S_{\ell,\lb}\}+\Rcal_{\ell,\lb},\nonumber
\end{eqnarray}
where we obtain that $\Rcal_{\ell+1,\lb}(x)={\cal O}(|x|^{\ell+1})$ because $\varphi_{\ell,\lb}^s(x)=x+{\cal O}(|x|^{\ell-1})$ given that $S_{\ell,\lb}$ is a polynomial of degree $\ell$. So we get the following equation in the space of homogeneous polynomials of degree $\ell$ for $N_{\ell,\lb}$ and $S_{\ell,\lb}$
\begin{equation}\label{HomolEll}
N_{\ell,\lb}-\{\Hcal_{\ell,\lb},S_{\ell,\lb}\}=P_{\ell,\lb}.
\end{equation} 

\vspace{1ex}

\noindent\textbf{Step 2. Existence of the family of Lie transforms and smoothness in term of $\lb$.} Suppose that we constructed a $\C^1$-family of Hamiltonian polynomial functions $S_\lb$ as above, satisfying moreover $S_\lb(x)={\cal O}(\lb|x|^2+|x|^3)$. Denote by $\varphi_{\lb}(t,x)$ the associated family of Hamiltonian flows, $\ie$ the fixed points of the Picard operator
\begin{equation}\label{Picard}
\varphi_\lb(t,x)=x+\int^{t}_{0}J_m\nabla S_\lb(\varphi_\lb(s,x))ds:={\cal G}(\varphi_\lb,\lb),
\end{equation}
where ${\cal G}:{\cal A}(D_{\rho,2})\times\Lambda\rightarrow{\cal A}(D_{\rho,2})$ with $D_{\rho,2}=\{(t,x) / |x|\leq\rho,|t|\leq 2\}$. Given that $\nabla S_\lb={\cal O}(\lb|x|+|x|^2)$, for $\lb$ and $\rho$ sufficiently small the fixed point theorem applies and we get the existence of the family of Lie transforms $\varphi_\lb(1,.)$ in ${\cal A}(B_{\R^{2m}}(0,\rho))$. 

To get the class $C^1$ of the family we apply the Implicit Function Theorem to the implicit equation \eqref{Picard} in the neighborhood of $(\varphi_0,0)$. The linear part at $(\varphi_0,0)$ reads
$$h\mapsto h-\int^{t}_{0}J_m { }^t(D^2 S_\lb(\varphi_0(s,x)).h(s,x))ds$$
and is invertible for $\rho$ sufficiently small given that $D^2S_0(x)={\cal O}(|x|)$ and $\varphi_0(s,x)=x+{\cal O}(|x|^2)$ for $s\leq 2$, so that the assumptions of the Implicit Function Theorem are verified and hence the family $\varphi_\lb$ is $C^1$ for $\lb$ sufficiently small.

\vspace{1ex}

\noindent{\bf Step 3. Study of the homological operator.} Observe that the equations \eqref{Homol2} and \eqref{HomolEll} are implicit equations of the form
$${\cal F}_\ell(S_{\ell,\lb},N_{\ell,\lb},\lb)=0$$
whose linear part for $\lb=0$ is
$${\cal L} : E_\ell\times E_\ell \to E_\ell  : (N, S)\to N-\{\Hcal_{2,0},S\}:=N-{\cal A}S$$
where $E_\ell$ is the space of the real-valued homogeneous polynomials of degree $\ell$.

{\it We look for $F_\ell$ and $G_\ell$ two subspaces of $E_\ell$ such that ${\cal L} : F_\ell\times G_\ell \to E_\ell$ is invertible. Observe that $F_\ell$ is the space where we look for the nornal form monomials of degree $\ell$. So, our aim is to choose $F_\ell$ as "small" as possible and in particular we would like to have $F_\ell=\{0\}$ when it is possible.}

\vspace{1ex}

On one hand, {\bf when ${\cal A}$ is invertible} from $E_\ell$ onto $E_\ell$, then ${\cal L}$ is invertible from $\{0\}\times E_\ell\to E_\ell$. In this case we can choose $F_\ell=\{0\}$ and $G_\ell=E_\ell$. So, in this case we chose $N_{\ell,\lb}=0$ and solve \eqref{HomolEll} with the implicit function theorem to get $S_{\ell,\lb}$ as a $C^1$ function of $\lb$. So when ${\cal A}$ is invertible from $E_\ell$ onto $E_\ell$, all the \mbox{$\ell-$th}order term of the Hamiltonian are removed  by the normalization procedure.

On the other hand, {\bf when ${\cal A}$ is not invertible}, $F_\ell$ must be chosen as a supplementary space of the range ${\rm Im} {\cal A}$ of ${\cal A}$ and $G_\ell$ must be chosen as a supplementary space of $\ker {\cal A}$, $\ie$
$$   E_\ell=F_\ell\oplus {\rm Im}_{_{E_\ell}} {\cal A}, \qquad E_\ell=G_\ell\oplus\ker_{_{E_\ell}} {\cal A},$$ 
so that ${\cal L}:F_\ell\times G_\ell \to E_\ell$ is invertible. In this case, neither $N_{\ell,\lb}$ nor $S_{\ell,\lb}$ are unique. They both depend of the choice of $F_\ell$ and $G_\ell$. A natural way to choose these two subspaces is to endow $E_\ell$ with an inner product and to chose 
$$ F_\ell= ({\rm Im}_{_{E_\ell}} {\cal A})^\bot, \qquad G_\ell=(\ker_{_{E_\ell}} {\cal A})^\bot.$$
In what follows we define an appropriate inner product such that $F_\ell:=({\rm Im}{\cal A})^\bot=\ker{\cal A^*}$ is given by
$$  F_\ell = \{ S\in E_\ell / \ S(e^{tL_0^*}x)=S(x) \ \mbox{ for all } t\in \R, z\in \R^m\}=\{S\in E_\ell / \ \{\Hcal_{2,0}\circ(-J),S\}=0\}.$$
For that purpose, let us define, for any pairs of polynomials $S,S':\R^{2m}\to \R$ lying in $E_\ell$ the inner product given  by
$$  \scal{S}{S'}_\ell= S(\partial_{x}) S'(x)|_{x=0}.$$
Endowed with this inner product $E_\ell$ is a finite dimensional Hilbert space. Observe that for any integers $\alpha_1,\cdots, \alpha_{2m}$ and $\beta_1,\cdots, \beta_{2m}$
$$\scal{x_1^{\alpha_1}\cdot \! \dots \! \cdot x_{2m}^{\alpha_{2m}}}{x_1^{\beta_1}\cdot \! \dots \! \cdot x_{2m}^{\beta_{2m}}}_\ell=\alpha_1!\cdots\alpha_{2m}!\ \delta_{\alpha_1,\beta_1}\cdots\delta_{\alpha_{2m},\beta_{2m}}$$
where $\delta_{\alpha_j,\beta_j} = 1$ if $\alpha_j=\beta_j$ and $0$ otherwise. Moreover, for any invertible linear operator on $\R^{2m}$, considering the change of coordinates $x=T^*y$ (for which $\partial_y=T\partial_x$ and observing that $x=0\Leftrightarrow y=0$), we get
$$  \scal{S\circ T}{S'}_\ell=\scal{S}{S'\circ T^*}_\ell.$$
Hence for every $t\in\R$ we have $\scal{S\circ \E^{tL_0}}{S'}_\ell=\scal{S}{S'\circ e^{tL_0^*}}_\ell $. Given that  and differentiating this identity with respect to $t$ we finally get for $t=0$ that
\begin{equation}\label{DiffFlotLin}
\scal{{\cal A}S}{S'}=  \scal{D_{z}S(z).L_0z}{S'}=\scal{S}{D_{z}S(z).L_0^*z}=   \scal{S}{{\cal A}^*S'}.
\end{equation}
\begin{rem}
Observe that $\varphi_{\Hcal_{2,0}}^t(x)=\E^{tL_0}x$ is the Hamiltonian flow of $\Hcal_{2,0}$ and $\varphi_{\Hcal_{2,0}\circ(-J)}^t(x)=\E^{tL_0^*}x$ is the Hamiltonian flow of $\Hcal_{2,0}\circ(-J)$. Thus \eqref{DiffFlotLin} can be seen as a consequence of \eqref{ResLie}. 
\end{rem}
This ensures that the adjoint ${\cal A}^*: E_\ell \to E_\ell$ of the homological operator ${\cal A}$ is given by 
$$({\cal A}^*S)(z)= D_{z}S(z).L_0^*z=\{\Hcal_{2,0}\circ(-J),S\}$$ 
and that 
\begin{equation}\label{EqKerAstar}  
\ker_{_{E_\ell}} {\cal A}^*=\{ S\in E_\ell / \ S(e^{tL_0^*}z)=S(z) \ \mbox{for all } t\in \R, z\in \R^m\}=\{S\in E_\ell / \ \{\Hcal_{2,0}\circ(-J),S\}=0\}. 
\end{equation} 
Finally, let us chose as claimed above
 $$F_\ell:=\ker_{_{E_\ell}} {\cal A}^*=({\rm Im}_{_{E_\ell}} {\cal A})^\bot, \quad G_\ell=(\ker_{_{E_\ell}}{\cal A})^\bot.$$
Denote by $\pi_\ell$ the orthogonal projection onto $\ker_{_{E_\ell}} {\cal A}^*$. Then,  ${\cal L} :  (N, S)\mapsto N-{\cal A} S$ is an isomorphism from $F_\ell\times G_\ell$ onto $E_\ell$ since, for any ${\cal G}\in E_\ell$,
$${\cal L}(N,S)={\cal G}\Leftrightarrow  N-{\cal A} S={\cal G}\Leftrightarrow  \left \lbrace\begin{array}{rl} N & = \pi_\ell {\cal G}\\-{\cal A} S & = (Id-\pi_\ell) {\cal G}\end{array}\right.$$
where ${\cal A}$ is an isomorphism from $G_\ell=(\ker_{_{E_\ell}} {\cal A} )^\bot$ onto $(Id-\pi_\ell)E_\ell = {\rm Im}_{_{E_\ell}} {\cal A}$.

\vspace{1ex}

{\bf \noindent Step 4.1. Implicit Function Theorem at order 2.} The functional equation \eqref{Homol2} reads 
\begin{equation}\label{EqLG2}   
{\cal F}_2(N_{2,\lb},S_{2,\lb},\lb)= 0
\end{equation}
where ${\cal F}_2$ is a $C^1$ function from $F_2\times G_2 \times \Lambda$ to $E_2$ satisfying 
$${\cal F}_2(0,0,0)=0, \qquad D_{N,S}{\cal F}_2(0,0,0)={\cal L},$$
where $F_2:=\ker_{_{E_2}} {\cal A}^*$ and  $G_2:=(\ker_{_{E_2}} {\cal A})^\bot$. Hence, Step 3 and the Implicit Function Theorem ensures that for $N,S,\lb$ sufficiently small, (\ref{EqLG2}) has a unique solution of the form
$$(N_{2,\lb},S_{2,\lb}) = \psi_2(\lb)$$
where $\psi_2 : \Lambda \to F_2 \times G_2$ is $C^1$ and satisfies $\psi_2(0)=(0,0)$. So we get that for $\lb$ sufficiently small, $N_{2,\lb}$ and $S_{2,\lb}$ are $C^{1}$ functions of $\lb$ satisfying $N_{2,0}\equiv 0$ and $N_{2,\lb}$ lies in $F_2=\ker_{_{E_2}}{\cal A}^*$ given by (\ref{EqKerAstar}).

\vspace{1ex}

{\bf \noindent Step 4.2. Implicit Function Theorem for order $\ell\geq 3$.} We start with the equation of order $\ell$ given by \eqref{HomolEll}. Let us denote $F_\ell := \ker_{_{E_\ell}} {\cal A}^*$ and $G_\ell=(\ker_{_{E_\ell}} {\cal A})^\bot$ . Since $P_{\ell,\lb}$ and $\Hcal_{2,\lb}$ are $C^{1}$ functions of $\lb$, \eqref{HomolEll} reads
$${\cal F}_\ell(N_{\ell,\lb},S_{\ell,\lb},\lb)=0$$
where ${\cal F}_\ell$ is a $C^1$ function from $F_\ell\times G_\ell\times \Lambda$ to $E_\ell$. For $\lb=0$ it reads
$${\cal L}(N_{\ell,0},S_{\ell,0})=P_{\ell,0}.$$
Step 3 ensures that this equation has a unique solution $(N_{\ell,0},S_{\ell,0})$ in $F_\ell\times G_\ell$. Moreover,
$$ D_{N,S}{\cal F}_{\ell}(N_{\ell,0},S_{\ell,0},0)={\cal L},$$
hence the Implicit Function Theorem applied at the point $(N_{\ell,0},S_{\ell,0},0)$ ensures that for $\lb$ sufficiently small, $N_{\ell,\lb},S_{\ell,\lb}$ are $C^{1}$ functions of $\lb$ and that $N_{\ell,\lb}$ lies in $F_\ell=\ker_{_{E_\ell}}{\cal A}^*$ given by (\ref{EqKerAstar}).
\cqfd

\section{Appendix. Definition of $\prec$ and technical lemmas}\label{tech}
In this appendix we define the relation $\prec$ on the set of formal power series, and state a few properties of this relation. The proofs are left to the reader (the complete proofs are in the Appendix B of \cite{these}).

\begin{defi} \label{DefPrec}\label{DefiMax}
Let $f$ and $g$ be two formal power series on $\CC$, with $d$ variables. We denote them
$$f(x_1,\cdots,x_d)=\sum\limits_{n\in\NN^d}a_n x_1^{n_1}\cdots x_d^{n_d}:=\sum_{n\in\NN^d}a_n x^n, \quad g(x_1,\cdots,x_d)=\sum\limits_{n\in\NN^d}b_n x^n.$$
We define $\prec$ by $f(x_1,\cdots,x_d) \prec g(x_1,\cdots,x_d)$ if and only if
$$\forall n\in\NN^d, \quad b_n\in \R^+ \quad\text{and}\quad |a_n|\leq b_n.$$
We define also
\begin{eqnarray}
\left|f\right|(x)&:=&\sum\limits_{n\in\NN^d}|a_n|x^n, \nonumber\\
\max\limits_{\prec}\{f,g\} (x)&:=&\sum_{n\in\NN^d}\max\{|a_n|,|b_n|\}x^n.\nonumber
\end{eqnarray}
\end{defi}

\begin{lem}\label{ProprPrec}
Let $f$ and $g$ be two formal power series from $\CC^d$ to $\CC$.
\begin{enumerate}
\item $f\prec g$ if and only if $|f|\prec g$.
\item If $f\prec g$ then
$$\forall i\leq d, \quad \partial_{x_i}f\prec \partial_{x_i}g.$$
\item Let $d'<d$ and $y=(y_1,\cdots,y_{d'})$. Consider $f$ and $g$ for
$$x=(x_1,\cdots,x_d)=(y_1,y_1,\cdots,y_1,\cdots,y_{d'},\cdots,y_{d'}),$$
that we denote $\tilde{f}(y)$ and $\tilde{g}(y)$. Then
$$ f(x)\prec g(x) \Rightarrow  \tilde{f}(y) \prec \tilde{g}(y).$$
\item When $d=1$. If $f$ is a convergent power series of order $n_0$, then there exist two positive constants $c$ and $\gamma$ such that
$$f(x)\prec c\Frac{x^{n_0}}{1-\gamma x} ;$$
and more precisely, if $f\in\An(\B_{\CC}(0,\rho),\CC)$ (see Definition \ref{DefAn}), then
$$f(x) \prec \frac{\NNorme{\An}{f}}{\rho^{n_0}}\frac{x^{n_0}}{1-\frac{1}{\rho}x}.$$
\item If $0\prec f(x) \prec g(x)$, then
$$\Frac{1}{1-f(x)}\prec\Frac{1}{1-g(x)}.$$
\item If $f\prec g$ and $g\in \An(\B_{\CC^d}(0,\rho),\CC)$, then $f\in\An(\B_{\CC^d}(0,\rho),\CC)$ and
$$\NNorme{\An}{f}\leq\NNorme{\An}{g}.$$
\end{enumerate}
\end{lem}

%
%
%
%
%
%
%
%
%
%

\begin{lem}\label{LemPrecTAF}
Let $F$ be a scalar formal power series of the variables $(x_1,\cdots,x_d)$, with positive coefficients, and $\Phi$, $\Psi$ be two vectorial formal power series
$$\Phi(x_1,\cdots,x_d)=(\phi_1,\cdots,\phi_d)(x_1,\cdots,x_d),\quad\Psi(x_1,\cdots,x_d)=(\psi_1,\cdots,\psi_d)(x_1,\cdots,x_d).$$
Then we have the upper bound
$$F(\Phi+\Psi)-F(\Phi)\prec |DF|(|\Phi|+|\Psi|).|\Psi|,$$
where we use the notation
$$|DF|(|\Phi|+|\Psi|).|\Psi|:=\sum\limits_{i=1}^d |\partial_{x_i}F|(|\Phi_1|+|\Psi_1|,\cdots,|\Phi_d|+|\Psi_d|).|\Psi_i|.$$
\end{lem}

\begin{lem}\label{LemPrecInv}
We consider a family of formal power series of the form
$$F_\nuu(Y):=\left(F_{\nuu,1}(y_1,\cdots,y_d),\cdots,F_{\nuu,d}(y_1,\cdots,y_d)\right)=(y_1,\cdots,y_d)+{\cal O}(|Y|^2)$$
such that $F_0$ is convergent and there exist $N\in\NN$ and a convergent power series $\M$ such that for all $j=1,\cdots,d$,
$$F_{\nuu,j}-F_{0,j}\prec\nuu(y_1+\cdots+y_d)^N\M(y_1+\cdots+y_d)$$
holds  ; we moreover suppose that $F_\nuu$ is invertible and denote
$$F_\nuu^{-1}(Y):=\left(F_{\nuu,1}^-,\cdots,F_{\nuu,d}^-\right)(Y).$$

Then there exists a convergent power series $\M_1$ such that for all $j=1,\cdots,d$ 
$$F_{\nuu,j}^- -F_{0,j}^-\prec\nuu(y_1+\cdots+y_d)^N\M_1(y_1+\cdots+y_d)$$
holds.
\end{lem}

\section{Appendix. Construction of a local canonical change of coordinates : proof of Proposition \ref{propFF}} \label{constr}
This appendix is devoted to the proof of the Proposition \ref{propFFcomplete} below, which is a more general and detailed version of Proposition \ref{propFF}.

This proposition gives some estimates on the dependence of a change of coordinates $\FFep$ in term of some parameters $\ep:=(\eps,\nuu,\muu)$. For fixed values of the parameter $\ep$ the existence of the change of coordinates $\FFep$ is already well-known, it is a theorem of Moser \cite{Moser2} together with a result of Russmann \cite{Russmann}. The main result here is that under some assumptions, {\it the singularity in $\eps$ in the quadratic part of the initial Hamiltonian does not affect the bounds $(i)-(ix)$} (they are all independent of $\eps$).

\vspace{1ex}

This proposition plays a {\it crucial part in the proof of Theorem \ref{TH}}, points $(i)$ and $(ii)$ being the pivotal results. These results allow in particular to get fine properties of many objects such as the center-stable and center-unstable manifolds or the energy level sets, given that these objects are very simple in the $\xiet$-coordinates. 

Indeed, these results play a big role in the proof of Proposition \ref{PropReta} : at each step of the proof, we use local properties, local being in the spatial sense or for small values of the parameter $\nuu$, and we have to verify that the local neighborhoods do not tends to an empty set or do not move too much when $\eps$ goes to 0. For instance, consider Lemma \ref{Lemgcsep} in which we claim that the center-stable manifold $\CS$ can be expressed locally as the graph of a map $\gcsep$. The results of Proposition \ref{propFFcomplete} ensures that this graph expression holds for $\eps,\nuu,\muu$ small, in a neighborhood (namely $\B(0,\dll)$) independent of $\eps, \nuu$ and $\muu$.

Knowing precisely these neighborhoods is important then to perform a perturbative method (namely KAM theorem with the perturbative parameter $\muu$, see Part \ref{SectionKAM}) in a fixed annulus. 

\begin{prop} \label{propFFcomplete}
Let us consider a family of real analytic Hamiltonians $\H$ in $\An(\B_{\R^4}(0,\rho),\R)$ of the form
$$\H(\qp,\ep)=-\all\ql\pl+\frac{\alz(\eps)}{\eps^2}(\qz^2+\pz^2)+\h(\ql,\pl)+\nuu\Rest'(\qp,\ep),$$
where $\ep=(\eps,\nuu,\muu)$. We suppose that $\alz$ is a continuous function verifying $\alz(0)\neq 0$ and that $\Rest(.,\ep)$ is a $C^1$-family of analytic functions, and
$$\h(\ql,\pl)={\cal O}(|(\ql,\pl)|^3), \quad  \Rest(\qp,\ep)={\cal O}(|\qp|^3).$$

Then, there exist $\epso$, $\roop$ and a family of canonical changes of coordinates 
$$\FFep=\left(\philep,\psilep,\phizep,\psizep\right)$$ 
defined for $|\ep|\leq\epso$ such that the Hamiltonian in the new coordinates $\xiet$ reads
\begin{equation}\label{DefHL}
\begin{array}{ll}
\H\left(\FFep\xiet,\ep\right)& =\HLep(\xil\etl,\xiz^2+\etz^2)\\
& =-\xil\etl+\frac{\om(\eps)}{2\eps^2}(\xiz^2+\etz^2)+{\cal O}(|(\xil\etl,\xiz^2+\etz^2)|^2),
\end{array}
\end{equation}
and such that for all $\ep$, $\FFep$ and $\FFep^{-1}$ belong to $\in\An(\B_{\R^4}(0,\roop),\R^4)$, and for $\ep:=(\eps,0,0)$, $\FF_{(\eps,0,0)}$ is independent of $\eps$ and is of the form 
\begin{equation}\label{EqFstaro}
\FF_{(\eps,0,0)}\xiet:=\FFo\xiet=(\philo(\xil,\etl),\psilo(\xil,\etl),\xiz,\etz).
\end{equation}

\vspace{1ex}

Moreover we have the following estimates (recall that $\prec$ is defined in Appendix \ref{tech}) : there exists a power series of one variable $\M$, convergent on $\B(0,4\roop)$ such that 
\begin{enumerate}
\item $(\FFep-\FFo)\xiet\prec\nuu(\xil\+\etl\+\xiz\+\etz)^2\M(\xil\+\etl\+\xiz\+\etz)$,
\item $\FFep\xiet-\xiet\prec(\xil\+\etl\+\xiz\+\etz)^2\M(\xil\+\etl\+\xiz\+\etz)$,
\item $\FFep^{-1}\qp-\qp\prec (\ql\+\pl\+\qz\+\pz)^2\M(\ql\+\pl\+\qz\+\pz)$.
\end{enumerate}
And there exists a real $\Mo$ satisfying
\begin{enumerate}
\setcounter{enumi}{3}
\item $|\FFep\xiet-\FFo\xiet|\leq\nuu\Mo$,
\item $|\FFep^{-1}\qp-\FFo^{-1}\qp|\leq\nuu\Mo$,
\item $|\FFep\xiet-\xiet|\leq \Mo|\xiet|^2$,
\item $|\FFep^{-1}\qp-\qp|\leq \Mo|\qp|^2$,
\item $|\phizep\xiet-\xiz|\leq \nuu\Mo|\xiet|^2$,
\item $|\psizep\xiet-\etz|\leq \nuu\Mo|\xiet|^2$.
\end{enumerate}
for all $\xiet$ in $\B(0,\roop)$ and all $\qp$ in $\B(0,\roop)$.
\end{prop}

\subsubsection*{Plan of the proof}

Proving this fair dependence in term of the parameters requires to perform again each step of the existence proofs of Moser and Russmann, to verify the effect of the singularity in term of $\eps$ at each step of the construction of $\FFep$. This is a very long and technical work given that $\FFep$ is firstly constructed as a formal power series, by induction on the coefficients, so that at each step we compute the estimates by the method of majorant series.

\begin{description}
\item[Section \ref{SubPassageComplexe} : complex variables.] We compute a change of coordinates which diagonalize the linear part of the vector field, but also leads to complex variables. We detail this classical change of coordinates in order to prove at the end of the proof (in Part \ref{RetourR}) that the final change of coordinate is real.
\item[Section \ref{FamilleMoser} : non canonical changes of coordinates.] We first consider the family of changes of coordinates $\Fep$, which are not canonical. For fixed values of $\ep$, they are the changes of coordinates constructed by Moser \cite{Moser2}. The construction of the final {\it canonical} changes of coordinates relies on the $\Fep$, and so do the estimates of Proposition \ref{propFFcomplete}. So, in Part \ref{FamilleMoser}, we compute estimates on the dependence of $\Fep$ with respect to $\ep$.
\item[Section \ref{IntroFstarep} : introduction of the canonical changes of coordinates.] We then introduce the family of {\it canonical} changes of coordinates $\Fstarep$. For fixed values of $\ep$, they nearly are the changes of coordinates constructed by Russmann \cite{Russmann}. We moreover prove in this part the particular form of $\Fstar_{(\eps,0,0)}$ (in order to get \eqref{EqFstaro}).
\item[Section \ref{EstimFstarep} : estimate $(i)$ on the canonical changes of coordinates.] We prove that the $\FFep$ satisfy an estimate of the form of $(i)$ of Proposition \ref{propFFcomplete}.
\item[Section \ref{RetourR} : back to $\R$ and proof of the final form of the Proposition.] Gathering the previous results and getting back in $\R$, we can finally prove the Proposition \ref{propFFcomplete}.

\end{description}

\subsection{Diagonalization / complex variables}\label{SubPassageComplexe}
\def\H{\mathbf{H}}
We first perform the following complex canonical change of coordinates 
{\small $$\left(\begin{matrix} \ql' \\ \pl' \\ \qz' \\ \pz' \end{matrix}\right) = \left(\begin{matrix} I & O \\ O & P \end{matrix}\right)\left( \begin{matrix} \ql \\ \pl \\ \qz \\ \pz \end{matrix}\right):=\mathcal{P}\left( \begin{matrix} \ql \\ \pl \\ \qz \\ \pz \end{matrix}\right), \quad \text{ where } P=\frac{1}{\sqrt{2}}\left(\begin{matrix} 1 & \I \\ \I & 1 \end{matrix}\right).$$}
Let us construct first the changes of coordinates $\Fstarep$ in $\CC^4$, and at the end (see subsection \ref{RetourR}) define $\FFep:=\mathcal{P}^{-1}\Fstarep\mathcal{P}$. The Hamiltonian in the new coordinates reads
\begin{eqnarray}
\H(\xy,\ep)=-\all\xl\yl-\I\frac{\alz(\eps)}{\eps^2}\xz\yz+\h(\xl,\yl)+\nuu\Rest(\xy,\ep).\label{DefHamConstr}
\end{eqnarray}
with $\all,\alz\in\R$. The associated Hamiltonian system reads
\begin{equation}
\left\{\begin{array}{rcl}
\dot{\xl}&=&\all\xl+\partial_{\yl}\h(\xl,\yl)+\nuu\partial_{\yl}\Rest(\xy,\ep) \\
\dot{\yl}&=&-\all\yl-\partial_{\xl}\h(\xl,\yl)-\nuu\partial_{\xl}\Rest(\xy,\ep) \\
\dot{\xz}&=&\I\frac{\alz(\eps)}{\eps^2}\xz+\nuu\partial_{\yz}\Rest(\xy,\ep)\\
\dot{\yz}&=&-\I\frac{\alz(\eps)}{\eps^2}\yz-\nuu\partial_{\xz}\Rest(\xy,\ep).
\end{array}\right.
\end{equation}
Let us denote it by
\begin{equation}\label{SystDiff}
\left\{\begin{array}{rcl}
\dot{\xl}&=&\all\xl+\flo(\xl,\yl)+\nuu\flep(\xy) \\
\dot{\yl}&=&-\all\yl+\glo(\xl,\yl)+\nuu\glep(\xy) \\
\dot{\xz}&=&\I\frac{\alz(\eps)}{\eps^2}\xz+\nuu\fzep(\xy)\\
\dot{\yz}&=&-\I\frac{\alz(\eps)}{\eps^2}\yz+\nuu\gzep(\xy).
\end{array}\right.
\end{equation}
Note that for $i=1,2$, the families $(\fiep)_{\ep}$ and $(\giep)_{\ep}$ are $\C^1$-families of the space $\left(\An(\B_{\CC^4}(0,\rho_1),\CC),\NNorme{\infty}{\cdot}\right)$, given that $\ep\mapsto\Rest(\cdot,\ep)$ is $\C^1$ on $\An(\B_{\CC^4}(0,\rho_1),\CC)$.

\subsection{First family $\F_{\protect\underline{\varepsilon}}$, not canonical}\label{FamilleMoser}

Let us firstly introduce the following notation :

\begin{defi}\label{NotCrochet}
Let $\phi\xy$ be a formal power series from $\CC^4$ to $\CC$, of the form
$$\phi\xy=\sum_{m,n\in\NN^2} c_{m,n}\xl^{m_1}\yl^{n_1}\xz^{m_2}\yz^{n_2}.$$
Let us define the "sub"-power series, from $\CC^2$ to $\CC$, made of the sum of all the monomials of $\phi$ satisfying $m_1=n_1$ and $m_2=n_2$,
$$\crochet{\phi}(\om_1,\om_2):=\sum_{n\in\NN^2} c_{n,n}\om_1^{n_1}\om_2^{n_2}.$$
\end{defi}

We first fix $\ep$ and state the following lemma
\begin{lem}\label{LemMoser}
Consider a Hamiltonian of the form (\ref{DefHamConstr}) with $\alpha_1,\alpha_2$ real satisfying the assumptions of Proposition \ref{propFFcomplete}. Then
\begin{description}
\item[$(i)$] for all $\ep$ there exist a formal power series
$$\xy=\tilde{\Fep}\xiet=(\phileptilde,\psileptilde,\phizeptilde,\psizeptilde)(\xiet),$$
and formal power series of 2 variables $\aieptilde$ and $\bieptilde$ for $i=1,2$ 
satisfying
$$\tilde{\Fep}(\x,\ep)=\x+{\cal O}(|\x|^2),$$  
\begin{equation}\nonumber
\begin{array}{rclrcl}
\aleptilde(\om_1,\om_2)&=&\all+{\cal O}((\om_1,\om_2)),& \quad\azeptilde(\om_1,\om_2)&=&\I\frac{\alz(\eps)}{\eps^2}+{\cal O}((\om_1,\om_2)),\\ \bleptilde(\om_1,\om_2)&=&-\all+{\cal O}((\om_1,\om_2)),&\quad \bzeptilde(\om_1,\om_2)&=&-\I\frac{\alz(\eps)}{\eps^2}+{\cal O}((\om_1,\om_2)),
\end{array}
\end{equation}
such that in the new coordinates $\xiet$, the hamiltonian system reads
\begin{equation} \label{SystDiffFin}
\left\{\begin{array}{rcl}
\dot{\xil}&=&\aleptilde(\xil\etl,\xiz\etz)\xil \\
\dot{\etl}&=&\bleptilde(\xil\etl,\xiz\etz)\etl \\ 
\dot{\xiz}&=&\azeptilde(\xil\etl,\xiz\etz)\xiz \\
\dot{\etz}&=&\bzeptilde(\xil\etl,\xiz\etz)\etz.
\end{array}\right.
\end{equation}
\item[$(ii)$] for all $\ep$, we have existence and uniqueness of the formal power series $\tilde{\Fep}$, $\aieptilde$ et $\bieptilde$ of $(i)$ if 
\begin{equation} \label{DegLiberte}
\crochet{\frac{\phileptilde}{\xil}}, \crochet{\frac{\psileptilde}{\etl}}, \crochet{\frac{\phizeptilde}{\xiz}}, \crochet{\frac{\psizeptilde}{\etz}}
\end{equation}
are any 4 given power series of the form $1+{\cal O}(\om_1,\om_2)$.
\item[$(iii)$] If $\tilde{\Fep}$ is a formal power series satisfying $(i)$, then $\tilde{\Fep}'$ also satisfies $(i)$ if and only there exist some formal power series $\Phi_i,\Psi_i$ such that
\begin{equation}\nonumber
\begin{array}{l}
\Phi_i(\om_1,\om_2)=1+{\cal O}((\om_1,\om_2)),\quad\Psi_i(\om_1,\om_2)=1+{\cal O}((\om_1,\om_2)),\\
\tilde{\Fep}'(\xi,\eta)=\tilde{\Fep}\left(\Phi_1(\xil\etl,\xiz\etz)\xil,\Psi_1(\xil\etl,\xiz\etz)\etl,\Phi_2(\xil\etl,\xiz\etz)\xiz,\Psi_2(\xil\etl,\xiz\etz)\etz\right).
\end{array}
\end{equation}
\item[$(iv)$] Necessarily, $\alep=-\blep$ et $\azep=-\bzep$.
\item[$(v)$] Moreover, if the 4 power series (\ref{DegLiberte}) are convergent power series, then $\tilde{\Fep}$, $\aieptilde$ and $\bieptilde$ are convergent, $\ie$ there exists a disk on which they are analytic.
\end{description}
\end{lem}

\textbf{Proof:} in the statement of this lemma, $(v)$ is in fact the direct consequence, in our case, of the {\it result stated by Moser \cite{Moser2}}. And $(i)-(iv)$ are the steps of his proof (still applied to our particular case) : $(i)$ and $(ii)$ correspond to his Part 2, $(iv)$ is the lemma stated in his Part 3 and $(iii)$ is the Step I of the proof of the latter lemma. Point $(v)$ is proved in his Part 4. We state here explicitly these steps of his proof because we use it below.

Morser's proof relies on the fact that system in the new coordinates is of the form (\ref{SystDiffFin}) if and only if $\tilde{\Fep}$ satisfies
\begin{equation}\label{Syst*}
\hspace{-1.3ex}\left\{\begin{array}{rcl}
\hspace{-1.5ex}\all\phileptilde\+\flo(\phileptilde,\psileptilde)\+\nuu\flep(\tilde{\F}_{\ep})\hspace{-1.8ex}&=&\hspace{-1.8ex}\left(\aleptilde(\xil\etl,\xiz\etz)D_1+\azeptilde(\xil\etl,\xiz\etz)D_2\right)\phileptilde \\
\hspace{-1.5ex}-\hspace{-0.3ex}\all\psileptilde\+\glo(\phileptilde,\psileptilde)\+\nuu\glep(\tilde{\F}_{\ep})\hspace{-1.8ex}&=&\hspace{-1.8ex}\left(\aleptilde(\xil\etl,\xiz\etz)D_1+\azeptilde(\xil\etl,\xiz\etz)D_2\right)\psileptilde \\
\hspace{-1.5ex}\I\frac{\alz(\eps)}{\eps^2}\phizeptilde+\nuu\flep(\tilde{\F}_{\ep})\hspace{-1.5ex}&=&\hspace{-1.5ex}\left(\aleptilde(\xil\etl,\xiz\etz)D_1+\azeptilde(\xil\etl,\xiz\etz)D_2\right)\phizeptilde \\
\hspace{-1.5ex}-\I\frac{\alz(\eps)}{\eps^2}\psizeptilde+\nuu\glep(\tilde{\F}_{\ep})\hspace{-1.5ex}&=&\hspace{-1.5ex}\left(\aleptilde(\xil\etl,\xiz\etz)D_1+\azeptilde(\xil\etl,\xiz\etz)D_2\right)\psizeptilde,
\end{array}\right.
\end{equation}
where we denote
\begin{equation}\label{DefD}
D_1:=\xil\partial_{\xil}-\etl\partial_{\etl}, \qquad D_2:=\xiz\partial_{\xiz}-\etz\partial_{\etz}.
\end{equation}
\cqfd

\begin{lem}\label{LemMajFep}
Let us denote by $\Fep=(\philep,\psilep,\phizep,\psizep)$, $\aiep$ the formal power series satisfying
\begin{equation}\label{DefFep}
\crochet{\frac{\philep}{\xil}}=1= \crochet{\frac{\psilep}{\etl}}=\crochet{\frac{\phizep}{\xiz}}=\crochet{\frac{\psizep}{\etz}},
\end{equation}
whose existence and uniqueness are asserted in Lemma \ref{LemMoser} $(ii)$. They have the following properties
\begin{description}
\item[$(i)$] $\F_{(\eps,0,0)}$ and $a_{1,(\eps,0,0)}$ are independent of $\eps$. Let us denote them by $a_{1,0}$ and $$\F_{(\eps,0,0)}:=\Fo=(\philo,\psilo,\phizo,\psizo),$$
Moreover, they read
\begin{eqnarray}
 &&\Fo\xiet=\left(\philo(\xil,\etl),\psilo(\xil,\etl),\xiz,\etz\right) \nonumber\\
 && a_{1,(\eps,0,0)}(\om_1,\om_2)=a_{1,0}(\om_1),\qquad a_{2,(\eps,0,0)}(\om_1,\om_2)=\I\frac{\alz(\eps)}{\eps^2}.\nonumber
\end{eqnarray}
\item[$(ii)$] There exist a convergent power series $\M$ such that for every $\ep\leq\eps_0$,
$$ \quad (\Fep-\Fo)\xiet \prec \nuu(\xil\+\etl\+\xiz\+\etz)^2\M(\xil\+\etl\+\xiz\+\etz).$$
\end{description}
\end{lem}

\textbf{Proof of $(i)$.} 
Let us consider System $(\ref{Syst*})$ with $\nuu=\muu=0$ and prove that there exists a solution $\F_{(\eps,0,0)}$ independent of $\eps$ and of the form required in $(i)$. Given that the uniqueness of such $\F_{(\eps,0,0)}$ and $a_{1,(\eps,0,0)}$ was established in Lemma \ref{LemMoser} $(ii)$, this will prove $(i)$.

With the system $(\ref{Syst*})$ together with the ansatz
\begin{eqnarray}
 &&\Fep\xiet=\left(\philep(\xil,\etl),\psilep(\xil,\etl),\xiz,\etz\right) \nonumber\\
 && \alep(\om_1,\om_2)=\alep(\om_1),\qquad \azep(\om_1,\om_2)=\I\frac{\alz(\eps)}{\eps^2},\nonumber
\end{eqnarray}
and using the particular form of the operators $D_1$ and $D_2$ (defined in (\ref{DefD})), we get
\begin{equation}
\hspace{-1.3ex}\left\{\begin{array}{rcl}
\hspace{-1.5ex}\all\philep(\xil,\etl)\+\flo(\philep(\xil,\etl),\psilep(\xil,\etl))\hspace{-1.5ex}&=&\hspace{-1.5ex}\alep(\xil\etl)D_1\philep(\xil,\etl) \\
\hspace{-1.5ex}-\hspace{-0.3ex}\all\psilep(\xil,\etl)\+\glo(\philep(\xil,\etl),\psilep(\xil,\etl))\hspace{-1.5ex}&=&\hspace{-1.5ex}\alep(\xil\etl)D_1\psilep(\xil,\etl) \\
\hspace{-1.5ex}\I\frac{\alz(\eps)}{\eps^2}\xiz&=&\I\frac{\alz(\eps)}{\eps^2}\xiz \\
\hspace{-1.5ex}-\I\frac{\alz(\eps)}{\eps^2}\etz&=&\I\frac{\alz(\eps)}{\eps^2}(-\etz).
\end{array}\right.
\end{equation}
Only the two first equations remain. To prove that there exist such $\philep,\psilep,\alep$ we moreover use $(iv)$ of Lemma \ref{LemMoser} : we are going to prove that there exist $\alep(\xil\etl)$ and $\blep(\xil\etl)$ solving the following system \ref{Syst*oReduit}, and we will conclude thank to $(iv)$ which asserts that $-\blep=\alep$.
\begin{equation}\label{Syst*oReduit}
\hspace{-1.3ex}\left\{\begin{array}{rcl}
\hspace{-1.5ex}\all\philep(\xil,\etl)\+\flo(\philep,\psilep)\hspace{-1.5ex}&=&\hspace{-1.5ex}\left(\alep(\xil\etl)\xil\partial_{\xil}+\blep(\xil\etl)\etl\partial_{\etl}\right)\philep(\xil,\etl) \\
\hspace{-1.5ex}-\hspace{-0.3ex}\all\psilep(\xil,\etl)\+\glo(\philep,\psilep)\hspace{-1.5ex}&=&\hspace{-1.5ex}\left(\alep(\xil\etl)\xil\partial_{\xil}+\blep(\xil\etl)\etl\partial_{\etl}\right)\psilep(\xil,\etl).
\end{array}\right.
\end{equation}
Our aim is to prove the existence of such formal power series $\philep,\psilep,\alep$ and $\blep$. Denote them by
\begin{eqnarray}
\philep(\xil,\etl):=\xil+\sum_{N=2}^{+\infty}\varphi_1^N(\xil,\etl), && \psilep(\xil,\etl):=\xil+\sum_{N=2}^{+\infty}\psi_1^N(\xil,\etl), \nonumber\\
\alep(\om_1):=\all+\sum_{N=1}^{+\infty}a_1^N(\om_1), &&  \blep(\om_1):=-\all+\sum_{N=1}^{+\infty}b_1^N(\om_1). \nonumber
\end{eqnarray}
where the $\varphi_1^N, \psi_1^N, a_1^N$ and $b_1^N$ are homogeneous polynomials of degree $N$. System (\ref{Syst*oReduit}) at degree $N\geq2$ read then
\begin{eqnarray}
\all\varphi_1^N-\all D_1\varphi_1^N-\xil a_1^{N-1}=\mathcal{F}_N (\varphi_1^M,\psi_1^M,a_1^{M-1}, M<N), \label{Syst*oRed1}\\
-\all\psi_1^N-\all D_1\psi_1^N-\etl b_1^{N-1}=\mathcal{G}_N (\varphi_1^M,\psi_1^M,a_1^{M-1}, M<N).\label{Syst*oRed2}
\end{eqnarray}
Observe that the kernel of the operator $(Id-D_1)$ is the vector space generated by the monomials of the form $\xil(\xil\etl)^{n_1}$, and is a bijection of the vector space generated by all the other monomials. We then prove the existence of the $\varphi_1^N, \psi_1^N, a_1^{N-1}$ et $b_1^{N-1}$ by induction: on one hand, for $\xil a_1^{N}(\xil\etl)$ we take the sum of all the monomials of the form $\xil(\xil\etl)^{n_1}$ in the polynomial $-\mathcal{F}_N (\varphi_1^M,\psi_1^M,a_1^{M-1}, M<N)$, and do the same for $b_1^{N}$. And on the other hand, for $\varphi_1^{N}(\xil,\etl)$ we chose the antecedent of
$$\mathcal{F}_N (\varphi_1^M,\psi_1^M,a_1^{M-1}, M<N)-\xil a_1^{N-1}(\xil\etl)$$
by the operator $(Id-D_1)$ in the vector space generated by the monomials which are note of the form $\xil(\xil\etl)^{n_1}$. We do the same with $\psi_1^{N}$, without monomials of the form $\etl(\xil\etl)^{n_1}$.

Then the solution $\Fep$ satisfies the properties stated $(i)$ ; indeed : on one hand $\crochet{\frac{\philep}{\xil}}=1= \crochet{\frac{\psilep}{\etl}}$ given that $\philep$ has only one monomial of the form $\xil(\xil\etl)^{n_1}$, which is $\xil$. And on the other hand $(\varphi_{1,(\eps,0,0)},\psi_{1,(\eps,0,0)}, a_{1,(\eps,0,0)}, b_{1,(\eps,0,0)})$ solves the system (\ref{Syst*oReduit}) which is independent of $\eps$. And by the uniqueness stated in Lemma \ref{LemMoser} $(ii)$, we get that $(\varphi_{1,(\eps,0,0)},\psi_{1,(\eps,0,0)}, a_{1,(\eps,0,0)}, b_{1,(\eps,0,0)})$ are independent of $\eps$.
\cqfd

\vspace{3ex}

\textbf{Proof of $(ii)$.}
Let us denote 
$$\Fep'=(\philep',\psilep',\phizep',\psizep'):=\frac{1}{\nuu}(\Fep-\Fo),\quad \aiep':=\frac{1}{\nuu}(\aiep-a_{i,0}).$$In the following proof, we use the {\it absolute value} and the {\it maximum} associated with the order relation $\prec$, which are introduced in Definition \ref{DefiMax}. We proceed in several steps.

\vspace{1ex}

\textbf{Step 1 : system (\ref{Syst*1''}) satisfied by the $\Fep'$.} $\Fep$ and $\Fo$ satisfy (\ref{Syst*}) (with $\ep$ for the former, and $(\eps,0,0)$) for the latter). We substract those 2 systems and divide by $\nuu$, and then use the explicit expression of $a_{2,0}$, $\philo, \psilo, \phizo, \psizo$, which implies in particular that 
$$D_2\philo=0=D_2\psilo=D_1\phizo=D_1\psizo.$$
We then get the following system satisfied by $\Fep'$
\begin{equation}\label{Syst*1}
\left\{\hspace{-1ex}\begin{array}{l}
\hspace{-0.5ex}\left(\all\moins\all D_1\moins\I\frac{\alz(\eps)}{\eps^2}D_2\right)\philep'\moins\xil\alep'\\
\hspace{13ex}=(a_{1,0}\moins\all)D_1\philep'\+\left(\alep'D_1\+\azep'D_2\right)(\philo\moins\xil\+\nuu\philep')\\
\hspace{15ex}\moins\flep(\Fo\+\nuu\Fep)\+\frac{1}{\nuu}\left(\flo(\philo,\psilo)\moins\flo(\philo\+\nuu\philep',\psilo\+\nuu\psilep')\right) \\
(\textit{the same with }\psi)\\
\hspace{-0.5ex}\left(\I\frac{\alz(\eps)}{\eps^2}\moins\all D_1\moins\I\frac{\alz(\eps)}{\eps^2}D_2\right)\phizep'\moins\xiz\azep'\\
\hspace{13ex}=(a_{1,0}\moins\all)D_1\phizep'\+\left(\alep'D_1\+\azep'D_2\right)\nuu\phizep'\moins\fzep(\Fo\+\nuu\Fep')\\
(\textit{the same with }\psi).
\end{array}\right.
\end{equation}

This system can be understood as an infinity of equalities of coefficients of power series, and then as equalities of the absolute values of coefficients. Moreover, observe that in the term of the left hand side, $\xil\alep'(\xil\etl,\xiz\etz)$ is made of monomials of the form $\xil(\xil\etl)^{n_1}(\xiz\etz)^{n_2}$ and $\left(\all(I\moins D_1)\moins\I\frac{\alz(\eps)}{\eps^2}D_2\right)\philep'$ does not have any monomial of the form $\xil(\xil\etl)^{n_1}(\xiz\etz)^{n_2}$.

Then, from each equality of (\ref{Syst*1}), we get two inequalities, and obtain the following system
{\small \begin{equation}\label{Syst*1''}
\left\{\begin{array}{l}
\left|\left(\hspace{-0.5ex}\all(I\moins D_1)\moins\I\frac{\alz(\eps)}{\eps^2}D_2\hspace{-0.5ex}\right)\philep'\right|\hspace{-0.5ex}\prec\hspace{-0.5ex}\bigg|(a_{1,0}\moins\all)D_1\philep'\+\left(\alep'D_1\+\azep'D_2\hspace{-0.5ex}\right)\hspace{-0.5ex}(\philo\moins\xil\+\nuu\philep')\\
\hspace{18ex}\left.\moins\flep(\Fo\+\nuu\Fep)\+\frac{1}{\nuu}\left(\flo(\philo,\psilo)\moins\flo(\philo\+\nuu\philep',\psilo\+\nuu\psilep')\right)\right|,\\ 
\left|\xil\alep'\right|\prec\bigg|(a_{1,0}\moins\all)D_1\philep'\+\left(\alep'D_1\+\azep'D_2\right)(\philo\moins\xil\+\nuu\philep')\\
\hspace{18ex}\left.\moins\flep(\Fo\+\nuu\Fep)\+\frac{1}{\nuu}\left(\flo(\philo,\psilo)\moins\flo(\philo\+\nuu\philep',\psilo\+\nuu\psilep')\right)\right|,\\
\text{the same with the 3 other equations of (\ref{Syst*1})}
\end{array}\right.
\end{equation}}

\vspace{2ex}

\textbf{Step 2 : an upper bound $\Mep$ of $\Fep'$.}
We introduce two families of formal power series 
\begin{eqnarray}
\Mep\xiet\hspace{-2ex}&:=&\hspace{-2ex}\max\limits_{\prec}\left\{\max\limits_{\prec}(D_1\philep',D_2\psilep')\+\left|\crochet{\philep'}\right|,\max\limits_{\prec}(D_1\psilep',D_2\psilep')+\left|\crochet{\psilep'}\right|\right.\nonumber\\
 & & \hspace{4ex}\left.\max\limits_{\prec}(D_1\phizep',D_2\psizep')\+\left|\crochet{\phizep'}\right|,\max\limits_{\prec}(D_1\psizep',D_2\psizep')+\left|\crochet{\psizep'}\right|\right\}\nonumber\\
\Aep(\om_1,\om_2)\hspace{-2ex}&:=&\hspace{-2ex}\max\limits_{\prec}\left\{\alep',\azep'\right\}\nonumber.
\end{eqnarray}

Let us prove that $\left|\philep'\right|\prec\Mep$. We denote
\begin{equation}\label{NotPhilep}
\philep'\xiet:=\sum\limits_{m,n\in\NN^2}c_{(m,n)}\xil^{m_1}\etl^{n_1}\xiz^{m_2}\etz^{n_2}.
\end{equation}
Then
\begin{eqnarray}
\hspace{-4ex}\max\limits_{\prec}(D_1\philep',D_2\psilep')\+\left|\crochet{\philep'}\right|\hspace{-2ex}&=&\hspace{-3.5ex}\sum\limits_{m,n\in\NN^2}\hspace{-2ex}|c_{(m,n)}\hspace{-0.5ex}|\max\left\{|m_1\moins n_1|,|m_2\moins n_2|,1\right\}\xil^{m_1}\etl^{n_1}\xiz^{m_2}\etz^{n_2}.\label{ValMep}\\
  &\prec& \sum\limits_{m,n\in\NN^2}|c_{(m,n)}\hspace{-0.5ex}|\xil^{m_1}\etl^{n_1}\xiz^{m_2}\etz^{n_2}=\left|\philep'\right|\xiet.\nonumber
\end{eqnarray}
So $\left|\philep'\right|\prec\Mep$. With the same method one can check that $\Mep$ is also an upper bound of $\left|\psilep\right|,\left|\phizep\right|$ and $\left|\psizep\right|$.


\vspace{2ex}

\textbf{Step 3 : lower bounds of the left hand side terms of (\ref{Syst*1''}) in term of $\Mep$.} Let us show that there is the following lower bound for the left hand part of equation (\ref{Syst*1''})
\begin{equation}\label{InegGauche}
\min\left\{\frac{|\all|}{4},\frac{\alz(\eps)}{2\eps^2}\right\}\cdot\left(\max\limits_{\prec}(D_1\philep',D_2\phizep')\+\left|\crochet{\philep'}\right|\right)\hspace{-0.5ex}\prec\hspace{-0.5ex}\left|\left(\all(I\moins D_1)\moins\I\frac{\alz(\eps)}{\eps^2}\right)\philep'\right|.
\end{equation}
Indeed, using notation (\ref{NotPhilep}), the term of the right hand side of (\ref{InegGauche}) satisfies
\begin{eqnarray}
&&\hspace{-6ex}\left|\left(\all(I-D_1)\moins\I\frac{\alz(\eps)}{\eps^2}\right)\philep'\right|\xiet \nonumber\\
&\succ & \sum\limits_{m,n\in\NN^2}\frac{|c_{(m,n)}|}{2}\left(|\all(1\moins(m_1\moins n_1))|\+|\frac{\alz(\eps)}{\eps^2}(m_2\moins n_2)|\right)\xil^{m_1}\etl^{n_1}\xiz^{m_2}\etz^{n_2}\nonumber
\end{eqnarray}
given that $\all$, $\alz$ are real. And we already computed the term of the left hand side of (\ref{InegGauche}) (see (\ref{ValMep})). So, (\ref{InegGauche}) will be proved if we show that, for all $m,n\in\NN^2$,
\begin{eqnarray}
&&\hspace{-5ex}\min\left\{\frac{|\all|}{4},\frac{\alz(\eps)}{2\eps^2}\right\}|c_{(m,n)}|\max\left\{|m_1\moins n_1|,|m_2\moins n_2|,1\right\}\nonumber\\
&\leq &\frac{|c_{(m,n)}|}{2}\left(|\all(1\moins(m_1\moins n_1))|\+|\frac{\alz(\eps)}{\eps^2}(m_2\moins n_2)|\right).\nonumber
\end{eqnarray}
holds. We consider two different cases : firstly, if $m_1=n_1+1$ and $m_2=n_2$, then $c_{(m,n)}=0$. Indeed, from the definition of $\Fep$ we get
$$\crochet{\frac{\philep'}{\xil}}=\frac{1}{\nuu}\left(\crochet{\frac{\philep}{\xil}}\moins\crochet{\frac{\philo}{\xil}}\right)=\frac{1}{\nuu}(1-1)=0.$$
Secondly, if $m_1\neq n_1+1$ or $m_2\neq n_2$, we then get the result using that $|1-(m_1-n_1)|\geq1$ or $|m_2-n_2|\geq1$. Hence
$$ |\all||1-(m_1-n_1)|+\left|\frac{\alz(\eps)}{\eps^2}\right||m_2-n_2|\geq \min\{|\all|,\left|\frac{\alz(\eps)}{\eps^2}\right|\}\cdot 1.$$

The proof of (\ref{InegGauche}) is then completed. And with the same method, we obtain the same type of inequalities for $\psilep',\phizep',\psizep'$, and finally get
\begin{eqnarray}
\min\left\{\frac{|\all|}{4},\frac{\alz(\eps)}{2\eps^2}\right\}\Mep\hspace{-1.5ex}&\prec&\hspace{-1.5ex} \max\limits_{\prec}\left(\left|\left(\all(I\moins D_1)\moins\I\frac{\alz(\eps)}{\eps^2}\right)\phiiep'\right|,\left|\left(\all(I\moins D_1)\moins\I\frac{\alz(\eps)}{\eps^2}\right)\psiiep'\right|\right)\nonumber
\end{eqnarray}

\vspace{2ex}

\textbf{Step 4 : upper bounds of the right hand side terms of (\ref{Syst*1''}) in term of $\Mep$.} From now on, taking $\xiet=(\om,\om,\om,\om)$, we consider all the power series as power series of one variable $\om$. The inequalities $\prec$ are preserved by this change of variables. Given that $\philo,\psilo,\phizo,\psizo$ and $a_{1,0}-\all$ are convergent power series (by Lemma \ref{LemMoser}) without term of degree $0$, and that $\flo,\glo$ are convergent without terms of degree $0$ and $1$, Lemma \ref{ProprPrec} ensures the existence of $c,\gamma$ such that
\begin{eqnarray}
&&\varphi_{i,0}(\om,\om,\om,\om),\psi_{i,0}(\om,\om,\om,\om) \prec c\frac{\om}{1-\gamma\om}, \qquad a_{1,0}(\om^2)-\all\prec c\frac{\om^2}{1-\gamma\om},\label{MajFo}\\
&&\flo,\glo(\om,\om)\prec c\frac{\om^2}{1-\gamma\om}.\label{Majflo}
\end{eqnarray}
$\flep,\glep,\fzep$ et $\gzep$ are $\C^0$ in term of $\ep$ in the space $\left(\An(\B_{\CC^4}(0,\rho_1),\CC),\NNorme{\An}{\cdot}\right)$, so that there exists an uniform upper bound in term of ${\An}{\cdot}$ uniform en $\ep$. Moreover, the $\flep,\glep,\fzep$ and $\gzep$ do not have any monomial of degree $0$ or $1$. Then we get (with larger $c$ and $gamma$ if necessary)
\begin{equation}\label{Majflep}
\flep,\glep,\fzep,\gzep(\om,\om,\om,\om)\prec c\frac{\om^2}{1-\gamma\om}.
\end{equation}
Let us compute upper bounds, in term of $\Mep(\om,\om,\om,\om)$ and $\Aep(\om^2,\om^2)$, of the terms of the right hand side of (\ref{Syst*1''}). We use Lemma \ref{LemPrecTAF} to obtain the last of the following inequalities.
\begin{eqnarray}
\big|(a_{1,0}-\all)D_1\philep'\big|(\om,\om,\om,\om)
 &\prec& c\frac{\om^2}{1-\gamma\om}\Mep(\om,\om,\om,\om);\nonumber\\
\left|\left(\alep'D_1\+\azep'D_2\right)\hspace{-0.5ex}(\philo\moins\xil\+\nuu\philep')\right|&\prec& \Aep(\om^2,\om^2)\left(2\nuu\Mep(\om,\om,\om,\om)+c\frac{\om^2}{1-\gamma\om}\right); \nonumber\\
\left|\flep(\Fep(\om,\om,\om,\om))\right|&\prec&\hspace{-2ex}|\flep|\hspace{-0.5ex}\left(|\philo|\+\cdots\+|\psizo|\+\nuu(|\philep'|\+\cdots\+|\psizep'|),\cdots\right)\nonumber\\
 &\prec& c \frac{\left(4c\frac{\om}{1-\gamma\om}\+4\nuu\Mep(\om,\om,\om,\om)\right)^2}{1-\gamma\left(4c\frac{\om}{1-\gamma\om}+4\nuu\Mep(\om,\om,\om,\om)\right)};\nonumber\\
\frac{1}{\nuu}\left(\flo(\philo,\psilo)-\flo(\philep,\psilep)\right)&& \nonumber\\
&&\hspace{-12ex}\prec \frac{1}{\nuu}\left|D\flo\right|\left(|\philo|\+|\psilo|\+\nuu(|\philep'|\+|\psilep'|)\right)\cdot\nuu(|\philep'|\+|\psilep'|)\nonumber\\
&&\hspace{-12ex} \prec c\frac{2c\frac{\om}{1-\gamma\om}\+2\nuu\Mep(\om,\om,\om,\om)}{1-\gamma\left(2c\frac{\om}{1-\gamma\om}\+2\nuu\Mep(\om,\om,\om,\om)\right)}2\Mep(\om,\om,\om,\om). \nonumber
\end{eqnarray}

\vspace{2ex}

\textbf{Step 5 : an inequation independent of $\ep$ satisfied by $\Mep$.}
We come back to system (\ref{Syst*1''}) and use the lower bounds of the terms of the left hand side (Step 3) and the upper bounds of the right hand side (Step 4). Taking the maximum of the terms of the left hand side, we obtain
{\small \begin{equation}\label{SystFin}
\left\{\begin{array}{rcl}
\Mep\hspace{-1ex}&\prec& \hspace{-1ex}c\frac{\om^2}{1-\gamma\om}\Mep\+\Aep\left(\hspace{-0.5ex}2\nuu\Mep\+c\frac{\om^2}{1-\gamma\om}\hspace{-0.5ex}\right)\+c \Frac{\left(4c\frac{\om}{1-\gamma\om}\+4\nuu\Mep\right)^2}{1\moins\gamma\left(\hspace{-0.5ex}4c\frac{\om}{1-\gamma\om}\+4\nuu\Mep\right)}\+c\Frac{2c\frac{\om}{1-\gamma\om}\+2\nuu\Mep}{1\moins\gamma\left(\hspace{-0.5ex}2c\frac{\om}{1-\gamma\om}\+2\nuu\Mep\right)}2\Mep\\
\om\Aep\hspace{-1ex}&\prec&\hspace{-1ex} c\frac{\om^2}{1-\gamma\om}\Mep\+\Aep\left(\hspace{-0.5ex}2\nuu\Mep\+c\frac{\om^2}{1-\gamma\om}\hspace{-0.5ex}\right)\+c \Frac{\left(4c\frac{\om}{1-\gamma\om}\+4\nuu\Mep\right)^2}{1\moins\gamma\left(\hspace{-0.5ex}4c\frac{\om}{1-\gamma\om}\+4\nuu\Mep\right)}\+c\Frac{2c\frac{\om}{1-\gamma\om}\+2\nuu\Mep}{1\moins\gamma\left(\hspace{-0.5ex}2c\frac{\om}{1-\gamma\om}\+2\nuu\Mep\right)}2\Mep.\\
\end{array}\right.
\end{equation}}
We introduce
$$\Nep(\om):=\Mep(\om,\om,\om,\om)+\om\Aep(\om^2,\om^2).$$
Moreover, given that $\philep,\psilep,\phizep$ and $\psizep$ do not have any monomial of degree $0$ or $1$ and the $\aiep$ of degree $0$, necessarily, $\Nep$ read $\Nep=\om^2\Nep'(\om)$, where $\Nep'$ is a formal power series. We aim now at showing an upper bound of $\Nep'$ independent of $\ep$. And from system (\ref{SystFin}) we get that there exist $c'$ and $\gamma'$ such that for $\nuu\leq 1$ $\Nep'$ satisfies
\begin{equation}\label{MajNep'}
\Nep'\prec c'\frac{c'\om}{1-\gamma\om}\Nep'+2\om\Nep'^2+c'\Frac{\left(\frac{c}{1-\gamma\om}+\om\Nep'\right)\left(\frac{c}{1-\gamma\om}+2\om\Nep'\right)}{1\moins\gamma'\left(\frac{c\om}{1-\gamma\om}+\om^2\Nep'\right)}.
\end{equation}
$\Nep'$ satisfies then a functional inequality independent of $\ep$. 

\vspace{2ex}

\textbf{Step 6 : construction of a convergent majorant series $Z$.}
We aim now at constructing a convergent power series $Z$ satisfying
\begin{equation}\label{EqImpZ}
Z(\om)= c'\frac{c'\om}{1-\gamma\om}Z(\om)+2\om Z^2+c'\Frac{\left(\frac{c}{1-\gamma\om}+\om Z\right)\left(\frac{c}{1-\gamma\om}+2\om Z\right)}{1\moins\gamma'\left(\frac{c\om}{1-\gamma\om}+\om^2Z\right)},
\end{equation}
and such that $Z(0)\geq \Nep'(0)$. Indeed, from (\ref{MajNep'}), by induction on the coefficients of $\Nep$ and $Z$ we will then obtain 
$$\Nep'(\om)\prec Z(\om),$$
so that $Z$ will be the convergent upper bound independent of $\ep$ that we are looking for.

And by the analytic implicit functions theorem applied to equation (\ref{EqImpZ}) in the neighborhood of $\om=0$, we get the existence of a such $Z$ :  indeed, for $\om=0$, equation (\ref{EqImpZ}) read
$$Z-c c'=0,$$
where $\partial_Z(Z-c c')=1\neq0$ for all $Z$. And, up to a choice of larger $c$ and $c'$ in (\ref{MajNep'}) if necessary, $c c'\geq\Nep'(0)$ is satisfied because $\Nep'(0)$ is uniformly bounded in term of $\ep$ for $\ep\leq\epso$ given that all the analytic functions defining $\Nep'$ are $\C^0$ in term of $\ep$. So we  get the existence of a convergent power series $Z$ such that $\Nep'(\om)\prec Z(\om)$. Then,
$$\Mep(\om,\om,\om,\om)\prec \om^2Z(\om), \quad \Aep(\om^2,\om^2)\prec\om Z(\om)$$
hold, and thus
$$\Mep\xiet\prec\Mep(\xil\+\etl\+\xiz\+\etz,\cdots,\xil\+\etl\+\xiz\+\etz)\prec (\xil\+\etl\+\xiz\+\etz)^2 Z(\xil\+\etl\+\xiz\+\etz)$$
and the same for $\Aep$. Finally thank to Step 2, we get
\begin{eqnarray}
|\philep'|,|\psilep'|,|\phizep'| \text{ and }|\psizep'|\xiet&\prec& (\xil\+\etl\+\xiz\+\etz)^2 Z(\xil\+\etl\+\xiz\+\etz),\nonumber\\
|\aiep'|\xiet&\prec& (\xil\+\etl\+\xiz\+\etz) Z(\xil\+\etl\+\xiz\+\etz).\nonumber
\end{eqnarray}
\cqfd


\subsection{The canonical family $\Fstar_{\protect\underline{\varepsilon}}$ : construction, properties of $\Fstaro$}\label{IntroFstarep}

In this part, the first Lemma \ref{LemFstar} states the existence of a family of canonical changes of coordinates $\Fstarep$ verifying a criteria, and then Lemma \ref{LemFstaro} gives the main properties of $\Fstaro$. Thank to the criteria, in next section we will prove that the $\Fstarep$ are convergent and we will show that $\Fstarep-\Fstaro$ admits upper bounds of the type claimed in Proposition \ref{propFFcomplete}. The appropriate choice of the criteria will also allows to prove in Part \ref{RetourR} that the $\FFep$ are real. Our criteria (\ref{Q}) is very close to the one used by R¸ssmann \cite{Russmann} but slightly different in order to obtain a real change of coordinates.

\begin{lem}\label{LemFstar}
For all $\ep$,
\begin{description}
\item[$(i)$] There exists some canonical formal power series
\begin{equation}\label{CloseId}
\Ftildep\xiet=\xiet+{\cal O}(|\xiet|^2),
\end{equation}
and a Hamiltonian formal power series of 2 variables $\HLep$ such that
\begin{equation}\label{EqHL}
\H(\Ftildep(\xiet),\ep)=\HLep(\xil\etl,\xiz\etz)=-\all\xil\etl-\I\frac{\alz(\eps)}{\eps^2}\xiz\etz+ \cdots
\end{equation}
\item[$(ii)$] If $\tilde{\Fep}$ satisfies $(i)$, then $\tilde{\Fep}'$ also satisfies $(i)$ if and only if there exists a formal power series $\S(\om_1,\om_2)$ such that
$$\tilde{\Fep}'(\xi,\eta)=\tilde{\Fep}\left(\Phi_1(\xil\etl,\xiz\etz)\xil,\Psi_1(\xil\etl,\xiz\etz)\etl,\Phi_2(\xil\etl,\xiz\etz)\xiz,\Psi_2(\xil\etl,\xiz\etz)\etz\right)$$
holds, where we denote
$$\Phi_i(\om_1,\om_2):=\E^{\partial_{\om_i}S(\om_1,\om_2)},\quad\Psi_i(\om_1,\om_2):=\E^{-\partial_{\om_i}S(\om_1,\om_2)}.$$
\item[$(iii)$] There exists a unique canonical formal power series $\Fstarep=(\philepstar,\psilepstar,\phizepstar,\psizepstar)$ satisfying $(i)$ and verifying the criteria
\begin{equation}\label{Q}
\om_1\left(\crochet{\frac{\philepstar}{\xil}}(\om_1,\om_2)-\crochet{\frac{\psilepstar}{\etl}}(\om_1,\om_2)\right)+\I\om_2\left(\crochet{\frac{\phizepstar}{\xiz}}(\om_1,\om_2)-\crochet{\frac{\psizepstar}{\etz}}(\om_1,\om_2)\right)=0.
\end{equation}
\end{description}
\end{lem}
\begin{rem}
Taking the values $\alep(\om_1,\om_2)=\partial_{\om_1}\HLep(\om_1,\om_2)$ and $\azep(\om_1,\om_2)=\partial_{\om_2}\HLep(\om_1,\om_2),$
we see that $(i)$ of this lemma is a particular case of Lemma \ref{LemMoser}.
\end{rem}

\textbf{Proof of $(i)$.} We construct $\Ftildep$ with the aid of a generatrix function, $\ie$ we construct a formal power series $\Wep$ and define $\Ftildep$ by the implicit equation
$$\Ftildep(\partial_{\etl}\Wep(x,\eta),\etl,\partial_{\etz}\Wep(x,\eta),\etz)=(\xl,\partial_{\xl}\Wep(x,\eta),\xz,\partial_{\xz}\Wep(x,\eta)).$$
In order to obtain $\Ftildep$ of the form (\ref{CloseId}) we look for $\Wep$ of the form
$$\Wep(x,\eta)=\Wep(\xl,\etl,\xz,\etz)=\xl\etl+\xz\etz+{\cal O}(|(x,\eta)|^3).$$
Then, the formal power series $\Wep$ and $\HLep$ verify (\ref{EqHL}) if and only if
\begin{equation}\label{EqWep}
\H((\xl,\partial_{\xl}\Wep(x,\eta),\xz,\partial_{\xz}\Wep(x,\eta)),\ep)=\HLep(\partial_{\etl}\Wep(x,\eta),\etl,\partial_{\etz}\Wep(x,\eta),\etz),
\end{equation}
with $\HLep$ of the form
$$\HLep\xiet=\HLep(\xil\etl,\xiz\etz)=-\all\xil\etl-\I\frac{\om(\eps)}{\eps^2}+{\cal O}(|\xiet|^3).$$
Let us denote
$$\Wep(x,\eta):=\sum\limits_{N\geq 2}\Wep^N(x,\eta), \quad \H(x,y,\ep):=\sum\limits_{N\geq 2}\H_{\ep}^N(x,y), \quad \HLep(\xi,\eta):=\sum\limits_{N\geq 2}\HLep^N(\xi,\eta);$$
where $\Wep^N,\H_{\ep}^N,\HLep^N$ are homogeneous polynomials of degree $N$. Then (\ref{EqWep}) is satisfied if and only if for all $N\geq2$,
$$\all D_1\Wep^N(x,\eta)+\I\frac{\alz(\eps)}{\eps^2}D_2\Wep^N(x,\eta)-\HLep^N(x,\eta)=\mathcal{F}_N(\Wep^M,\HLep^M,M<N),$$
where $\mathcal{F}_N$ only depends of the $\Wep^M, \HLep^M$ for $M<N$ and where we denote $D_1=\xl\partial_{\xl}-\etl\partial_{\etl}$, and $D_2=\xz\partial_{\xz}-\etz\partial_{\etz}$.

Given that $\all$ and $\alz$ are real, observing the images of the monomials, we get that $(\all D_1+\I\frac{\alz(\eps)}{\eps^2}D_2)$ is a bijection of the set of formal power series without monomials of the form $(\xl\etl)^m(\xz\etz)^n$. Thus, we construct the $\Wep^N$ and $\HLep^N$ by induction, choosing at the step $N$ :
\begin{eqnarray}
\HLep^N(\xl\etl,\xz\etz)&\hspace{-1ex}=&\hspace{-1ex}\left[\mathcal{F}_N(\Wep^M,\HLep^M,M<N)\right](\xil\etl,\xiz\etz) \nonumber\\
\Wep^N(\xl,\etl,\xz,\etz)&\hspace{-1ex}=&\hspace{-1ex}(\all D_1+\I\frac{\alz(\eps)}{\eps^2}D_2)^{-1}\left(\mathcal{F}_N(\Wep^M,\HLep^M,M<N)-\HLep^N(\xl\etl,\xz\etz)\right).\nonumber
\end{eqnarray}

\vspace{1ex}

\textbf{Proof of $(ii)$. } This is nearly the result proved by R¸ssmann \cite{Russmann} in its part 3 : he shows  that the change of coordinates
$$\xiet\mapsto\left(\Phi_1(\xil\etl,\xiz\etz)\xil,\Psi_1(\xil\etl,\xiz\etz)\etl,\Phi_2(\xil\etl,\xiz\etz)\xiz,\Psi_2(\xil\etl,\xiz\etz)\etz\right)$$
is canonical if and only if $\Phi_i$ and $\Psi_i$ are of the form stated in $(ii)$. Using the result of \cite{Moser2} here above in $(iii)$ of Lemma \ref{LemMoser} we then achieve the proof of $(ii)$.

\vspace{1ex}

\textbf{Proof of $(iii)$. } Let us fix one $\Ftildep$ satisfying $(i)$. From $(ii)$, we get that $\Fstarep$ solves $(iii)$ if and only if
on one hand there exists a formal power series $\S(\om_1,\om_2)$ such that
\begin{equation}\label{S}
\hspace{-3ex}\Fstarep=\Ftildep\left(\E^{\partial_{\om_1}\S(\xil\etl,\xiz\etz)}\xil,\E^{-\partial_{\om_1}\S(\xil\etl,\xiz\etz)}\etl,\E^{\partial_{\om_2}\S(\xil\etl,\xiz\etz)}\xiz,\E^{-\partial_{\om_2}\S(\xil\etl,\xiz\etz)}\etz\right),
\end{equation}
and on the other hand, denoting $\Fstarep=(\philepstar,\psilepstar,\phizepstar,\philepstar)$, the criteria (\ref{Q}) is verified, $\ie$
\begin{equation}\nonumber
\om_1\left(\crochet{\frac{\philepstar}{\xil}}(\om_1,\om_2)-\crochet{\frac{\psilepstar}{\etl}}(\om_1,\om_2)\right)+\I\om_2\left(\crochet{\frac{\phizepstar}{\xiz}}(\om_1,\om_2)-\crochet{\frac{\psizepstar}{\etz}}(\om_1,\om_2)\right)=0.
\end{equation}
Let us show that for this $\Ftildep$ there exists a unique $\S$ such that $\Fstarep$ satisfies those two conditions (except $S(0,0)$ that does not need to be unique). Denoting $\Ftildep:=(\phileptilde,\psileptilde,\phizeptilde,\psizeptilde)$, from (\ref{S}) we get for $i=1,2$
\begin{eqnarray}
\crochet{\frac{\phiiepstar}{\xi_i}}\hspace{-0.5ex}(\om_1,\om_2)\hspace{-0.5ex}=\E^{\partial_{\om_i}\hspace{-0.5ex}\S(\om_1,\om_2)}\hspace{-0.5ex}\crochet{\frac{\phiieptilde}{\xi_i}}\hspace{-0.5ex}(\om_1,\om_2),&& \hspace{-1.5ex}
\crochet{\frac{\psiiepstar}{\eta_i}}\hspace{-0.5ex}(\om_1,\om_2)\hspace{-0.5ex}=\E^{-\partial_{\om_i}\hspace{-0.5ex}\S(\om_1,\om_2)}\hspace{-0.5ex}\crochet{\frac{\psiieptilde}{\eta_i}}\hspace{-0.5ex}(\om_1,\om_2),\nonumber
\end{eqnarray}
Recall that $\Ftildep\xiet=\xiet+{\cal O}(|\xiet|^2)$, thus $\crochet{\frac{\phiieptilde}{\xi_i}}=1+\cdots$ and $\crochet{\frac{\psiieptilde}{\eta_i}}=1+\cdots$. Then (\ref{Q}) holds if and only if
\begin{equation}\label{Q'}
\begin{array}{l}
-2(\om_1\partial_{\om_1}\+\I\om_2\partial_{\om_2})\S\\
=\om_1\left[(\E^{\partial_{\om_1}\hspace{-0.5ex}\S}\moins(1\+\partial_{\om_1}\hspace{-0.5ex}\S))\+\left(\hspace{-0.5ex}\crochet{\frac{\phileptilde}{\xil}}\moins1\hspace{-0.5ex}\right)\E^{\partial_{\om_1}\hspace{-0.5ex}\S}\moins(\E^{-\partial_{\om_1}\hspace{-0.5ex}\S}\moins(1\moins\partial_{\om_1}\hspace{-0.5ex}\S))\moins\left(\hspace{-0.5ex}\crochet{\frac{\psileptilde}{\etl}}\moins1\hspace{-0.5ex}\right)\E^{-\partial_{\om_1}\hspace{-0.5ex}\S}\right]\\
+\I\om_2\left[(\E^{\partial_{\om_2}\hspace{-0.5ex}\S}\moins(1\+\partial_{\om_2}\hspace{-0.5ex}\S))\+\left(\hspace{-0.5ex}\crochet{\frac{\phizeptilde}{\xiz}}\moins1\hspace{-0.5ex}\right)\E^{\partial_{\om_2}\hspace{-0.5ex}\S}\moins(\E^{-\partial_{\om_2}\hspace{-0.5ex}\S}\moins(1\moins\partial_{\om_2}\hspace{-0.5ex}\S))\moins\left(\hspace{-0.5ex}\crochet{\frac{\psizeptilde}{\etz}}\moins1\hspace{-0.5ex}\right)\E^{-\partial_{\om_2}\hspace{-0.5ex}\S}\right].
\end{array}
\end{equation}
Recall that $S(0,0)$ does not need to be unique to complete the proof. Let us choose $\S(0,0)=0$ and denote
$$\S(\om_1,\om_2):=\sum\limits_{N\geq1}\S^N(\om_1,\om_2),$$
where the $\S^N$ are homogeneous polynomials of degree $N$. Then (\ref{Q'}) reads
\begin{equation}\label{Qinduction}
\left\{\begin{array}{rcl}
-2(\om_1\partial_{\om_1}\+\I\om_2\partial_{\om_2})\S^1&=&0,\\
-2(\om_1\partial_{\om_1}\+\I\om_2\partial_{\om_2})\S^N&=&\mathcal{F}_N(\S^M,1\leq M<N)\quad\text{for }N\geq2.
\end{array}\right.
\end{equation}
Computing the images of the monomials through the operator $(\om_1\partial_{\om_1}\+\I\om_2\partial_{\om_2})$, we observe that it is invertible on the set of formal power series without monomial of degree 0. Then \ref{Qinduction} allows to construct $\S$ without monomial of degree 0 in a unique way.
\cqfd

\vspace{1ex}

The following Lemma gives a more precise description of $\Fstar_{(\eps,0,0)}$.

\begin{lem}\label{LemFstaro}
$\Fstar_{(\eps,0,0)}$ is independent of $\eps$ ; then we denote it $\Fstar_{(\eps,0,0)}=\Fstaro$. 
Moreover $\Fstaro$ is of the form
\begin{equation}\label{FormFstaro}
\Fstaro\xiet=\left(\philostar(\xil,\etl),\psilostar(\xil,\etl),\xiz,\etz\right).
\end{equation}
\end{lem}

\textbf{Proof. } We proceed in 2 main steps.
\def\ep{(\eps,0,0)}

\textbf{Step 1. Construction of a $\Ftildep$ independent of $\eps$, of the form (\ref{FormFstaro}).} 

Following the strategy of proof of $(i)$ of Lemma \ref{LemFstar}, it is sufficient to prove that we can construct $\Wep$ and $\HLep$ of the form
$$\Wep(x,\eta)=\xl\etl+\xz\etz+W'(\xl,\etl),\quad \HLep(\om_1,\om_2)=-\all\om_1-\I\frac{\alz(\eps)}{\eps^2}+\HL'(\om_1).$$
We use the particular form of $\H$ when $\nuu=0$ (see (\ref{DefHamConstr})). Denoting $D_1=\xl\partial_{\xl}-\etl\partial_{\etl}$, we obtain that $\Wep$ and $\HLep$ satisfy $(i)$ of Lemma \ref{LemFstar} if and only if
\begin{equation}\label{EqFtildo}
-\all D_1W'(\xl,\etl)=\HL'((\xl+\partial_{\etl}W'(\xl,\etl))\etl)-\h(\xl,\etl+\partial_{\xl}W'(\xl,\etl)).
\end{equation}
Note that this equation is independent of $\eps$. Let us show that we can construct such $W'$ and $\HL'$. Denote
$$W'(\xl,\etl)=\sum\limits_{N\geq3}W'^N(\xl,\etl),\qquad \HL'(\om_1)\sum\limits_{N\geq2}\HL'^N(\om_1),$$
where the $W'^N, \HL'^N$ are homogeneous polynomials of degree $N$. At degree $N$, (\ref{EqFtildo}) reads
$$-\all D_1W'^N(\xl,\etl)-\HL'^N(\xl\etl)=\mathcal{F}_N(W'^M, \HL'^M, M<N).$$
Recall that $\all\neq0$. As in the proof of $(i)$ of Lemma \ref{LemFstar}, $\all D_1$ in invertible on the set of formal power series without monomials of the form $(\xl\etl)^m$, and we get the existence proceeding by induction with the same strategy.

\vspace{1ex}

\textbf{Step 2. Construction of $\Fstarep$ independent of $\eps$ and of the form (\ref{FormFstaro}).}

Following the strategy of construction of the proof of $(iii)$ of Lemma \ref{LemFstar}, it is sufficient to prove that we can construct 
$$\Fstarep=\Ftildep\left(\E^{\partial_{\om_i}\S_{\ep}(\xil\etl,\xiz\etz)}\xi_i,\E^{-\partial_{\om_i}\S_{\ep}(\xil\etl,\xiz\etz)}\eta_i ; i=1,2\right),$$
with the $\Ftildep=\Ftild_0$ of Step 1 above, and 
$$\S_{\ep}(\om_1,\om_2)=\S_0(\om_1).$$
Given that $\Ftild_0$ is of the form (\ref{FormFstaro}), the criteria (\ref{Q}) reads
$$\om_1\left(\E^{\partial_{\om_1}\S_0(\om_1)}\crochet{\frac{\tilde{\varphi}_{1,0}}{\xil}}(\om_1)-\E^{-\partial_{\om_1}\S_0(\om_1)}\crochet{\frac{\tilde{\psi}_{1,0}}{\etl}}(\om_1)\right)=0$$
\begin{eqnarray}
\Leftrightarrow -2\om_1\partial_{\om_1}\S_0(\om_1)=\om_1&\hspace{-8ex}\left[(\E^{\partial_{\om_1}\hspace{-0.5ex}\S_0(\om_1)}\moins(1\+\partial_{\om_1}\hspace{-0.5ex}\S_0(\om_1)))\+\left(\hspace{-0.5ex}\crochet{\frac{\tilde{\varphi}_{1,0}}{\xil}}\moins1\hspace{-0.5ex}\right)\E^{\partial_{\om_1}\hspace{-0.5ex}\S_0(\om_1)}\right.\nonumber\\
&\left.-(\E^{-\partial_{\om_1}\hspace{-0.5ex}\S_0(\om_1)}\moins(1\moins\partial_{\om_1}\hspace{-0.5ex}\S_0(\om_1)))\moins\left(\hspace{-0.5ex}\crochet{\frac{\tilde{\psi}_{1,0}}{\etl}}\moins1\hspace{-0.5ex}\right)\E^{-\partial_{\om_1}\hspace{-0.5ex}\S_0(\om_1)}\right].\nonumber
\end{eqnarray}
This equation is independent of $\eps$ and admits a unique solution $\S_0(\om_1)$ : the proof is similar to that of $(iii)$ of Lemma $\ref{LemFstar}$. 
\cqfd

\def\ep{\underline{\eps}}
\subsection{Upper bound of $\F^*_{\protect\underline{\varepsilon}}-\Fstaro$}\label{EstimFstarep}

In this part, we state in the Lemma \ref{LemMajFstar} that $\Fstarep-\Fstaro$ admits an upper bound of the form of those of Proposition \ref{propFFcomplete}. But firstly, we prove the technical Lemma \ref{LemCompo} which will be useful to prove Lemma \ref{LemMajFstar}.


\begin{lem}\label{LemCompo}
Let $\PPhi$ be a vectorial formal power series of the form
$$\PPhi :\xiet\hspace{-0.5ex} \mapsto\hspace{-0.5ex}\left(\Phil(\xil\etl,\xiz\etz)\xil,\Psil(\xil\etl,\xiz\etz)\etl,\Phiz(\xil\etl,\xiz\etz)\xiz,\Psiz(\xil\etl,\xiz\etz)\etz\right),$$
where $\Phil,\Psil,\Phiz,\Psiz$ are formal power series satisfying
$$\Phil(0,0)=1=\Psil(0,0)=\Phiz(0,0)=\Psiz(0,0).$$
Then there exist some formal power series $\Phil^-,\Psil^-,\Phiz^-,\Psiz^-$ such that
$$\PPhi^{-1}\hspace{-1ex}:\hspace{-0.5ex}\xiet \hspace{-0.5ex}\mapsto\hspace{-0.5ex} \left(\Phil^-\hspace{-0.5ex}(\xil\etl,\xiz\etz)\xil,\Psil^-\hspace{-0.5ex}(\xil\etl,\xiz\etz)\etl,\Phiz^-\hspace{-0.5ex}(\xil\etl,\xiz\etz)\xiz,\Psiz^-\hspace{-0.5ex}(\xil\etl,\xiz\etz)\etz\right)\hspace{-0.5ex}.$$
\end{lem}

\textbf{Proof.}
Let us denote the inverse of $\PPhi$
$$\PPhi^{-1}:\xiet\mapsto\left(\Philt\xiet,\Psilt\xiet,\Phizt\xiet,\Psizt\xiet\right).$$
We want to prove that $\Philt$ reads $\xil\Phi^-(\xil\etl,\xiz\etz)$, $\ie$ that all the monomials of the formal power series $\Philt$ are of the form
$(\xil\etl)^{n_1}(\xiz\etz)^{n_2}\xil$. This is equivalent to show that
$$(\xil\partial_{\xil}-\etl\partial_{\etl})\Philt=\Philt, \quad\text{and}\quad (\xiz\partial_{\xiz}-\etz\partial_{\etz})\Philt=0.$$
Writing similar conditions on $\Psilt$, $\Phizt$, $\Psizt$, we get that Lemma \ref{LemCompo} holds if and only if $\PPhi$ satisfies
\begin{equation}\label{eq1PPhi-}
D\PPhi^{-1}\xiet\left(\begin{matrix}
											1&0&0&0\\
											0&-1&0&0\\
											0&0&0&0\\
											0&0&0&0
											\end{matrix}\right)\left(\begin{matrix}\xil\\ \etl\\ \xiz\\ \etz\end{matrix}\right)
=\left(\begin{matrix}
											1&0&0&0\\
											0&-1&0&0\\
											0&0&0&0\\
											0&0&0&0
											\end{matrix}\right)\PPhi^{-1}\xiet,
\end{equation}											
and
\begin{equation}\label{eq2PPhi-}
D\PPhi^{-1}\xiet\left(\begin{matrix}
											0&0&0&0\\
											0&0&0&0\\
											0&0&1&0\\
											0&0&0&-1
											\end{matrix}\right)\left(\begin{matrix}\xil\\ \etl\\ \xiz\\ \etz\end{matrix}\right)
=\left(\begin{matrix}
											0&0&0&0\\
											0&0&0&0\\
											0&0&1&0\\
											0&0&0&-1
											\end{matrix}\right)\PPhi^{-1}.
\end{equation}											
We know that $\PPhi$ verifies those two criteria (\ref{eq1PPhi-}) and (\ref{eq2PPhi-}). Let us prove then that $\PPhi^{-1}$ satisfies (\ref{eq1PPhi-}) (same proof for (\ref{eq2PPhi-})). From (\ref{eq1PPhi-}) verified by $\PPhi$ at $\xiet=\PPhi^{-1}(\xil',\etl',\xiz',\etz')$ and using $$D\PPhi^{-1}\xiet=(D\PPhi(\PPhi^{-1}\xiet))^{-1},$$ 
we get
$$D\PPhi\circ\PPhi^{-1}(\xil',\etl',\xiz',\etz')\left(\begin{matrix}
											1&0&0&0\\
											0&-1&0&0\\
											0&0&0&0\\
											0&0&0&0
											\end{matrix}\right)\PPhi^{-1}(\xil',\etl',\xiz',\etz')$$
$$=\left(\begin{matrix}
											1&0&0&0\\
											0&-1&0&0\\
											0&0&0&0\\
											0&0&0&0
											\end{matrix}\right)\PPhi\circ\PPhi^{-1}(\xil',\etl',\xiz',\etz').$$
And this implies the result (\ref{eq1PPhi-}).
\cqfd

\begin{lem}\label{LemMajFstar}
For all $\ep\in]-\epso,\epso[$, the formal power series $\Fstarep$ is a convergent series.
Moreover there exists a convergent power series $\M$ such that
\begin{equation}\label{MajFstar}
\Fstarep-\Fstaro\prec\nuu(\xil\+\etl\+\xiz\+\etz)^2\M(\xil\+\etl\+\xiz\+\etz).
\end{equation}
\end{lem}

\textbf{Proof.} To show that the $\Fstarep$ are convergent, the strategy is the same as that of R¸ssmann \cite{Russmann}, even if here we use a criterium (\ref{Q}) different from his. {\it Here we only prove the upper bound (\ref{MajFstar})}.  

The first key ideas of this proof are the similar upper bound proved for $\Fep$ in Lemma \ref{LemMajFep} $(ii)$ and the result $(iii)$ of Lemma \ref{LemMoser} which allows to express $\Fstarep$ in term of $\Fep$ : there exist a vectorial formal power series 
{\small \begin{equation}\nonumber
\PPhiep : \xiet \mapsto \left(\Philep(\xil\etl,\xiz\etz)\xil,\Psilep(\xil\etl,\xiz\etz)\etl,\Phizep(\xil\etl,\xiz\etz)\xiz,\Psizep(\xil\etl,\xiz\etz)\etz\right)
\end{equation}}
where $\Phiiep(\om_1,\om_2)=1+\cdots$ and $\Psiiep(\om_1,\om_2)=1+\cdots,$ such that
$$\Fep=\Fstarep\circ\PPhiep.$$
Then, in this proof we will split $\PPhiep$ into a product of maps of the form
\begin{eqnarray}
\Fep\hspace{-0.6ex}&=&\Fstarep\hspace{2ex}\circ\hspace{5ex}\overbrace{\left((\xii,\eti)\mapsto(\Phiieptild\xii,\Psiieptild\eti)\right)}^{:=\PPhieptild}\hspace{5ex}\circ\hspace{2ex}\overbrace{\left((\xii,\eti)\hspace{-0.6ex}\mapsto\hspace{-0.6ex}(\Phiiep\Psiiep\xii,\eti)\right)}^{:=\h_{\ep}}\nonumber\\
&=&\hspace{-0.6ex}\underbrace{\Fstarep\hspace{-0.6ex}\circ\hspace{-0.6ex}\left(\hspace{-0.6ex}(\xii,\eti)\hspace{-0.6ex}\mapsto\hspace{-0.6ex}(\E^{\partial_{\om_i}\S}\xii,\E^{\moins\partial_{\om_i}\S}\eti)\hspace{-0.6ex}\right)}_{:=(\phileptch,\psileptch,\phizeptch,\psizeptch)}\hspace{-0.6ex}
\circ\hspace{-0.6ex}\underbrace{\left(\hspace{-0.6ex}(\xii,\eti)\hspace{-0.6ex}\mapsto\hspace{-0.6ex}(\E^{\Schapi}\xii,\E^{\moins\Schapi}\eti)\hspace{-0.6ex}\right)}_{:=\fep}\hspace{-0.6ex}
\circ\hspace{-0.6ex}\left((\xii,\eti)\hspace{-0.6ex}\mapsto\hspace{-0.6ex}(\Phiiep\Psiiep\xii,\eti)\right),\nonumber
\end{eqnarray}
where $\Sep, \Schapi, \Phiieptild,\Psiieptild, \Phiiep, \Psiiep$ are power series of two variables $(\om_1,\om_2)=(\xil\etl,\xiz\etz).$ 

\begin{rem}\label{RemRussman}
This factorization of $\Fstarep$ is the same as that used by R¸ssman \cite{Russmann} in his proof : he proves that $\Fstarep$ converges by showing that every map of the factorization converges. In our proof below, we recall how these maps are constructed, but \textbf{we consider that their convergence for all fixed $\ep$ is a known fact}.
\end{rem}
Here is an outline of the proof :
\begin{description}
\item[Step 1.] We define $\h_{\ep}$ and compute an estimate of $\h_{\ep}-\h_0$ ; 
\item[Step 2.] We introduce $\PPhieptild$ and prove that it reads
$$\PPhieptild: (\xii,\eti)\mapsto(\E^{\S_{i,\ep}(\xil\etl,\xiz\etz)}\xii,\E^{-\S_{i,\ep}(\xil\etl,\xiz\etz)}\eti) ;$$
\item[Step 3.] We show that the $\S_{i\ep}$ read 
$$\S_{i,\ep}=\Schapi+\partial_{\om_i}\Sep,$$
and prove an upper bound of $\fep-f_0$ ;
\item[Step 4.] We introduce $(\phileptch,\psileptch,\phizeptch,\psizeptch)$ and show an estimate of $\Sep-\S_0$.
\end{description}
The key argument of Steps 1, 2 and 3 is that $\Fstarep$ is symplectic. Step 4 relies on the criterium satisfied by $\Fstarep$. Let us now detail these steps.

\vspace{2ex}

\textbf{Step 1. } Let us define 
$$\h_{\ep}:\xiet\mapsto\left((\Philep\hspace{-0.5ex}\cdot\hspace{-0.5ex}\Psilep)(\xil\etl,\xiz\etz)\xil,\hspace{0.5ex}\etl,\hspace{0.5ex}(\Phizep\hspace{-0.5ex}\cdot\hspace{-0.5ex}\Psizep)(\xil\etl,\xiz\etz)\xiz,\hspace{0.5ex}\etz\right)\hspace{-0.5ex}.$$
Given that $\Fstarep$ is symplectic, , $ ^tD\Fstarep \Om D\Fstarep=\Om$ holds and thus
$$^tD\Fep \Om D\Fep= {}^tD\PPhiep {}^tD\Fstarep\Om D\Fstarep D\PPhiep= {}^tD\PPhiep\Om D\PPhiep.$$
Considering two of the coefficients of these matrix, we get
\begin{equation}\label{Eq*1}
\begin{array}{rcl}
\partial_{\etl}\philep\partial_{\xil}\psilep\moins\partial_{\etl}\psilep\partial_{\xil}\philep\+\partial_{\etl}\phizep\partial_{\xil}\psizep\moins\partial_{\etl}\psizep\partial_{\xil}\phizep\hspace{-1ex}&=&\hspace{-1ex}-\partial_{\om_1}\hspace{-0.5ex}\left(\Philep\Psilep\om_1\right)\\
\partial_{\etz}\philep\partial_{\xiz}\psilep\moins\partial_{\etz}\psilep\partial_{\xiz}\philep\+\partial_{\etz}\phizep\partial_{\xiz}\psizep\moins\partial_{\etz}\psizep\partial_{\xiz}\phizep\hspace{-1ex}&=&\hspace{-1ex}-\partial_{\om_2}\hspace{-0.5ex}\left(\Phizep\Psizep\om_2\right),
\end{array}
\end{equation}
with $\Fep:=(\philep,\psilep,\phizep,\psizep)$ and where $\Phiiep,\Psiiep$ are formal power series of two variables $(\om_1,\om_2)=(\xil\etl,\xiz\etz)$. 

On one hand, from the result $(ii)$ of Lemma \ref{LemMajFep} we get the existence of a convergent power series $\M_1$ such that
$$\partial_{\xil}(\Fep-\Fo)\xiet \prec \nuu (\xil\+\etl\+\xiz\+\etz)\M_1(\xil\+\etl\+\xiz\+\etz),$$
and obtain similar upper bounds for $\partial_{\xiz}(\Fep-\Fo)$ and for the $\partial_{\eti}(\Fep-\Fo)$. On the other hand, Lemma \ref{LemMoser}$(v)$ ensures that $\Fo$ converges. 

Thus, substracting (\ref{Eq*1}) for $\ep$  from (\ref{Eq*1}) for $\ep=(\eps,0,0)$, we get the existence of a convergent power series $\M_2$ such that
$$\partial_{\om_i}(\Phiiep\Psiiep\om_i)-\partial_{\om_i}(\Phiio\Psiio\om_i)\prec\nuu\M_2(\om_1\+\om_2).$$
From which we finally obtain the existence of a convergent power series $\M_{\h}$ such that $$\h_{\ep}-\h_0\prec\nuu(\xil\+\etl\+\xiz\+\etz)\M_{\h}(\xil\+\etl\+\xiz\+\etz).$$

\vspace{2ex}

\textbf{Step 2.} Let us define
$$\PPhieptild:=\PPhiep\circ\h_{\ep}^{-1}.$$
Given that $\PPhiep$ and $\h_{\ep}$ are both of the form
$$(\xii,\eti)\mapsto\left(f_i(\xil\etl,\xiz\etz)\xii,g_i(\xil\etl,\xiz\etz)\eti\right) \quad \text{with }f_i(\om_1,\om_2)=1+\cdots, g_i(\om_1,\om_2)=1+\cdots,$$
Lemma \ref{LemCompo} ensures that $\PPhieptild$ is also of this form, $\ie$ there exist $\Phiieptild=1+\cdots,\Psiieptild=1+\cdots$ such that 
$$\PPhieptild : (\xii,\eti)\mapsto\left(\Phiieptild(\xil\etl,\xiz\etz)\xii,\Psiieptild(\xil\etl,\xiz\etz)\eti\right).$$

{\bf Let us show that for $i=1,2$, $\Phiieptild$ and $\Psiieptild$ satisfy  
\begin{equation}\label{EqPhitild}
\Phiieptild(\om_1,\om_2)\Psiieptild(\om_1,\om_2)=1.
\end{equation}}
From the definition $\PPhieptild\circ\h_{\ep}=\PPhiep$ of $\PPhieptild$, we get that for $i=1,2$
\begin{equation}
\left\{\begin{array}{l} \Phiieptild\left(\Philep\Psilep(\om_1,\om_2)\om_1,\Phizep\Psizep(\om_1,\om_2)\om_2\right)\Phiiep\Psiiep(\om_1,\om_2)=\Phiiep(\om_1,\om_2)\\
\Psiieptild\left(\Philep\Psilep(\om_1,\om_2)\om_1,\Phizep\Psizep(\om_1,\om_2)\om_2\right) =\Psiiep(\om_1,\om_2).
\end{array} \right.\nonumber
\end{equation}
and then, for $i=1,2$
$$\Phiieptild\left(\Philep\Psilep\om_1,\Phizep\Psizep\om_2\right)\Psiieptild\left(\Philep\Psilep\om_1,\Phizep\Psizep\om_2\right)\Phiiep\Psiiep\om_i=\Phiiep\Psiiep\om_i.$$
Given that $\Philep\Psilep(\om_1,\om_2)=1+\cdots$ and $\Phizep\Psizep(\om_1,\om_2)=1+\cdots,$ the map
$$(\om_1,\om_2)\mapsto(\Philep\Psilep\om_1,\Phizep\Psizep\om_2)$$
is invertible, and we then obtain 
$$\Phiieptild(\om_1,\om_2)\Psiieptild(\om_1,\om_2)\om_i=\om_i,$$
which achieves the proof of (\ref{EqPhitild}). Introducing $\Siep:=\ln(\Phiieptild)$, $\Phiieptild$ reads
$$\Phiieptild(\om_1,\om_2)=\E^{\Siep(\om_1,\om_2)},\qquad \Psiieptild(\om_1,\om_2)=\E^{-\Siep(\om_1,\om_2)}.$$

\vspace{2ex}

\textbf{Step 3. } Given that $\Fstarep$ is symplectic and $\Fep\circ\h_{\ep}^{-1}=\Fstarep\circ\PPhieptild,$ we get
\begin{equation}\label{EgMatr}
^tD(\Fep\circ\h_{\ep}^{-1})\Om D(\Fep\circ\h_{\ep}^{-1})=^tD\PPhieptild \Om D\PPhieptild.
\end{equation}
Considering the coefficients of this matrix and denoting
$$\Ftildep=\left(\phileptilde,\psileptilde,\phizeptilde,\psizeptilde\right)=\Fep\circ\h_{\ep}^{-1}$$
we obtain
\begin{equation}\label{Eqstar1}
\begin{array}{rcl}
\partial_{\xiz}\phileptilde\partial_{\xil}\psileptilde\moins\partial_{\xiz}\psileptilde\partial_{\xil}\phileptilde\+\partial_{\xiz}\phizeptilde\partial_{\xil}\psizeptilde\moins\partial_{\xiz}\psizeptilde\partial_{\xil}\phizeptilde\hspace{-2ex}&=&\hspace{-2ex}\etl\etz(\partial_{\om_2}\Slep\moins\partial_{\om_1}\Szep).
\end{array}
\end{equation}
From Step 1 above and $(ii)$ of Lemma \ref{LemMajFep}, together with Lemma \ref{LemPrecInv}, we get that $\Ftildep$ reads
$$\Ftildep=\Ftild_0+\nuu\Ftildep',$$
where $\Ftildep'$ has an convergent upper bound for $\prec$, uniform in term of $\ep$. Then Equation (\ref{Eqstar1}) is of the form
$$R_0(\om_1,\om_2)+\nuu R_{\ep}(\om_1,\om_2)=\partial_{\om_2}\Slep(\om_1,\om_2)\moins\partial_{\om_1}\Szep(\om_1,\om_2),$$
where $R_0$ is convergent and $R_{\ep}$ is uniformly bounded for $\prec$. Let us define then
$$\Schapl(\om_1,\om_2):=\int^{\om_2}_{0}(R_0+\nuu R_{\ep})(\om_1,t)dt, \qquad \Schapz(\om_1,\om_2)\equiv 0,$$
and

\vspace{-3ex}

$$\fep:\xiet\mapsto\left(\E^{\Schapl(\xil\etl,\xiz\etz)}\xil,\E^{-\Schapl(\xil\etl,\xiz\etz)}\etl,\E^{\Schapz(\xil\etl,\xiz\etz)}\xiz,\E^{-\Schapz(\xil\etl,\xiz\etz)}\etz\right).$$ 
From these definition, the $\Schapi$ and $\fep$ are also of the form
$$\Schapi:=\hat{S}_{i,0}+\nuu\hat{S'}_{i,\ep},\qquad \fep=f_0+\nuu \fep'$$
with the $\hat{S'}_{i,\ep}$ and $\fep'$ bounded for $\prec$ by a convergent power series independent of $\ep$.

\vspace{1ex}

Moreover, the $\S_{i,\ep}-\Schapi$ verify then
\begin{equation}\label{EqSep}
\partial_{\om_2}(\S_{1,\ep}-\Schapl)-\partial_{\om_1}(\S_{2,\ep}-\Schapz)=0.
\end{equation}
Let us define
$$\Sep(\om_1,\om_2):=\int^{\om_1}_{0}(\S_{1,\ep}-\Schapl)(t,\om_2)dt+\int^{\om_2}_{0}(\S_{2,\ep}-\Schapz)(0,t)dt ;$$
thank to (\ref{EqSep}), $\Sep$ satisfies
$$\partial_{\om_1}\Sep=\S_{1,\ep}-\Schapl, \qquad \partial_{\om_2}\Sep=\S_{2,\ep}-\Schapz.$$

\vspace{2ex}

\textbf{Step 4.} Gathering the results of Steps 1 to 3 and using Lemma \ref{LemPrecInv}, we get that $\Fstarep$ reads
\begin{eqnarray}
\Fstarep&=&\Fep\circ\h_{ep}^{-1}\circ\fep^{-1}\left(\E^{-\partial_{\om_1}\Sep}\xil,\E^{\partial_{\om_1}\Sep}\etl,\E^{-\partial_{\om_2}\Sep}\xiz,\E^{\partial_{\om_2}\Sep}\etz\right)\nonumber\\    &:=&(\phileptch,\psileptch,\phizeptch,\psizeptch)\left(\E^{-\partial_{\om_1}\Sep}\xil,\E^{\partial_{\om_1}\Sep}\etl,\E^{-\partial_{\om_2}\Sep}\xiz,\E^{\partial_{\om_2}\Sep}\etz\right)\label{DefPhileptch}
\end{eqnarray}
where
\begin{equation}\label{EqDecompPhitch}
(\phileptch,\psileptch,\phizeptch,\psizeptch)=(\philotch,\psilotch,\phizotch,\psizotch)+\nuu(\phileptch',\psileptch',\phizeptch',\psizeptch'),
\end{equation}
with $(\phileptch',\psileptch',\phizeptch',\psizeptch')$ have a convergent upper bound for $\prec$ independent of $\ep$. 

Let us now study $\Sep$. We proceed in several steps.

\vspace{1ex}

\textbf{Step 4.1. Consequence of Criterium (\ref{Q}) : equation satisfied by $\Sep':=\frac{1}{\nuu}(\Sep-\S_0)$.}
With the notation (\ref{DefPhileptch}), observe that $\Fstarep$ satisfies the criterium (\ref{Q}) if and only if
\begin{eqnarray}
 &&\om_1 \left(\E^{-\partial_{\om_1}\Sep}\crochet{\frac{\phileptch}{\xil}}(\om_1,\om_2)-\E^{\partial_{\om_1}\Sep}\crochet{\frac{\psileptch}{\etl}}(\om_1,\om_2)\right)\nonumber\\
 &+&\I\om_2 \left(\E^{-\partial_{\om_2}\Sep}\crochet{\frac{\phileptch}{\xil}}(\om_1,\om_2)-\E^{\partial_{\om_2}\Sep}\crochet{\frac{\psileptch}{\etl}}(\om_1,\om_2)\right)=0.\nonumber
\end{eqnarray}
The same calculations as those of the proof of Lemma \ref{LemFstar}$(iii)$ (to obtain equation (\ref{Q'})) lead to the equivalent equation
\begin{equation}\nonumber
\begin{array}{l}
-2(\om_1\partial_{\om_1}\+\I\om_2\partial_{\om_2})\Sep\\
=\om_1\left[(\E^{-\partial_{\om_1}\hspace{-0.5ex}\Sep}\moins(1\moins\partial_{\om_1}\hspace{-0.5ex}\Sep))\+\left(\hspace{-0.5ex}\crochet{\frac{\phileptch}{\xil}}\moins1\hspace{-0.5ex}\right)\E^{-\partial_{\om_1}\hspace{-0.5ex}\Sep}\moins(\E^{\partial_{\om_1}\hspace{-0.5ex}\Sep}\moins(1\+\partial_{\om_1}\hspace{-0.5ex}\Sep))\moins\left(\hspace{-0.5ex}\crochet{\frac{\psileptch}{\etl}}\moins1\hspace{-0.5ex}\right)\E^{\partial_{\om_1}\hspace{-0.5ex}\Sep}\right]\\
+\I\om_2\left[(\E^{-\partial_{\om_2}\hspace{-0.5ex}\Sep}\moins(1\moins\partial_{\om_2}\hspace{-0.5ex}\Sep))\+\left(\hspace{-0.5ex}\crochet{\frac{\phizeptch}{\xiz}}\moins1\hspace{-0.5ex}\right)\E^{-\partial_{\om_2}\hspace{-0.5ex}\Sep}\moins(\E^{\partial_{\om_2}\hspace{-0.5ex}\Sep}\moins(1\+\partial_{\om_2}\hspace{-0.5ex}\Sep))\moins\left(\hspace{-0.5ex}\crochet{\frac{\psizeptch}{\etz}}\moins1\hspace{-0.5ex}\right)\E^{\partial_{\om_2}\hspace{-0.5ex}\Sep}\right].
\end{array}
\end{equation}
Writing this latter equation for $\ep$ and for $\ep=(\eps,0,0)$, we get that $\Sep':=\frac{1}{\nuu}(\Sep-\S_0)$ satisfies
\begin{equation}\label{Q*}
\begin{array}{l}
\hspace{-4ex}-2(\om_1\partial_{\om_1}\+\I\om_2\partial_{\om_2})\Sep'\\
=\om_1\left[\frac{1}{\nuu}(\E^{-\partial_{\om_1}\hspace{-0.5ex}\Sep}\moins\E^{-\partial_{\om_1}\hspace{-0.5ex}\S_0})\+\partial_{\om_1}\hspace{-0.5ex}\Sep'\+\crochet{\frac{\phileptch'}{\xil}}\E^{-\partial_{\om_1}\hspace{-0.5ex}\Sep}\+\left(\crochet{\frac{\philotch}{\xil}}\moins1\right)\frac{1}{\nuu}(\E^{-\partial_{\om_1}\hspace{-0.5ex}\Sep}\moins\E^{-\partial_{\om_1}\hspace{-0.5ex}\S_0})\right.\\
\hspace{6ex}\left.-\frac{1}{\nuu}(\E^{\partial_{\om_1}\hspace{-0.5ex}\Sep}\moins\E^{\partial_{\om_1}\hspace{-0.5ex}\S_0})\+\partial_{\om_1}\hspace{-0.5ex}\Sep'\moins\crochet{\frac{\psileptch'}{\etl}}\E^{\partial_{\om_1}\hspace{-0.5ex}\Sep}\moins\left(\crochet{\frac{\psilotch}{\etl}}\moins1\right)\frac{1}{\nuu}(\E^{\partial_{\om_1}\hspace{-0.5ex}\Sep}\moins\E^{\partial_{\om_1}\hspace{-0.5ex}\S_0})\right]\\

\hspace{3ex}+\I\om_2\left[\frac{1}{\nuu}(\E^{-\partial_{\om_2}\hspace{-0.5ex}\Sep}\moins\E^{-\partial_{\om_2}\hspace{-0.5ex}\S_0})\+\partial_{\om_2}\hspace{-0.5ex}\Sep'\+\crochet{\frac{\phizeptch'}{\xiz}}\E^{-\partial_{\om_2}\hspace{-0.5ex}\Sep}\+\left(\crochet{\frac{\phizotch}{\xiz}}\moins1\right)\frac{1}{\nuu}(\E^{-\partial_{\om_2}\hspace{-0.5ex}\Sep}\moins\E^{-\partial_{\om_2}\hspace{-0.5ex}\S_0})\right.\\
\hspace{6ex}\left.-\frac{1}{\nuu}(\E^{\partial_{\om_2}\hspace{-0.5ex}\Sep}\moins\E^{\partial_{\om_2}\hspace{-0.5ex}\S_0})\+\partial_{\om_2}\hspace{-0.5ex}\Sep'\moins\crochet{\frac{\psizeptch'}{\etz}}\E^{\partial_{\om_2}\hspace{-0.5ex}\Sep}\moins\left(\crochet{\frac{\psizotch}{\etz}}\moins1\right)\frac{1}{\nuu}(\E^{\partial_{\om_2}\hspace{-0.5ex}\Sep}\moins\E^{\partial_{\om_2}\hspace{-0.5ex}\S_0})\right].
\end{array}
\end{equation}

In order to obtain a majorant series of the $\Sep'$, we define
$$V_{i,\ep}(\om_1,\om_2):=|\partial_{\om_i}\Sep|(\om_1,\om_2), \qquad V_{\ep}(\om):=V_{1,\ep}(\om,\om)+V_{2,\ep}(\om,\om),$$
and compute a majorant series of the $V_{\ep}$ with the aid of (\ref{Q*}).

\vspace{1ex}

\textbf{Step 4.2. Lower bound of the left hand term of (\ref{Q*}).}

Let us denote
$$\Sep:=\sum\limits_{n_1,n_2\in\NN}a_{n_1,n_2} \om_1^{n_1}\om_2^{n_2}.$$
Then
\begin{eqnarray}
\om_i V_{i,\ep}(\om_1,\om_2)&=&\sum\limits_{n_1,n_2\in\NN}n_i|a_{n_1,n_2}|\om_1^{n_1}\om_2^{n_2},\nonumber\\
|\om_1\partial_{\om_1}\Sep+\I\om_2\partial_{\om_2}\Sep|&=& \sum\limits_{n_1,n_2\in\NN} |n_1+\I n_2||a_{n_1,n_2}|\om_1^{n_1}\om_2^{n_2}.\nonumber
\end{eqnarray}
Given that $n_i\leq|n_1+\I\n_2|$, we get
\begin{equation}\label{EqMajVi}
\om_i V_{i,\ep}\prec |\om_1\partial_{\om_1}\Sep+\I\om_2\partial_{\om_2}\Sep|(\om_1,\om_2).
\end{equation}

\vspace{1ex}

\textbf{Step 4.3. Existence of an upper bound of the right hand term of (\ref{Q*})} of the form
\begin{equation}\label{FormeMaj}
\sum\limits_{i=1}^2 \om_i\left(G_i^1(\om_1,\om_2)+G_i^2(\om_1,\om_2,V_{1,\ep},V_{2,\ep})\right),
\end{equation}
where $G_i^1, G_i^2$ are convergent power series independent of $\ep$ and $G_i^1$ is an homogeneous polynomial of degree 1 and $G_i^2$ is of order more than 2. To show that, let us compute successively upper bounds of each right hand side term of (\ref{Q*}):

\vspace{1ex}

\begin{eqnarray}
&&\frac{1}{\nuu}(\E^{-\partial_{\om_i}\hspace{-0.5ex}\Sep}\moins\E^{-\partial_{\om_i}\hspace{-0.5ex}\S_0})-\frac{1}{\nuu}(\E^{-\partial_{\om_i}\hspace{-0.5ex}\Sep}\moins\E^{-\partial_{\om_i}\hspace{-0.5ex}\S_0})+2\partial_{\om_i}\hspace{-0.5ex}\Sep' \nonumber\\
&=& -\frac{2}{\nuu}(\sinh(\partial_{\om_i}\Sep)-\sinh(\partial_{\om_i}\Sep))+2\partial_{\om_i}\hspace{-0.5ex}\Sep'=-2\partial_{\om_i}\Sep'\sum\limits_{p=1}^{\infty}\frac{\sum\limits_{k=0}^{2p}(\partial_{\om_i}\Sep)^k(\partial_{\om_i}\S_0)^{2p-k}}{(2p+1)!}\nonumber\\
&\prec& 2 V_{i,\ep} \sum\limits_{p=1}^{\infty}\frac{(2p+1)(|\partial_{\om_i}\S_0|+\nuu V_{i,\ep})^{2p}}{(2p+1)!} \prec 2 V_{i,\ep} (|\partial_{\om_i}\S_0|+V_{i,\ep})\sum\limits_{p=0}^{\infty} \frac{(|\partial_{\om_i}\S_0|+V_{i,\ep})^{2p}}{(2p+2)!} ; \nonumber
\end{eqnarray}
(the last inequality holds for $\nuu\leq 1$), where the last upper bound is of the form $G_i^2$ given that $\S_0$ is convergent and independent of $\ep$.

\vspace{1ex}

Given that $\phileptch$ is of the form (\ref{EqDecompPhitch}),  there exists a convergent series $\M$ independent of $\ep$ such that
\begin{eqnarray}
&&\left|\crochet{\frac{\phileptch'}{\xil}}\E^{-\partial_{\om_1}\Sep}\right|(\om_1,\om_2)=\left|\crochet{\frac{\phileptch'}{\xil}}+\crochet{\frac{\phileptch'}{\xil}}(\E^{-\partial_{\om_1}\Sep}-1)\right|(\om_1,\om_2)\nonumber\\
 &\prec& (\om_1\+\om_2)\M(\om_1,\om_2)(1+(\E^{|\partial_{\om_1}\Sep|(\om_1,\om_2)}-1))\nonumber\\
&\prec&(\om_1\+\om_2)\M(0,0)\+(\om_1\+\om_2)(\M(\om_1,\om_2)\moins\M(0,0)\+\M(\om_1,\om_2)(\exp-1)\circ(|\partial_{\om_1}\S_0|\+V_{1,\ep})\nonumber
\end{eqnarray}
which is of the form (\ref{FormeMaj}) given that for all $\ep$, $|\partial_{\om_i}\Sep|(0,0)=0$ holds because
$$\partial_{\om_i}\Sep(0,0)=\S_{i,\ep}(0,0)-\Schapi(0,0)=\S_{i,\ep}(0,0),$$
where
$$\E^{\S_{i,\ep}(\om_1,\om_2)}=\PPhiep(\om_1,\om_2)=1+\cdots.$$
We get similarly an upper bound of $\crochet{\frac{\psileptch'}{\etl}}\E^{\partial_{\om_1}\Sep}$.

\vspace{1ex}

Finally,
\begin{eqnarray}
&&\left|\crochet{\frac{\philotch}{\xil}}-1\right|\frac{1}{\nuu}\left|\E^{-\partial_{\om_1}\Sep}-\E^{-\partial_{\om_1}\S_0}\right|\nonumber\\
&\prec& \left|\crochet{\frac{\philotch}{\xil}}\moins1\right| \frac{1}{\nuu} |\moins\partial_{\om_1}\Sep\moins(\moins\partial_{\om_1}\S_0)| \ \E^{|\partial_{\om_1}\Sep|+|\partial_{\om_1}\S_0|} \quad\text{by Lemma \ref{LemPrecTAF}}\nonumber\\
&\prec& \left|\crochet{\frac{\philotch}{\xil}}\moins1\right|V_{1,\ep} \ \E^{2|\partial_{\om_1}\S_0|+V_{1,\ep}} \quad\text{for }\nuu\leq 1,\nonumber
\end{eqnarray}
which is of the form $G_i^2$. 

\vspace{1ex}

\textbf{Step. Equation of construction of the majorant series $Z$.}

Gathering the results of Steps 4.1, 4.2 and 4.3, we obtain that $V_{\ep}$ satisfies
\begin{eqnarray}
&&\om V_{\ep}(\om)\prec \om\left(G_1^1(\om,\om)\+G_2^1(\om,\om)\+G_1^2(\om,\om,V_{\ep}(\om),V_{\ep}(\om))\+G_2^2(\om,\om,V_{\ep}(\om),V_{\ep}(\om))\right)\nonumber\\
&\Rightarrow & V_{\ep}(\om)\prec c\om+ c\Frac{(\om+V_{\ep}(\om))^2}{1-\gamma(\om\+V_{\ep}(\om))} \quad \text{by $(iv)$ of Lemma \ref{ProprPrec}}.\label{IneqVep}
\end{eqnarray}
\textit{Let us construct a convergent power series $Z(\om)$ such that
 $$Z(\om)= c\om+ c\Frac{(\om+Z(\om))^2}{1-\gamma(\om\+Z(\om))} ,$$
 and prove by induction on the coefficients of the series that $V_{\ep}\prec Z(\om)$.}
 
\vspace{1ex}

\textbf{Step 4.5. Initialization of the induction (if existence of $Z$).}

Thank to (\ref{IneqVep}) the proof by induction works if we prove the initialization. At degree 1 in $(\om_1,\om_2)$, (\ref{Q*}) reads
$$2(\om_1\partial_{\om_1}+\I\om_2\partial_{\om_2})(\Sep')_1=0$$
given that the right hand side (\ref{Q*}) is of the form (\ref{FormeMaj}). Thus the degree 1 term of $\Sep'$, $(\Sep')_1$ vanishes.
At degree 2 in $(\om_1,\om_2)$, (\ref{Q*}) gives
$$2(\om_1\partial_{\om_1}+\I\om_2\partial_{\om_2})(\Sep')_2=\sum\limits_{i=1}^2\om_i\left(\crochet{\frac{\phileptch'}{\xil}}_1-\crochet{\frac{\psileptch'}{\etl}}_1\right) ;$$
where the right hand side term is uniformly bounded. Thus, there exists $\alpha$ independent of $\ep$ such that
$$\om (V_{\ep})_1(\om) \prec \alpha\om^2.$$

\vspace{1ex}

\textbf{Step 4.6. Existence of $Z$.}
To satisfy the initialization assumption, we look for $Z$ satisfying
$$Z(\om)=c\om + c \Frac{(\om+Z(\om))^2}{1-\gamma (\om+Z(\om))}
\quad \text{and } Z(\om)=\om Z_1(\om), \quad \text{with } Z_1(0)\geq \alpha ;$$
$\ie$ we look for $Z_1$ such that
$$Z_1\cdot(1-\gamma\om(1+Z_1))-c(1-\gamma\om(1+Z_1))-c\om(1+Z_1)^2=0 \Longleftrightarrow  \mathcal{F}(Z_1,\om)=0. $$
This equation, independent of $\ep$, satisfies the assumptions of the analytic Implicit Functions Theorem in the neighborhood of $(\om,Z_1)=(0,\alpha)$. Thus we get the existence of the analytic function $Z_1$, $\ie$ a convergent power series in the neighborhood of $0$. 

\vspace{1ex}

Choosing $\Sep'(0,0)=0$ (possible because only the derivatives of $\Sep'$ are in the definition of $\Fstarep$), we obtain the upper bound
$$\Sep'(\om_1,\om_2)\prec \Sep'(0,0)+\om_1 V_{1,\ep}(\om_1,\om_2)+\om_2 V_{2,\ep}(\om_1,\om_2)\prec (\om_1\+\om_2) V_{\ep}(\om_1\+\om_2) \prec (\om_1\+\om_2) Z_1(\om_1\+\om_2),$$
where $Z_1$ is a convergent power series independent of $\ep$.

\vspace{2ex}

\textbf{Step 5. Conclusion of the Lemma's proof.}

$\Fstarep$ reads
$$\Fstarep=\Fep\circ\h_{\ep}^{-1}\circ\fep^{-1}\circ\left((\xii,\eti)\mapsto(\E^{-\partial_{\om_i}\Sep}\xii,\E^{\partial_{\om_i}\Sep}\eti))\right).$$
And thank to Lemma \ref{LemPrecInv} and to the results of the previous steps, this is a product of formal power series of the form $R_{\ep}=R_0+\nuu R_{\ep}'$ where $R_0$, $R_{\ep}$ are convergent and admit a convergent upper bound $\M$ independent of $\ep$. Thus $\Fstarep$ is also of this form $R_0+\nuu R_{\ep}'$. Moreover, $\Fstarep$ is defined  as
$$\Fstarep\xiet=\xiet+\cdots, \qquad \Fstaro\xiet=\xiet+\cdots,$$
so that we can chose the upper bound of $\Fstarep-\Fstaro$ without monomials of degree 0 and 1 ; $\ie$ there exists a convergent power series $\M$ such that
$$(\Fstarep-\Fstaro)\xiet\prec\nuu(\xil\+\etl\+\xiz\+\etz)^2\M(\xil\+\etl\+\xiz\+\etz).$$
\cqfd
\def\H{H}
\subsection{Back in $\R^4$, proof of Proposition \ref{propFFcomplete}}\label{RetourR}

\begin{lem}\label{LemEtapeIV}
\begin{enumerate}
\item Define $\FFep:=\mathcal{P}^{-1}\Fstarep\mathcal{P}$, where $\Fstarep$ was introduced in Lemma \ref{LemFstar}. Then
\begin{equation}\label{FFepreal}
		\overline{\FFep(\overline{\xil},\overline{\etl},\overline{\xiz},\overline{\etz})}=\FFep\xiet ;
		\end{equation}
				$\ie$ $\FFep$ is a real power series with real.
\item There exists $\HLep$ such that
$$\H\left(\FFep\xiet,\ep\right)=\HLep(\xil\etl,\xiz^2+\etz^2).$$				
\end{enumerate}
\end{lem}

\textbf{Proof of $(i)$.} In this proof we use the following notation if the power series $f$ reads $f(x)=\sum\limits_{n\in\NN^d} a_n x^n$, then we denote
$$\overline{f}(x):=\sum\limits_{n\in\NN^d} \overline{a_n} x^n.$$
Let us define
$$\begin{matrix} J_1:& \CC^4&\rightarrow&\CC^4\\
																&\xiet&\mapsto&(\overline{\xil},\overline{\etl},\overline{\xiz},\overline{\etz}) \end{matrix},$$
and denote $J_2:=\mathcal{P}J_1\mathcal{P}^{-1}$. Then  $J_2(\xil',\etl',\xiz',\etz')=(\overline{\xil'},\overline{\etl'},\I \ \overline{\etz'},\I \ \overline{\xiz'}).$
And we get that the equality \eqref{FFepreal} is satisfied if and only if
\begin{equation}
J_2\Fstarep J_2=\Fstarep. \label{EqJ2}
\end{equation}
To prove that (\ref{EqJ2}) holds we use the uniqueness of the $\Fstarep$ of $(iii)$ in Lemma \ref{LemFstar} : let us prove that $\check{\Fep}:=J_2 \Fstarep \ J_2$ verifies the conditions of $(iii)$ of Lemma \ref{LemFstar}. 

\vspace{1ex}

\textbf{Condition 1 : $\check{\Fep}$ is symplectic.}
From part \ref{SubPassageComplexe} we know that $\mathcal{P}$ is symplectic, so it is equivalent to prove that $\mathcal{P}^{-1}\check{\Fep} \mathcal{P}$ is symplectic, what holds if and only if
$$\overline{^tD\FFep \Om D\FFep} = \Om$$
holds. And this holds given that $\FFep$ is symplectic and $\Om$ is real.

\vspace{1ex}

\textbf{Condition 2 : $\Hprime(\check{\Fep}\xiet,\ep)$ is a function of $\xil\etl$, $\xiz\etz$ and $\ep$.}
Recall that lemmas \ref{LemMoser}, \ref{LemMajFep}, \ref{LemFstar}, \ref{LemFstaro} and \ref{LemMajFstar} were proved for a general Hamiltonian $\mathbf{H}$ introduced in (\ref{DefHamConstr}), and observe that the Hamiltonian 
$$\Hprime((\ql',\pl',\qz',\pz'),\ep):=\H(\mathcal{P}^{-1}(\ql',\pl',\qz',\pz'),\ep)$$
is of the form (\ref{DefHamConstr}). Moreover, the Hamiltonian $\H$ is a real power series and from Lemma \ref{LemFstar}, we know that 
$$\Hprime(\Fstarep\xiet,\ep)=\HLep(\xil\etl,\xiz\etz).$$
Then a short computation leads to 
$$\Hprime(\check{\Fep}\xiet,\ep)=\overline{\HLep}(\xil\etl,-\xiz\etz).$$

\vspace{1ex}

\textbf{Condition 3 : $\check{\Fep}$ satisfies the criteria (\ref{Q}).}
Fix $\ep$, and denote
\begin{eqnarray}
\check{\Fep}\xiet:=(\check{\varphi_1},\check{\psi_1},\check{\varphi_2},\check{\psi_2})\xiet ; \nonumber\\
\Fstarep\xiet := (\varphi_1^*,\psi_1^*,\varphi_2^*,\psi_2^*)\xiet. \nonumber
\end{eqnarray}

Then we get
$$\check{\varphi_1}\xiet=\overline{\varphi_1^*}(\xil,\etl,-\I\etz,-\I\xiz), \quad \check{\varphi_2}\xiet=\I \ \overline{\psi_2^*}(\xil,\etl,-\I\etz,-\I\xiz),$$
$$\check{\psi_1}\xiet=\overline{\psi_1^*}(\xil,\etl,-\I\etz,-\I\xiz), \quad \check{\psi_2}\xiet=\I \ \overline{\varphi_2^*}(\xil,\etl,-\I\etz,-\I\xiz).$$



Let us express $\crochet{\Frac{\check{\varphi_1}}{\xil}}$ in terms of $\crochet{\Frac{\varphi_1^*}{\xil}}$. Let us denote
$$\varphi_1^*:=\sum\limits_{m,n\in\NN^2}\alpha_{m,n}\xil^{m_1}\etl^{n_1}\xiz^{m_2}\etz^{n_2},$$
then we obtain successively
\begin{eqnarray}
&&\crochet{\frac{\varphi_1^*}{\xil}}=\sum\limits_{\dindice{m_2=n_2}{m_1=n_1+1}}\alpha_{m,n}(\xil\etl)^{n_1}(\xiz\etz)^{n_2},\nonumber\\
&&\overline{\varphi_1^*}(\xil,\etl,-\I\etz,-\I\xiz)=\sum\limits_{m,n\in\NN^2}\overline{\alpha_{m,n}}(-\I)^{m_2+n_2}\xil^{m_1}\etl^{n_1}\etz^{m_2}\xiz^{n_2} \nonumber\\
&&\crochet{\frac{\check{\varphi_1}}{\xil}}\xiet=\sum\limits_{\dindice{m_2=n_2}{m_1=n_1+1}}\overline{\alpha_{m,n}}(-\I)^{2n_2}(\xil\etl)^{n_1}(\xiz\etz)^{n_2}=\overline{\crochet{\frac{\varphi_1^*}{\xil}}}(\xil\etl,-\xiz\etz).\nonumber\\
\end{eqnarray}
Similarly, we get

\begin{eqnarray}
\crochet{\frac{\check{\psi_1}}{\etl}}(\xil\etl,\xiz\etz)&=&\overline{\crochet{\frac{\psi_1^*}{\etl}}}(\xil\etl,-\xiz\etz),\nonumber\\
\crochet{\frac{\check{\varphi_2}}{\xiz}}(\xil\etl,\xiz\etz)&=&-\overline{\crochet{\frac{\psi_2^*}{\etz}}}(\xil\etl,-\xiz\etz), \quad
\crochet{\frac{\check{\psi_2}}{\etz}}(\xil\etl,\xiz\etz)=-\overline{\crochet{\frac{\varphi_2^*}{\xiz}}}(\xil\etl,-\xiz\etz).\nonumber
\end{eqnarray}
Then, using the fact that $\Fstarep$ satisfies (\ref{Q}), we obtain that for all $(\om_1,\om_2)$, 
\begin{eqnarray}
&&\om_1\left(\crochet{\frac{\check{\varphi_1}}{\xil}}-\crochet{\frac{\check{\psi_1}}{\etl}}\right)(\om_1,\om_2)+\I\om_2\left(\crochet{\frac{\check{\varphi_2}}{\xiz}}-\crochet{\frac{\check{\psi_2}}{\etz}}\right)(\om_1,\om_2)\nonumber\\
&=&\overline{\om_1\left(\crochet{\frac{\varphi_1^*}{\xil}}-\crochet{\frac{\psi_1^*}{\etl}}\right)+\I\om_2\left(\crochet{\frac{\varphi_2^*}{\xiz}}-\crochet{\frac{\psi_2^*}{\etz}}\right)}(\om_1,-\om_2)=0,\nonumber
\end{eqnarray}
$i.e.$ $\check{\Fep}$ satisfies the criteria (\ref{Q}).
\cqfd

\vspace{2ex}

\textbf{Proof of $(ii)$.} From Lemma \ref{LemFstar}, there exists $\HLep*$ such that $\Fstarep$ verifies
$$\Hprime(\Fstarep\xiet,\ep)=\HLep^*(\xil\etl,\xiz\etz).$$
Then we obtain
\begin{eqnarray}
&&\H(\FFep\xiet,\ep)=\H(\mathcal{P}^{-1}\Fstarep\mathcal{P}\xiet,\ep)\nonumber\\
&=& \Hprime(\Fstarep(\xil,\etl,\frac{1}{\sqrt{2}}(\xiz+\I\etz),\frac{1}{\sqrt{2}}(\etz+\I\xiz)),\ep)=\HLep^*(\xil\etl,\frac{\I}{2}(\xiz^2+\etz^2)).\nonumber
\end{eqnarray}
Thus $(ii)$ holds with $\HLep(\xil\etl,\xiz^2+\etz^2):=\HLep^*(\xil\etl,\frac{\I}{2}(\xiz^2+\etz^2))$.
\cqfd

\vspace{2ex}

\textbf{Proof of Proposition \ref{propFFcomplete}.} Let us prove that the family $\FFep$ defined in Lemma \ref{LemEtapeIV} satisfies the results claimed in Proposition \ref{propFFcomplete}. Let $\roop$ be such that $\roop<\roo$ and $4\roop$ is a radius of convergence of the power series $\M$ of Lemma \ref{LemMajFstar} and such that $\roop$ is a radius of convergence of $\Fstaro$ (recall that Lemma \ref{LemMajFstar} ensures that $\Fstaro$ is convergent). And define $\Mo:=\NNorme{\An(\B(0,4\roop))}{\M}$.

Then (\ref{DefHL}) holds thank to $(ii)$ of Lemma \ref{LemEtapeIV} and (\ref{EqFstaro}) is a consequence of Lemma \ref{LemFstaro}.

$(i)$ is a consequence of Lemma \ref{LemMajFstar} ; $(ii)$ and $(iii)$ hold thank to Lemma \ref{LemMajFstar} and because the monomial of degree 1 of $\Fstaro\xiet$ is $\xiet$ and the monomial of degree 1 of $(\Fstaro)^{-1}\qp$ is $\qp$.

We get $(iv)$ from $(i)$ above and $(vi)$ of Lemma \ref{ProprPrec} ; we obtain $(v)$ from $(i)$ above together with Lemma \ref{LemPrecInv} and $(vi)$ of Lemma \ref{ProprPrec} ; $(vi)$ and $(vii)$ are a consequence of $(i)$ and $(ii)$ above and of Lemma \ref{ProprPrec}-$(vi)$ ; we get $(viii)$ and $(ix)$ from $(vi)$ and Lemma \ref{LemFstaro}.
\cqfd

\end{document}